\documentclass[11pt]{article}

\usepackage{amsmath, amssymb, amsthm}
\usepackage{mathtools}
\usepackage{graphicx}
\usepackage{placeins}   
\usepackage{booktabs}
\usepackage{siunitx}
\usepackage{enumitem}
\usepackage{hyperref}
\usepackage{xcolor}
\usepackage{float}
\usepackage{comment}
\usepackage{xfrac}
\usepackage{fullpage}

\newcommand{\R}{\mathbb{R}}

\newcommand{\muvec}{\boldsymbol{\mu}}

\newcommand{\normLtwo}[1]{\left\lVert #1 \right\rVert_{2}}

\newcommand{\bU}{\mathbf{U}} 
\newcommand{\bV}{\mathbf{V}}
\newcommand{\bF}{\mathbf{F}} 
\newcommand{\bX}{\mathbf{X}} 

\newcommand{\calP}{\mathcal{P}}
\newcommand{\calQ}{\mathcal{Q}}
\newcommand{\calG}{\mathcal{G}}
\newcommand{\calC}{\mathcal{C}}
\newcommand{\calS}{\mathcal{S}}
\newcommand{\calN}{\mathcal{N}}

\newcommand{\Nt}{\tilde{N}}

\newcommand{\oh}{\sfrac{1}{2}}

\newcommand{\rev}[1]{%
  \begingroup
  #1%
  \endgroup
}

\newenvironment{revblock}
{\begingroup 
}
{\endgroup}

\title{Comparison of a Parametric Physics-Informed Neural Network and a Tensorial Reduced-Order Model for the Shallow-Water Dam-Break Problem}

\author{
Anton Myshak \thanks{Dept. of Mathematics, University of Houston, Houston, TX 77204,
amyshak@uh.edu} \quad
Md Rezwan Bin Mizan \thanks{Dept. of Mathematics, University of Houston, Houston, TX 77204,
mbinmizan@uh.edu} \quad
Ilya Timofeyev\thanks{Dept. of Mathematics, University of Houston, Houston, TX 77204, itimofey@cougarnet.uh.edu}
}
\date{} 

\begin{document}
\maketitle

\begin{abstract}
We develop two parametric data-driven reduced models: a physics-informed neural network (PINN) and a non-intrusive tensorial reduced-order model (TROM), and apply both approaches to the parametrized one-dimensional shallow-water dam-break problem.
Neither reduced model requires time integration: both learn a direct parameter-to-solution map from space, time, and dam-break parameters to the physical state, with the PINN providing predictions at arbitrary times and the TROM reconstructing solutions at the stored snapshot times.
In addition, we demonstrate that it is essential to introduce shock-aware collocation to improve the robustness of the PINN model. 
\end{abstract}

\smallskip


\smallskip


\section{Introduction}
\label{sec:intro}
Hyperbolic balance laws are ubiquitous in many areas of applied mathematics, including fluid flows. One of the main characteristics of hyperbolic problems is the appearance of sharp transients and discontinuities, which are intrinsic to the dynamics.
One of the 
canonical benchmarks is the one-dimensional (1D) shallow-water dam-break problem (e.g., \cite{stoker2019water,toro2024swe,castro2017ritter}).
Despite its conceptual simplicity, the solution typically develops nonlinear wave interactions, moving fronts, and (depending on wet/dry configurations) sharp transitions (shocks) that are challenging numerically for both classical discretizations and surrogate modeling. 
Reliable and efficient prediction of such regimes is central to many practical problems, including uncertainty quantification, design and control, and real-time forecasting, motivating the development of reduced models that remain accurate in the presence of transport and shocks \cite{LeVeque2002FVM,Toro2009Riemann,ohlberger2015reduced,fluids11030076,abgrall2019,chan2020}.

Projection-based reduced-order models (ROMs) (see survey \cite{benner2015survey} and references therein) constitute a well-established approach to accelerating numerical simulations by projecting the state onto a low-dimensional approximation space constructed from high-fidelity snapshots.
However, it is known that for transport-dominated problems that contain shocks, the solution manifold often exhibits slow decay of the Kolmogorov $n$-width \cite{melenk2000n,lassila2013generalized} and its approximations \cite{cohen2010convergence,buffa2012priori,cohen2016kolmogorov,hesthaven2016certified,bachmayr2017kolmogorov} computed by greedy POD \cite{binev2011convergence}, limiting the efficiency of linear-subspace compression and leading to a loss of fidelity near discontinuities.
Consequently, substantial research has focused on ROM variants that mitigate these issues, for instance, through localization, nonlinear manifolds, or feature alignment (e.g., \cite{barnett2022quadratic,amsallem2012nonlinear,reiss2018shifted,taddei2020registration,Rim2021TransportedSubs}).
Recently, a parametric tensor Reduced-Order Model (TROM) was proposed to address some of these issues \cite{mamonov2022interpolatory,mamonov2024tensorial,olshanskii2025approximating}. This approach combines a low-rank tensor representation of the parameter-to-solution map with localized reduced bases to improve accuracy across regimes. 
The TROM is particularly suitable for the development of reduced models for parametric problems.
In \cite{Mizan2025TROM} and 
\cite{Mizan2026TROM}, this approach was successfully applied to the dam-break problem for the shallow-water equations for Newtonian and visco-plastic fluids, respectively. In this paper, we consider the non-intrusive TROM, where solutions for new parameter values are computed directly from the low-rank compressed representation without time integration of the reduced equations. Details are provided in Section \ref{sec:trom}. 

In parallel, physics-informed neural networks (PINNs) \cite{Raissi2019PINN} have emerged as a promising paradigm for constructing differentiable surrogate models by embedding governing equations directly into the learning objective.
In their classical form, PINNs optimize a neural network representing the solution fields using a composite loss that penalizes violations of the partial differential equations (PDEs) as well as mismatches to initial and boundary data \cite{Raissi2019PINN}.
This hybridization of data and physics has shown a strong potential for learning PDE solution operators from relatively sparse supervision, offering a flexible route to surrogate modeling that is mesh-free in space and time.
Nevertheless, hyperbolic conservation laws pose additional difficulties: solutions may not be smooth, and, thus, the differential form of the PDE is not valid at discontinuities, which can degrade training stability and shock localization. 
This has been addressed by introducing various modifications such as 
conservative PINNs  (cPINNs) \cite{Jagtap2020cPINN}, 
weak PINNs (wPINNs) and 
Integral PINNs (IPINNs)
\cite{deRyck2024wpinn, chaumet2024improvingwpinn, RajvanshiKetcheson2024IPINN} that utilize weak formulation of PDEs, variational PINNs (VPINNs)
\cite{kharazmi2019vpinn,kharazmi2021hpvpinn} that adopt a Petrov--Galerkin viewpoint by pairing neural-network trial functions with polynomial test spaces.
Architectural choices have also been explored; for example, physics-informed attention-based models have demonstrated improved shock-front representation and generalization behavior in hyperbolic settings \cite{RodriguezTorrado2022PIANN}.

Recently,
PINNs have been applied to free-surface flows governed by the SWEs, motivated by the need for fast surrogates in forecasting and scenario exploration. Recent studies developed PINNs for SWEs in one and two spatial dimensions, including dam-break scenarios \cite{pinnswe-bihlo2022,pinnswe-qi2024,mumtaz2025dambreakPINN}. These studies also discuss practical issues such as stabilization, loss weighting, and the trade-off between smoothness and accuracy relative to shock-capturing numerical schemes. However, these studies primarily address fixed PDE instances or individual forward scenarios rather than learning a single parametric surrogate over a family of dam-break initial conditions. 
Thus, while they support the feasibility of physics-informed learning for SWE dynamics, they do not resolve the parametric generalization problem studied here.

This work constructs and compares two parametric reduced-order models for the one-dimensional shallow-water dam-break problem: a physics-informed neural network (PINN) and a non-intrusive tensorial reduced-order model (TROM). Both models are constructed using the same Full-Order Model (FOM) dataset and evaluated using consistent metric across dry-bed, wet-bed, and near-dry-bed regimes, following benchmark settings used in previous ROM studies \cite{chang2011numerical,magdalena2022numerical,muchiri2024numerical,peng20121d,podromswe2012,zoppou2000numerical} and recent TROM work \cite{Mizan2025TROM}. We compare their predictive accuracy for the water depth and the discharge, their performance in parameter interpolation and extrapolation, and their respective offline and online computational costs. \rev{It is important to emphasize that both surrogate models are evaluated against the same FOM reference solution, rather than against the exact solution of the SWE, so that the surrogate approximation error is not conflated with the discretization error of the underlying numerical FOM solver.} An important component of the study is the construction of a robust parametric PINN for this shock-dominated problem. Through ablation experiments, we demonstrate that shock-aware collocation, residual scaling, and PDE gating substantially improve the accuracy and robustness of the PINN. Overall, the comparison identifies the complementary strengths and limitations of the two approaches for parametric hyperbolic problems.

The rest of the paper is organized as follows.
Section~\ref{sec:swe} introduces the shallow-water model and the parametrized dam-break setup.
Sections \ref{sec:trom} and \ref{sec:pinn}
give an overview of the tensorial Reduced-Order Model (TROM) and parametric Physics-Informed Neural Networks (PINNs).
The numerical results 
for the interpolation and extrapolation in the parameter space
are presented in 
Sections \ref{sec:num} and \ref{sec:extrap}, respectively.
Ablation studies for the PINN model are presented in Section \ref{sec:ablation}.
The discussion and conclusions are presented in Section 
\ref{sec:conc}.

\section{Methods and Setup}
\label{sec:swe}

In this section, we briefly outline the dam-break problem for the shallow-water equations and the numerical method.

\subsection{Governing equations}

We consider the one-dimensional shallow water equations (SWEs), posed on a spatial interval $\Omega := [0,L]$ and time interval $[0,T]$. 
The equations can be written as a system of conservation laws 
\begin{equation}
\label{eq:swe}
  \partial_t \bU + \partial_x \bF(\bU) = \mathbf{0},
  \qquad (x,t)\in \Omega\times(0,T],
\end{equation}
with 
\[
\bU(x,t;\muvec) :=
\begin{bmatrix}
    h(x,t;\muvec) \\
    q(x,t;\muvec)
\end{bmatrix}, \qquad
\bF(\bU) = 
\begin{bmatrix}
q \\
\dfrac{q^2}{h} + \dfrac12 g h^2 
\end{bmatrix},
\]   
where $h(x,t;\muvec)\ge 0$ denotes the water depth, $q(x,t;\muvec)=h(x,t;\muvec)\,u(x,t;\muvec)$
is the discharge (momentum), $u$ is the depth-averaged velocity,  and $g>0$ is the gravitational acceleration.
Here $\muvec$ is a vector of parameters describing the initial conditions and $\bF(\bU)$ is the flux function. 
The system in \eqref{eq:swe} is a strictly hyperbolic system for $h>0$; solutions may develop discontinuities (shocks) even from smooth initial data. Consequently, the relevant notion of solution is that of a (weak) entropy solution; see, e.g., \cite{toro2024swe} for a detailed discussion.

\paragraph{Parametrized dam-break initial data.}
We study a parametrized dam-break configuration with a depth discontinuity located at
$x_{\mathrm{dam}} = L/2.$
The parameter vector is $\muvec := (h_L,h_R)\in\mathcal{P}\subset\R^2$, where $h_L$ and $h_R$
denote the upstream and downstream depths, respectively, and satisfy
$h_L > h_R \ge 0.$
Thus, the initial dam condition is
\begin{equation}\label{eq:ic_dambreak}
  h(x,0;\muvec) =
  \begin{cases}
    h_L, & x < x_{\mathrm{dam}},\\
    h_R, & x \ge x_{\mathrm{dam}},
  \end{cases}
  \qquad
  u(x,0;\muvec) = 0,
\end{equation}
which implies $q(x,0;\muvec)=0$.

The case $h_R=0$ is commonly referred to as the \emph{dry-bed} dam-break problem and yields a solution dominated by a rarefaction wave, while the case $h_R>0$ is the \emph{wet-bed} dam-break problem and is more challenging due to the presence of a downstream shock for all $h_R>0$ and for all $t>0$. Although the solution of the wet-bed case converges uniformly to the solution of the dry-bed case on finite time intervals $[0,T]$ such that $x_{\mathrm{dam}} + 2cT < L$ with $c = \sqrt{gh_L}$, $\partial h / \partial h_R \sim O(h_R^{-1/2})$
in a small interval near the shock for small $h_R$
(see discussion in \cite{Mizan2025TROM} Appendix A).
This complicates the use of projection-based reduction techniques whose performance depends on the smoothness of the parametric solution manifold. 

We consider a standard benchmark setup with 
$L=100$, boundary conditions
\begin{equation}\label{eq:outflow_bc}
  \partial_x h(0,t;\muvec)=\partial_x u(0,t;\muvec)=
  \partial_x h(L,t;\muvec)=\partial_x u(L,t;\muvec)=0,
  \quad t\in(0,T],
\end{equation}
and parameter domain
\begin{equation}
\label{eq:param_domain_box}
  [h_L^{\min},h_L^{\max}] \times [h_R^{\min},h_R^{\max}]
  = [10,28]\times[0,8].
\end{equation}

\paragraph{Parametric solution operator and learning objective.}
Let $\bU^{\mathrm{FOM}}(x,t;\muvec)$ denote the \rev{FOM numerical reference} solution of
\eqref{eq:swe}.
The main goal of this paper is to examine two types of surrogate data-driven models  
for the \emph{parameter-to-solution map}
\begin{equation*}
  \calS :\ (x,t,\muvec)\ \mapsto\ \bU^{\mathrm{FOM}}(x,t;\muvec),
  \qquad (x,t)\in\Omega\times[0,T],\ \muvec\in\mathcal{P}.
\end{equation*}
\rev{We emphasize that the parameter-to-solution map $\calS$ is defined by the FOM solver and therefore inherits the discretization errors associated with the underlying numerical method.}
The first model is a parametrized 
physics-informed neural network
$\calN_\theta := (h_\theta(x,t;\muvec),u_\theta(x,t;\muvec))$ where 
$\theta$ are trainable weights.
This model should be consistent with the induced conservative quantity $q_\theta := h_\theta u_\theta$.
The surrogate should provide accurate predictions for previously unseen parameter values $\muvec$ and
times $t$ (with particular interest in generalization within the training parameter domain). The
learning problem is therefore to approximate $\mathcal{S}$ on $\Omega\times[0,T]\times\mathcal{P}$
while respecting the governing PDE \eqref{eq:swe}, initial condition
\eqref{eq:ic_dambreak}, boundary condition \eqref{eq:outflow_bc}, and admissibility constraint
$h(x,t) \ge 0$.

The second model is a non-intrusive tensorial reduced-order model (TROM). Unlike neural networks, this model does not require any iterative training, but it is still a data-driven model since it requires a dataset generated by numerical simulations of \rev{the FOM solver.}
The TROM uses a reduced representation which is obtained by performing a higher-order singular value decomposition (HOSVD)
of the full snapshot matrix, nonlinear interpolation in parameter space, and reconstruction of the solution matrix from the reduced tensor representation. Essentially, this model learns the 
parameter-to-solution map $\calS$ by computing the Tucker reduced tensor representation during the offline phase.

\subsection{Numerical Method for the Full-Order Model}  
We use a Local Lax-Friedrichs finite-volume discretization to discretize the SWE in space. 
The semi-discrete scheme becomes
\[
\frac{d}{dt} \bU_i(t)  = -\frac{\bF_{i+ \oh} - \bF_{i-\oh}}{\Delta x}, \qquad i=1,\ldots,N
\]
with
\[ 
\bF_{i+\oh}
  =
  \frac{\bF(\bU_i)+\bF(\bU_{i+1})}{2}
  -\frac{\lambda_{i+\oh}}{2}
  \left(\bU_{i+1}-\bU_i\right),
\]
with $\lambda_{i+\oh} = \max\left(|u_i|+\sqrt{g h_i},\; |u_{i+1}|+\sqrt{g h_{i+1}}\right)$
and $\Delta x = L/N$, where $N$ is the number of discrete points within the domain.
For time-stepping, we use the Euler discretization with the time-step $\Delta t$ satisfying the CFL condition 
$\max_{i}\lambda_{i+\oh} \, \Delta t / \Delta x < 1.$
We use two ghost cells $\bU_0 = \bU_1$ and
$\bU_{N+1}=\bU_N$ to satisfy the boundary conditions.
It is known that the above discretization belongs to the class of 
strong stability preserving 
methods \cite{shu1988efficient,gottlieb2001strong,gottlieb2011strong}.
The scheme is also height-positivity-preserving.
Thus, the non-negativity of the water depth $h_i$ is preserved at the discrete level for all times.

\rev{To assess the influence of the FOM discretization, we also implemented a Lax-Wendroff finite-volume scheme equipped with the van Leer limiter (LW-vL) and compared both numerical solutions with the exact Riemann solution. As expected, the LW-vL scheme is more accurate, with errors in the water depth \(h\) and discharge \(q\) typically smaller by a factor of approximately 2-3 for the representative regimes considered. Nevertheless, the errors obtained with the two schemes remain of the same order of magnitude, and both LLF and LW-vL reproduce the same overall dam-break solution structure, including the rarefaction and shock propagation. The main differences are concentrated near discontinuities, where the LLF solution is more diffusive. In the present study, the exact solution is used only to assess the discretization error of the FOM, whereas the accuracy of the PINN and TROM is evaluated relative to the common LLF-generated FOM reference.}

\subsection{Dataset construction}
\label{sec:dataset}
\rev{In the present study, we use a reduced parameter grid consisting of \((N_L\times N_R)=(5\times13)\) parameter instances, with sampled parameter values}
\begin{equation}\label{eq:paramgrid}
\begin{aligned}
  \calP_1 &= \{10, \,14.5,\,19,\,23.5,\,28\},\\
  \calP_2 &= \{0,\,0.1363,\,0.5359,\,1.1716,\,2,\,2.9647,\,4,
\,5.0353,\,6,\,6.8284,\,7.4641,\,7.8637,\,8\}.
\end{aligned}
\end{equation}
Therefore, the total sampled parameter set $\calP = \calP_1 \times \calP_2$
contains 
$N_L N_R=65$ sampled parameter values.

We integrate the full-order model (FOM) with the time-step $\Delta t = 10^{-4}$
to a final time $T=2.5$ and store snapshots at times $t_k=k\Delta t_s$, $k=1,\dots,N_T$ with $\Delta t_s=0.1$. We use $N_x=402$, which includes two
ghost points (one on each boundary) for implementing the boundary conditions.
\rev{Note that the snapshot times $t_k$ do not include $t=0$. The initial state is given analytically by \eqref{eq:ic_dambreak} for every parameter value, so no surrogate is required at $t=0$, and all comparisons reported below are made at $t>0$. Should reconstruction at $t=0$ be desired, the analytically known initial state can be appended to the snapshot data as an additional temporal slice without any further full-order computation.}
For each parameter value $\muvec^{(m)}$, $m=1,\ldots,N_L N_R$, 
we store fields $h^{(m)}(x_i,t_k)$ and $u^{(m)}(x_i,t_k)$ with $i=1,\ldots,N_x$ and $k=1,\ldots,N_T$.
These snapshots can be combined into a fourth-order tensor, as discussed in Section \ref{sec:trom}.

\begin{figure}[H]
  \centering
  \includegraphics[width=0.95\linewidth]{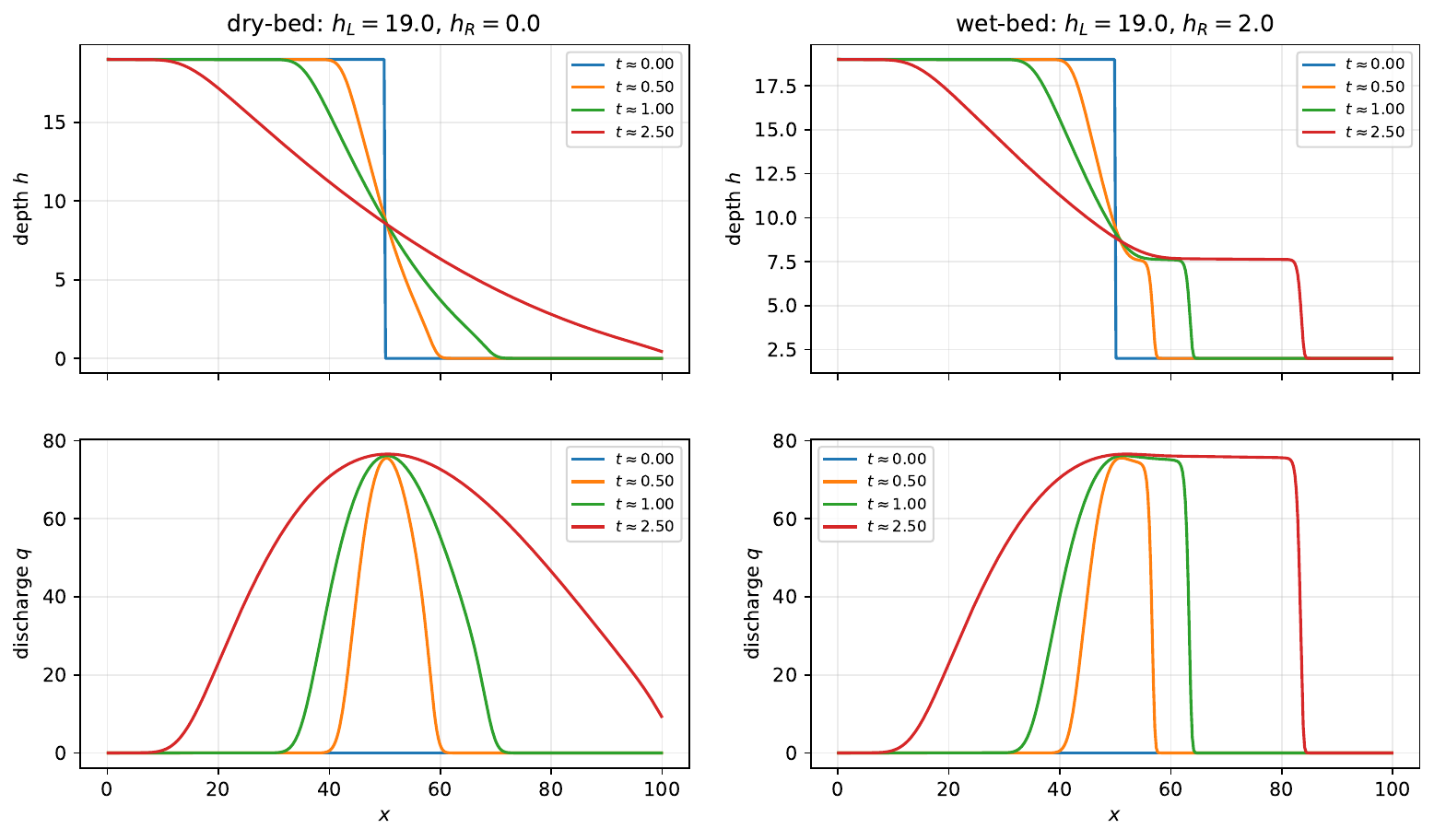}
  \caption{Representative FOM snapshots (depth $h(t,x)$ and discharge $q(t,x)$) for a dry-bed and a wet-bed choice of parameters $(h_L, h_R)=(19,0)$ (left) and $(h_L, h_R)=(19,2)$ (right), at times $t=0, 0.5, 1.0, 2.5$. Top part - water depth $h(t,x)$, Bottom part - discharge $q(t,x)$.}
  \label{fig:example_snapshots}
\end{figure}

Figure~\ref{fig:example_snapshots} illustrates representative space--time structure in the generated FOM trajectories for a dry-bed case $(h_L,h_R)=(19,0)$ and a wet-bed case $(h_L,h_R)=(19,2)$.
The depth profiles $h(t,x)$ (top row) for $t=0, 0.5, 1, 2.5$ highlight the qualitative difference between regimes: the dry-bed solution is dominated by a rarefaction wave that spreads from the dam location, whereas the wet-bed solution exhibits a persistent downstream shock for all $t>0$. The corresponding discharge profiles $q(t,x)=h(t,x)u(t,x)$ (bottom row) depict the local flow rate (volume flux per unit width) and highlight where the solution carries significant momentum. 
\FloatBarrier

\paragraph{Train/validation/test splits.}
\rev{We design the evaluation protocol to probe generalization in parameter space, covering both interpolation within the sampled parameter domain (Section~\ref{sec:num}) and extrapolation beyond it (Section~\ref{sec:extrap}); in either case, predictions are compared with the FOM at the stored snapshot times $t_k$. The split used for training the PINN is specified in Section~\ref{sec:pinn}. The TROM requires no such split, as discussed at the beginning of Section~\ref{sec:trom}.}


\section{Tensorial reduced-order Model}
\label{sec:trom}
In this section, we present the non-intrusive Tensorial Reduced-Order Model (TROM); additional details can be found in \cite{rezwan2026thesis}.
As we discussed in Section \ref{sec:dataset}, the solution of the full-order model depends on parameters 
$\muvec^{(m)} \in \R^2$ and we consider 
the parameter grid $\calP$ defined in \eqref{eq:paramgrid}.
Using the full-order model, 
we generate solution snapshots $\{h^{(m)}(x_i, t_k), \, u^{(m)}(x_i, t_k)\}$, where $\muvec^{(m)} \in \R^2$ is a particular parameter value from \eqref{eq:paramgrid}. 
For each field variable $\phi \in \{h, u\}$, these snapshots can be 
organized into a fourth-order tensor 
\begin{equation*}
\calQ_\phi \in \R^{N_x \times N_L \times N_R \times N_T},
\end{equation*}
where the first mode corresponds to the spatial discretization ($N_x = 402$), modes 2 and 3 correspond to the parameter dimensions ($N_L=5$ and $N_R=13$), and the last mode corresponds to time ($N_T = 25$). 
\rev{We emphasize that, unlike the parametric PINN of Section~\ref{sec:pinn}, the TROM requires no iterative training.
The low-rank representation is computed once, from the snapshots at all $N_LN_R=65$ parameter values of the grid $\calP$ in \eqref{eq:paramgrid}, all $N_T=25$ snapshot times, and all $N_x=402$ spatial grid points, without subsampling in space, time, or parameter space. The threshold $\epsilon$ in \eqref{eq:rank_selection} below is the only tunable parameter of the construction. All parameter values at which the TROM is evaluated in Sections~\ref{sec:num} and~\ref{sec:extrap} lie outside the grid \eqref{eq:paramgrid}, so every TROM result reported in this paper is out of sample with respect to the parameter grid.}
Next, the Tucker decomposition \cite{tucker1966some} can be used to develop a compressed tensor format.
In particular, the tensor 
$\calQ_\phi$ can be represented as
\begin{equation}
\label{eq:tucker}
\calQ_\phi \approx 
\widetilde{\calQ_\phi} = 
\calG_\phi \times_1 \bV^{(x)}_\phi \times_2 \bV^{(L)}_\phi 
\times_3 \bV^{(R)}_\phi
\times_4 \bV^{(T)}_\phi,
\end{equation}
where $\calG_\phi \in \R^{\Nt_x \times \Nt_L \times \Nt_R \times \Nt_T}$ is the \textit{core tensor} and $\bV^{(k)}_\phi \in \R^{N_k \times \Nt_k}$ are factor matrices with $k \in \{x,L,R,T\}$.
Here $" \cdot \times_n \cdot "$ denotes the \textit{mode-$n$ product} of a tensor
with a matrix (see, e.g., \cite{kolda2009tensor} for a review of tensor algebra).
The tuple $(\Nt_x, \Nt_L, \Nt_R, \Nt_T)$ represents the multilinear rank of the tensor $\widetilde{\calQ_\phi}$. 
The Tucker format allows independent compression ratios for each mode, and orthonormality of the factor matrices ensures numerical stability and direct error control. 
The Higher-Order SVD (HOSVD) \cite{de2000multilinear} can be used to compute the Tucker decomposition by performing the SVD of each mode-$n$ unfolding independently and selecting the truncation rank for each mode. 
The multilinear ranks $\Nt_k$ with $k\in \{x,L,R,T\}$ are selected using a criterion based on the singular value decay:
\begin{equation}
\Nt_k = \min \left\{ r : \frac{\sum_{i=1}^r \sigma_{i,k}}{\sum_{i=1}^{N_k} \sigma_{i,k}} \geq 1 - \epsilon \right\},
\label{eq:rank_selection}
\end{equation}
where $\sigma_{i,k}$ are the singular values of the $k$-mode unfolding of $\calQ_\phi$. 
We use $\epsilon = 10^{-3}$ as a threshold for all modes.
\begin{revblock}
The criterion \eqref{eq:rank_selection} is based on the normalized cumulative sum of the singular values rather than of their squares.
The two variants are related: for a non-increasing sequence $\sigma_1 \ge \sigma_2 \ge \dots \ge \sigma_N \ge 0$ and any $r$,
\[
\Big(\sum_{i\le r}\sigma_i^{2}\Big)\Big(\sum_{i=1}^{N}\sigma_i\Big)
-\Big(\sum_{i\le r}\sigma_i\Big)\Big(\sum_{i=1}^{N}\sigma_i^{2}\Big)
=\sum_{i\le r}\sum_{j>r}\sigma_i\sigma_j\big(\sigma_i-\sigma_j\big)\ \ge\ 0 ,
\]
since $\sigma_i \ge \sigma_j$ whenever $i \le r < j$, so the normalized partial sums of the squared singular values dominate those of the singular values themselves. Consequently, for a fixed $\epsilon$, the criterion \eqref{eq:rank_selection} selects multilinear ranks no smaller than those produced by the standard squared--singular-value (energy) criterion, i.e., it is the more conservative of the two.
\end{revblock}
The multilinear ranks obtained from the truncated HOSVD using \rev{the criterion above} are
\(
(\Nt_x, \Nt_L, \Nt_R, \Nt_T)
= (\text{94}, \text{4}, \text{5}, \text{21}).
\)

\rev{We next quantify the storage reduction achieved by the Tucker representation at these multilinear ranks. For each field, the snapshot tensor $\calQ_\phi$ contains $N_xN_LN_RN_T = 402\times5\times13\times25 = 653\,250$ entries, whereas the Tucker representation \eqref{eq:tucker} requires $\Nt_x\Nt_L\Nt_R\Nt_T = 39\,480$ entries for the core tensor and $402\times94 + 5\times4 + 13\times5 + 25\times21 = 38\,398$ entries for the four factor matrices, i.e., $77\,878$ entries in total. This corresponds to a storage reduction by a factor of approximately $8.4$, or $11.9\%$ of the uncompressed snapshot data. Note that the ranks in the spatial and temporal modes are limited by the size of those modes themselves: with only $N_T=25$ stored snapshots the temporal mode is intrinsically low-dimensional and close to full rank by construction, so for a parametric surrogate the essential compression takes place in the parameter modes and in the core tensor. We also emphasize that the non-intrusive TROM is not intended to circumvent the slow decay of the Kolmogorov $n$-width that is characteristic of transport-dominated problems and was discussed in Section~\ref{sec:intro}; the ranks above are those returned by the conservative criterion \eqref{eq:rank_selection} at $\epsilon = 10^{-3}$ and are not claimed to be minimal. The accuracy and the online cost of the resulting representation are reported in Sections~\ref{sec:num} and~\ref{sec:extrap}, and its limitations near the moving front and outside the sampled parameter range are discussed in Section~\ref{sec:conc}.}
For the non-intrusive TROM,
solution snapshots for a new parameter value $\muvec^* \in \R^2$ are constructed directly from the low-order representation \eqref{eq:tucker} without any time integration
(see also \cite{Mizan2026TROM}).
First, interpolation coefficients for the new parameter values are computed using cubic spline interpolation on the grid $\calP_k$
\begin{equation*}
	\chi_k(\mu_k^*) = \text{interp. coefficients}(\mathcal{P}_k,  \mu_k^*), \quad k = 1,2.
\end{equation*}
The vector $\chi_k \in \R^{N_k}$ contains the interpolation coefficients corresponding to the parameter value $\mu_k^*$ on the grid $\mathcal{P}_k$. 
Next, the Tucker tensor from \eqref{eq:tucker} is contracted along all parameter modes with the interpolation vectors:
\begin{equation}
	\calC_\phi(\muvec^*) = \mathcal{G}_\phi \times_2 \left(\bV_\phi^{(L)\top}\chi_1(\mu_1^*)\right) \times_3 
     \left(\bV_\phi^{(R)\top}\chi_2(\mu_2^*)\right).
\end{equation}
This yields the low-dimensional matrix $\calC_\phi(\muvec^*) \in \R^{\Nt_x \times \Nt_T}$.
The  solution  is then reconstructed via
\begin{equation*}
\bX_\phi(\muvec^*) =\bV^{(x)}_\phi \, \calC_\phi(\muvec^*) \, \left(\bV_\phi^{(T)}\right)^{\!\top}  \in \R^{N_x \times N_T},
\end{equation*}
where the $k$-th column of $\bX_\phi(\muvec^*)$ approximates $\phi(\mathbf{x}, t_k; \muvec^*)$.
\rev{The reconstruction above returns the solution at the stored snapshot times $t_k$, $k=1,\dots,N_T$, which are the columns of $\bX_\phi(\muvec^*)$; the non-intrusive TROM as formulated here performs interpolation in parameter space only. Predictions at times between the stored snapshots can be obtained within the same framework by contracting the temporal mode of \eqref{eq:tucker} with spline interpolation coefficients $\chi_T(t^*)$ computed on the grid $\{t_k\}$, exactly as is done above for the parameter modes. We do not pursue this extension here, and all TROM results reported below are evaluated at the snapshot times $t_k$.}

\section{Parametric Physics-Informed Neural Network}
\label{sec:pinn}

Recall from Section~\ref{sec:swe} that the parametric solution operator maps
\(
(x,t,\muvec)\in \Omega\times[0,T]\times\mathcal{P}
\)
to the state of the flow.
Here, $\Omega = [0,L]$ and $\mathcal{P} \subset \R^2$ is the parameter domain.
Accordingly, we seek a neural network approximation 
$\calN_\theta$
of the form
\begin{equation*}
  \calN_\theta:\R^4\to\R^2,
  \qquad
  \calN_\theta: (x,t,\muvec)\ \mapsto\ (h_\theta(x,t;\muvec),\,u_\theta(x,t;\muvec)),
\end{equation*}
where \(\theta\) denotes the trainable network parameters. Additional details about the PINN
presented in this Section can be found in \cite{myshak2026thesis}.

Rather than feeding
\(\muvec=(h_L,h_R)\) directly into the network, we use the equivalent parameterization
\begin{equation}
  \bar h := \frac{h_L+h_R}{2},
  \qquad
  \gamma := \frac{h_L-h_R}{h_L+h_R},
  \qquad
  \bar \muvec := (\bar h, \gamma)
  \label{eq:param_reparam}
\end{equation}
Thus, the actual network input and output are
\begin{equation*}
  z := (x,t,\bar \muvec) = (x,t,\bar h,\gamma)\in\mathbb{R}^4,
  \qquad
  \mathcal{N}_\theta(z)=
  \begin{bmatrix}
    h_\theta(z)\\
    u_\theta(z)
  \end{bmatrix}\in\mathbb{R}^2.
\end{equation*}
This reparameterization separates the mean water-depth scale \(\bar h\) from the relative jump magnitude \(\gamma\). Empirically, this provides a more balanced description of the regime space and is natural for dam-break problems.

\paragraph{Network architecture}
We approximate the parametric solution map by a fully connected feed-forward neural network without skip connections or residual blocks. 
The network with $K$ hidden layers is defined recursively by
\begin{align}
  a^{(0)} &= z, \\
  a^{(k+1)} &= \sigma\!\left(W^{(k)}a^{(k)}+b^{(k)}\right),
  \qquad k=0,\dots,K-2, \\
  \mathcal{N}_\theta(z) &= W^{(K-1)}a^{(K-1)}+b^{(K-1)},
\end{align}
where \(W^{(k)}\) and \(b^{(k)}\) are the trainable weights and biases of layer \(k\), and \(\sigma\) is the hidden-layer activation function.
Throughout this work, we use the \emph{sigmoid linear unit} (SiLU) activation function
\(
\sigma(x)=x\,\mathrm{sigmoid}(x)={x}/(1+e^{-x}).
\)
The SiLU activation is smooth and non-polynomial, which is advantageous in PINNs because the PDE residual requires repeated differentiation of the network output with respect to the input coordinates. Compared with piecewise-linear activations, smooth activations produce smooth derivatives and therefore more stable automatic-differentiation-based residuals.

All linear layers are initialized with \emph{Xavier normal initialization} \cite{glorot2010understanding}. 
If a fully connected layer has weight matrix \(W\), then each entry is initialized as
\(
  W_{ij}\sim \mathcal{N}\!\left(0,\;2/(m+n)\right).
\)
This initialization is designed to keep activation magnitudes at a reasonable scale across layers and to improve optimization stability in deep feed-forward networks. Biases are initialized to zero.


\paragraph{Physics-informed residual}
Although the network predicts the primitive variables \((h_\theta,u_\theta)\), the governing equations from Section~\ref{sec:swe} are enforced through the conservative form of the PDE, using the induced discharge \(q_\theta=h_\theta u_\theta\). Then the two strong-form residuals are
\begin{align}
  \mathcal{R}_1(x,t;\muvec)
  &:= \partial_t h_\theta + \partial_x(h_\theta u_\theta),
  \label{eq:r1_method}
  \\
  \mathcal{R}_2(x,t;\muvec)
  &:= \partial_t(h_\theta u_\theta)
  + \partial_x\!\left(h_\theta u_\theta^2 + \frac{1}{2}g h_\theta^2\right).
  \label{eq:r2_method}
\end{align}
These residuals are exactly the mass and momentum residuals corresponding to the conservative shallow-water system when expressed through the primitive variables.
In practice, the derivatives in \eqref{eq:r1_method}--\eqref{eq:r2_method} are computed by automatic differentiation in PyTorch, by differentiating the network outputs with respect to the input coordinates \(x\) and \(t\) using \texttt{autograd}.

This formulation keeps the exposition consistent with the rest of the paper: the network predicts \((h,u)\), training data are stored in \((h,u)\), and only the induced conservative quantity \(q=hu\) is formed internally when the governing equations require it.

\subsection{Training data and collocation sets}
Training uses four types of information: (1) supervised solution snapshots, (2) collocation points for PDE residual enforcement, (3) boundary points for enforcing the boundary conditions, and (4) points at \(t=0\) for enforcing the initial condition.
First, let
\[
\mathcal{D}_{\mathrm{data}}
=
\{(x_i,t_i,\mu_i,h_i,u_i)\}_{i=1}^{N_d}
\]
denote supervised solution snapshots generated by the reference finite-volume solver. These points are sampled from multiple dam-break regimes in the training set. The spatial and parametric coordinates of boundary points extracted from these snapshots form the set
\[
\mathcal{D}_{\mathrm{bc}}
=
\{(x_k^{\mathrm{bc}},t_k^{\mathrm{bc}},\mu_k^{\mathrm{bc}})\}_{k=1}^{N_{\mathrm{bc}}},
\]
where \(x_k^{\mathrm{bc}}\in[0,L]\). These points are used to enforce the zero-gradient boundary condition.

In addition, we construct the initial-condition set
\[
\mathcal{D}_{\mathrm{ic}}
=
\{(x_m^{\mathrm{ic}},0,\mu_m^{\mathrm{ic}},
h_m^{\mathrm{ic}},u_m^{\mathrm{ic}})\}_{m=1}^{N_{\mathrm{ic}}},
\]
where the targets are defined analytically by
\[
h_m^{\mathrm{ic}}
=
\begin{cases}
h_L, & x_m^{\mathrm{ic}}<L/2,\\
h_R, & x_m^{\mathrm{ic}}\ge L/2,
\end{cases}
\qquad
u_m^{\mathrm{ic}}=0.
\]

\rev{
For the final model, the supervised training set uses all $N_x=400$ spatial grid points (excluding the two ghost cells) from each snapshot at $t_k=k\Delta t_s$, $k=1,\ldots,25$, for every $\mu\in P$, without spatial subsampling. Thus, $N_d=65\times25\times400=650{,}000$. The initial-condition set contains all spatial points for every training parameter instance, giving $N_{\mathrm{ic}}=65\times400=26{,}000$, while the boundary set contains both boundary points at every stored positive-time snapshot, giving $N_{\mathrm{bc}}=65\times25\times2=3{,}250$.}

Second, let
\[
\mathcal{D}_{\mathrm{col}}
=
\{(x_j^{(f)},t_j^{(f)},\mu_j^{(f)})\}_{j=1}^{N_f}
\]
denote collocation points for enforcing the PDE residual. These points are sampled in the interior of the space-time-parameter domain. Depending on the experiment, collocation may be sampled uniformly or with a shock-aware strategy. In all cases, the method itself is unchanged: the PDE residual is penalized at collocation points, while data and boundary mismatches are penalized at supervised points.

\paragraph{Shock-aware collocation strategy.}
For wet-bed regimes, strong-form PDE residuals are most difficult to enforce in a narrow neighborhood of the moving shock, where the entropy solution is non-smooth and not classically differentiable. Because strong-form PINNs rely on pointwise spatial and temporal derivatives, PDE enforcement in this region becomes particularly difficult due to high gradients. 
To avoid over-penalizing the network in this region, we optionally replace uniform collocation sampling in \(x\) by a shock-aware strategy that excludes a small band around the shock.
Thus, we first estimate the location of the shock, \(x_{\mathrm{sh}}(t;\muvec)\), from the reference FOM dataset (see Appendix \ref{app:pinn}).
Given the estimated shock position, we then define a spatial exclusion region for collocation points
\begin{equation}
  B_{\mathrm{excl}}(t;\muvec)
  :=
  \left\{
    x\in\Omega:\ |x-x_{\mathrm{sh}}(t;\muvec)|<w_{\mathrm{excl}}
  \right\},
  \label{eq:shock_excl_band}
\end{equation}
where \(w_{\mathrm{excl}}>0\) is a parameter that defines the exclusion half-width. 
The spatial collocation points for the PDE residual are then sampled from
\(
x \sim\mathrm{Unif}\!\left(\Omega\setminus B_{\mathrm{excl}}(t;\muvec)\right).
\)
Consequently, no collocation points are placed inside the excluded shock core \(B_{\mathrm{excl}}\).

For dry-bed or near-dry regimes, shock exclusion is disabled. Specifically, if \(h_R<\varepsilon_{\mathrm{dry}}\) for a prescribed threshold \(\varepsilon_{\mathrm{dry}}>0\), then \(x\) is sampled uniformly on the full spatial domain \(\Omega\). This design reflects the fact that the dry-bed problem is dominated by a rarefaction structure rather than a persistent downstream shock of the wet-bed type. Values of the parameters are reported at the \rev{end of the next subsection}.

\subsection{Loss function}
\label{sec:loss}

The loss function is a weighted sum
\begin{equation}
  \mathcal{L}(\theta)
  =
  \lambda_{\mathrm{data}}\,\mathcal{L}_{\mathrm{data}}
  +
  \lambda_{\mathrm{pde}}\,\mathcal{L}_{\mathrm{pde}}
  +
  \lambda_{\mathrm{bc}}\,\mathcal{L}_{\mathrm{bc}}
  +
  \lambda_{\mathrm{ic}},\mathcal{L}_{\mathrm{ic}}
  +
  \lambda_{\mathrm{nn}}\,\mathcal{L}_{\mathrm{nn}},
  \label{eq:total_loss}
\end{equation}
where \(\mathcal{L}_{\mathrm{data}}\) is the supervised data loss, \(\mathcal{L}_{\mathrm{pde}}\) is the PDE residual loss, \(\mathcal{L}_{\mathrm{bc}}\) is the boundary-condition loss, \(\mathcal{L}_{\mathrm{ic}}\) is the initial-condition loss, and \(\mathcal{L}_{\mathrm{nn}}\) is an optional non-negativity penalty for the depth.

\paragraph{\underline{Data loss.}}
The data loss is computed using 
\begin{equation}
  \mathcal{L}_{\mathrm{data}}
  =
  \frac{1}{N_d}
  \sum_{i=1}^{N_d}
  \left[
    \ell_h\left(e_i^{(h)}\right)
    +
\lambda_u\ell_u\left(e_i^{(u)}\right)
  \right], \quad
  e_i^{(h)} = \frac{h_{\theta,i}-h_i}{\sigma_h},
  \quad
  e_i^{(u)} = \frac{u_{\theta,i}-u_i}{\sigma_u},
  \label{eq:data_loss}
\end{equation}
where $\sigma_h = \text{StdDev}(h_i)$,  $\sigma_u = \text{StdDev}(u_i)$ computed over the entire training dataset $\{(h_i,u_i)\}$, and 
the scalar \(\lambda_u>0\) balances the velocity contribution relative to the depth contribution. Here, \(\ell_h,\ell_u\) denote either \(|\cdot|\) or \((\cdot)^2\), depending on the experiment.

\paragraph{\underline{PDE residual loss.}}
The PDE residual \eqref{eq:r1_method}--\eqref{eq:r2_method} is penalized at chosen collocation points.
Since the two residual components may differ substantially in magnitude, especially during early training, we rescale them component-wise before aggregation:
\begin{equation}
\widetilde{\mathcal{R}}_i={\mathcal{R}_i}/{s_i},
\label{eq:pde_residual_scaling}
\end{equation}
where \(s_i\) are exponential-moving-average root-mean-square (EMA-RMS) scales (see Appendix \ref{app:pinn}).
Hence, the strong-form PDE loss is
\begin{equation}
\mathcal{L}_{\mathrm{pde}}
=
\frac{1}{N_f}
\sum_{j=1}^{N_f}
w_j
\left[
\psi\!\bigl(\widetilde{\mathcal{R}}_1(z_j^{(f)};\theta)\bigr)
+
\psi\!\bigl(\widetilde{\mathcal{R}}_2(z_j^{(f)};\theta)\bigr)
\right],
\label{eq:pde_loss_method}
\end{equation}
where \(w_j\in(0,1]\) is a point-dependent gating weight used to downweight strong-form residual enforcement in regions where the solution is less regular, such as near steep gradients or very small predicted depths. 
Details are presented in Appendix \ref{app:pinn}.

The penalty function \(\psi\) 
in \eqref{eq:pde_loss_method}
is either the quadratic penalty \(\psi(r)=r^2\) or the Charbonnier penalty
\(
  \psi(r)=\sqrt{r^2+\varepsilon^2},
\)
with \(\varepsilon=10^{-5}\), which behaves like \(|r|\) for large residuals but remains smooth near zero. This makes it a differentiable approximation to the absolute value and can improve optimization stability when residual outliers occur.

\paragraph{\underline{Boundary loss.}}
At the left and right boundaries, we impose a zero-gradient Neumann condition on the primitive variables,
\rev{
\[
\partial_x h_\theta=0,
\qquad
\partial_x u_\theta=0,
\qquad x\in\{0,L\}.
\]
}
The corresponding boundary loss is
\begin{equation}
\mathcal{L}_{\mathrm{bc}}
=
\frac{1}{N_{\mathrm{bc}}}
\sum_{k=1}^{N_{\mathrm{bc}}}
\left[
\left(\frac{\partial_x h_\theta(z_k^{\mathrm{bc}})}{\sigma_h}\right)^2
+
\lambda_u
\left(\frac{\partial_x u_\theta(z_k^{\mathrm{bc}})}{\sigma_u}\right)^2
\right].
\label{eq:bc_loss_method}
\end{equation}
The spatial derivatives are evaluated using automatic differentiation. This term stabilizes the learned solution near the edges of the spatial domain and ensures consistency with the reference trajectories used in training.

\paragraph{\underline{Initial-condition loss.}}
The initial dam-break state is imposed explicitly using the analytically constructed dataset \(\mathcal{D}_{\mathrm{ic}}\). The corresponding loss is
\begin{equation*}
\mathcal{L}_{\mathrm{ic}}
=
\frac{1}{N_{\mathrm{ic}}}
\sum_{m=1}^{N_{\mathrm{ic}}}
\left[
\ell_h\left(
\frac{h_\theta(z_m^{\mathrm{ic}})-h_m^{\mathrm{ic}}}
{\sigma_h}
\right)
+
\lambda_u
\ell_u\left(
\frac{u_\theta(z_m^{\mathrm{ic}})-u_m^{\mathrm{ic}}}
{\sigma_u}
\right)
\right],
\end{equation*}
where the penalty functions \(\ell_h\) and \(\ell_u\) are chosen consistently with the supervised data loss. This term enforces the discontinuous initial depth profile and the initially stationary velocity field directly at (t=0).

\paragraph{\underline{Non-negativity penalty.}}
Physical admissibility requires \(h\ge 0\). The implementation therefore supports an additional penalty term that penalizes negative depth predictions on supervised points
\begin{equation*}
  \mathcal{L}_{\mathrm{nn}}
  =
  \frac{1}{N_d}\sum_{i=1}^{N_d}\max(0,-h_{\theta,i}).
\end{equation*}

\paragraph{Optimization.}
The network is trained by minimizing \eqref{eq:total_loss} with Adam optimizer \cite{Kingma2015Adam}. 
We also implemented a cosine-annealing learning-rate schedule and gradient clipping to control occasional large updates. In addition, the PDE, boundary-condition, and initial-condition weights are ramped up during the early training epochs. This warm-up strategy is designed to avoid forcing the network too aggressively toward the physics term before it has learned a reasonable coarse representation of the data. Such staged weighting can improve optimization stability for parametric hyperbolic problems, where the residual landscape may be highly nonuniform across regimes and times.

\paragraph{Final training configuration.}
The final model is a fully connected parametric PINN with  four-dimensional input $(x,t,\bar h, \gamma)$, where $(\bar h, \gamma)$ is a transformation of $\muvec$ in
\eqref{eq:param_reparam}.
The network has five hidden layers of width \(120\), SiLU activation, and a two-dimensional output for \((h,u)\). The network is trained for \(20{,}000\) epochs with Adam using learning rate \(10^{-3}\). A cosine scheduler is used with \(T_{\max}=20000\) and \(\eta_{\min}=10^{-5}\), and gradients are clipped at \(1.0\). The final loss weights are $\lambda_{\mathrm{data}}=1.0,\, \lambda_{\mathrm{pde}}=0.3,\,\lambda_{\mathrm{bc}}=0.6,\,\lambda_{\mathrm{ic}}=1.0,\,\lambda_{\mathrm{nn}}=0.01,$ with an additional velocity weight $\lambda_u = 5.0$. The supervised data and initial-condition terms use \(L^1\)-type penalties, while the PDE residual and boundary-condition terms use \(L^2\)-type penalties. The PDE, boundary-condition, and initial-condition weights are initialized at $10\%$ of their final values and increased linearly to their target values over the first quarter of training. For strong-form residual enforcement, we use $N_{\mathrm{colloc}}^{\mathrm{train}}=80\,000,\, N_{\mathrm{colloc}}^{\mathrm{val}}=20\,000$ collocation points. Shock-aware collocation is enabled with exclusion width \(w_{\mathrm{excl}}=2.5\) \rev{and dry-bed threshold $\epsilon_{\mathrm{dry}}=0.5$.} EMA-RMS residual scaling is used with \(\beta=0.99\) and \(\varepsilon_{\mathrm{scale}}=10^{-8}\). The PDE gating parameters are set to \(\alpha=8.0\), \(h_{\min}=0.05\), and \(w_{\min}=0.05\). The model contains $58\,922$ trainable parameters.

\rev{
For validation, we use the held-out parameter value $\mu_{\mathrm{val}}=(17,4)$, which does not belong to $P$, and generate FOM snapshots with sampling interval $\Delta t_s^{\mathrm{val}}=0.25$. The validation PDE residual is evaluated at $20{,}000$ uniformly sampled collocation points. All training losses are evaluated in full-batch mode. NumPy and PyTorch random-number generators are initialized with seed 42 for ablation study and seeds 42-46 for numerical results in Tables \ref{tab:l2errt1} and \ref{tab:l2errt2}. Visual comparison is presented for the network initialized with seed 42. The network at the final training epoch ($20{,}000$) is used for all reported test evaluations.}

\rev{
All PINN computations were performed in single precision on the CPU of a MacBook Pro equipped with an Apple M1 Max chip with 10 CPU cores and 64~GB of unified memory, using 10 PyTorch threads.}

\section{Numerical Results}
\label{sec:num}

To quantify our results, we use relative $L^2$ errors for $h(x,t)$ and $q(x,t)$. In particular, we define
\begin{equation}
\varepsilon_h^{\text{Pred}}(t;\muvec)
=
\frac{\normLtwo{h^{\text{Pred}}(t;\muvec)-h^{\text{FOM}}(t;\muvec)}}
{\normLtwo{h^{\text{FOM}}(t;\muvec)}},
\qquad
\varepsilon_q^{\text{Pred}}(t;\muvec)
=
\frac{\normLtwo{q^{\text{Pred}}(t;\muvec)-q^{\text{FOM}}(t;\muvec)}}
{\normLtwo{q^{\text{FOM}}(t;\muvec)}},
\label{eq:exp_rel_l2}
\end{equation}
where $h^{\text{FOM}}$ and 
$q^{\text{FOM}}$ stand for \rev{numerical solutions} of the full-order model, and 
$h^{\text{Pred}}$ and 
$q^{\text{Pred}}$ denote reduced model predictions using TROM or PINN.
\rev{Prediction errors in \eqref{eq:exp_rel_l2} measure the fidelity of each surrogate relative to the FOM numerical solution, which defines the parameter-to-solution map $\calS$ that both reduced models are constructed to approximate. The discretization error of the full-order scheme itself is therefore not included in \eqref{eq:exp_rel_l2}. Since the same reference data are used for both surrogates, this error is common to both surrogate evaluations and is distinct from the surrogate approximation error measured here.}
We consider six representative parameter regimes
\begin{equation}
\label{eq:params}
(h_L, h_R) = 
(12, 0), \,\, (12, 7), \,\, (15, 4), \,\, (18, 0), \,\, (26, 0.14), \,\, (26, 7)
\end{equation}
evaluated at $t\in\{0.5, 1, 1.5, 2\}$. These cases include both dry-bed and wet-bed settings and are used to compare qualitative field predictions against the FOM simulations.
All parameter regimes in \eqref{eq:params} are out-of-sample with respect to the parameter pairs used for training.
\rev{To assess the sensitivity of the PINN predictions to random initialization and stochastic sampling, the final PINN configuration was trained independently using five random seeds, 42-46. Accordingly, the PINN errors reported in Tables \ref{tab:l2errt1} and \ref{tab:l2errt2} are given as the mean $\pm$ sample standard deviation over these five runs. The relatively small variation across the independent PINN runs indicates that the principal accuracy trends are stable with respect to the random seed.}

\begin{table}[H]
\centering
\caption{Relative $L^2$ Errors in \eqref{eq:exp_rel_l2} for the TROM and the PINN predictions at time $t=1$. \rev{PINN results are reported as mean $\pm$ sample standard deviation over five independent training runs.}}
\label{tab:l2errt1}
\begin{tabular}{lcc|cc}
\toprule
$\muvec$ & $\varepsilon_h^{\text{TROM}}$ & $\varepsilon_h^{\text{PINN}}$ &
$\varepsilon_q^{\text{TROM}}$ &
$\varepsilon_q^{\text{PINN}}$ \\
\midrule
(12, 0)      & $1.4027\times 10^{-3}$ & \rev{$(6.30 \pm 1.22)\times 10^{-3}$} 
            & $1.1859\times 10^{-2}$  & \rev{$(4.24 \pm 0.66)\times 10^{-2}$} \\
(12, 7)      & $4.0641\times 10^{-3}$ & \rev{$(2.66 \pm 0.60)\times 10^{-3}$} 
            & $4.1335\times 10^{-2}$ & \rev{$(1.409 \pm 0.270)\times 10^{-2}$} \\
(15, 4)      & $3.2675\times 10^{-3}$ & \rev{$(3.18 \pm 1.63)\times 10^{-3}$}
            & $1.9264\times 10^{-2}$ & \rev{$(7.16 \pm 1.15)\times 10^{-3}$} \\
(18, 0)      & $3.5761\times 10^{-4}$ & \rev{$(7.16 \pm 1.18)\times 10^{-3}$} 
            & $2.8654\times 10^{-3}$ & \rev{$(5.11 \pm 0.85)\times 10^{-2}$} \\
(26, 0.14)   & $5.6640\times 10^{-3}$ & \rev{$(1.494 \pm 0.209)\times 10^{-2}$}
            & $3.3705\times 10^{-2}$ & \rev{$(8.11 \pm 1.22)\times 10^{-2}$} \\
(26, 7)      & $6.0776\times 10^{-3}$ & \rev{$(2.44 \pm 1.59)\times 10^{-3}$}
            & $2.9323\times 10^{-2}$ & \rev{$(5.88 \pm 2.88)\times 10^{-3}$} \\
\bottomrule
\end{tabular}
\end{table}
\begin{table}[H]
\centering
\caption{Relative $L^2$ Errors in \eqref{eq:exp_rel_l2} for the TROM and the PINN predictions at time $t=2$. \rev{PINN results are reported as mean $\pm$ sample standard deviation over five independent training runs.}}
\label{tab:l2errt2}
\begin{tabular}{lcc|cc}
\toprule
$\muvec$ & $\varepsilon_h^{\text{TROM}}$ & $\varepsilon_h^{\text{PINN}}$ &
$\varepsilon_q^{\text{TROM}}$ &
$\varepsilon_q^{\text{PINN}}$ \\
\midrule
(12, 0)      & $2.5749\times 10^{-3}$ & \rev{$(8.55 \pm 1.27)\times 10^{-3}$}
            & $1.5302\times 10^{-2}$ & \rev{$(4.80 \pm 0.71)\times 10^{-2}$} \\
(12, 7)      & $1.1370\times 10^{-2}$ & \rev{$(2.59 \pm 0.48)\times 10^{-3}$} 
            & $8.1167\times 10^{-2}$ & \rev{$(9.66 \pm 0.78)\times 10^{-3}$} \\
(15, 4)      & $9.3945\times 10^{-3}$ & \rev{$(4.64 \pm 2.34)\times 10^{-3}$} 
            & $3.7016\times 10^{-2}$ & \rev{$(5.95 \pm 1.93)\times 10^{-3}$} \\
(18, 0)      & $7.2675\times 10^{-4}$ & \rev{$(9.56 \pm 1.36)\times 10^{-3}$} 
            & $4.1002\times 10^{-3}$ & \rev{$(4.98 \pm 0.71)\times 10^{-2}$} \\
(26, 0.14)   & $1.5224\times 10^{-2}$ & \rev{$(2.624 \pm 0.308)\times 10^{-2}$} 
            & $6.0324\times 10^{-2}$ & \rev{$(1.008 \pm 0.106)\times 10^{-1}$} \\
(26, 7)      & $2.0034\times 10^{-2}$ & \rev{$(3.86 \pm 1.66)\times 10^{-3}$} 
            & $6.4883\times 10^{-2}$ & \rev{$(4.66 \pm 2.34)\times 10^{-3}$} \\
\bottomrule
\end{tabular}
\end{table}

We depict a comparison of the water depth $h(x,t;\muvec)$
and the discharge $q(x,t;\muvec)$
in Figures \ref{fig:profh1}, \ref{fig:profh2} and 
\ref{fig:profq1}, \ref{fig:profq2}, respectively.
\def\figwidth{0.4}
\def\figheight{2in}
\begin{figure}[H]  
\centerline{\includegraphics[width=\figwidth\textwidth,height=\figheight]{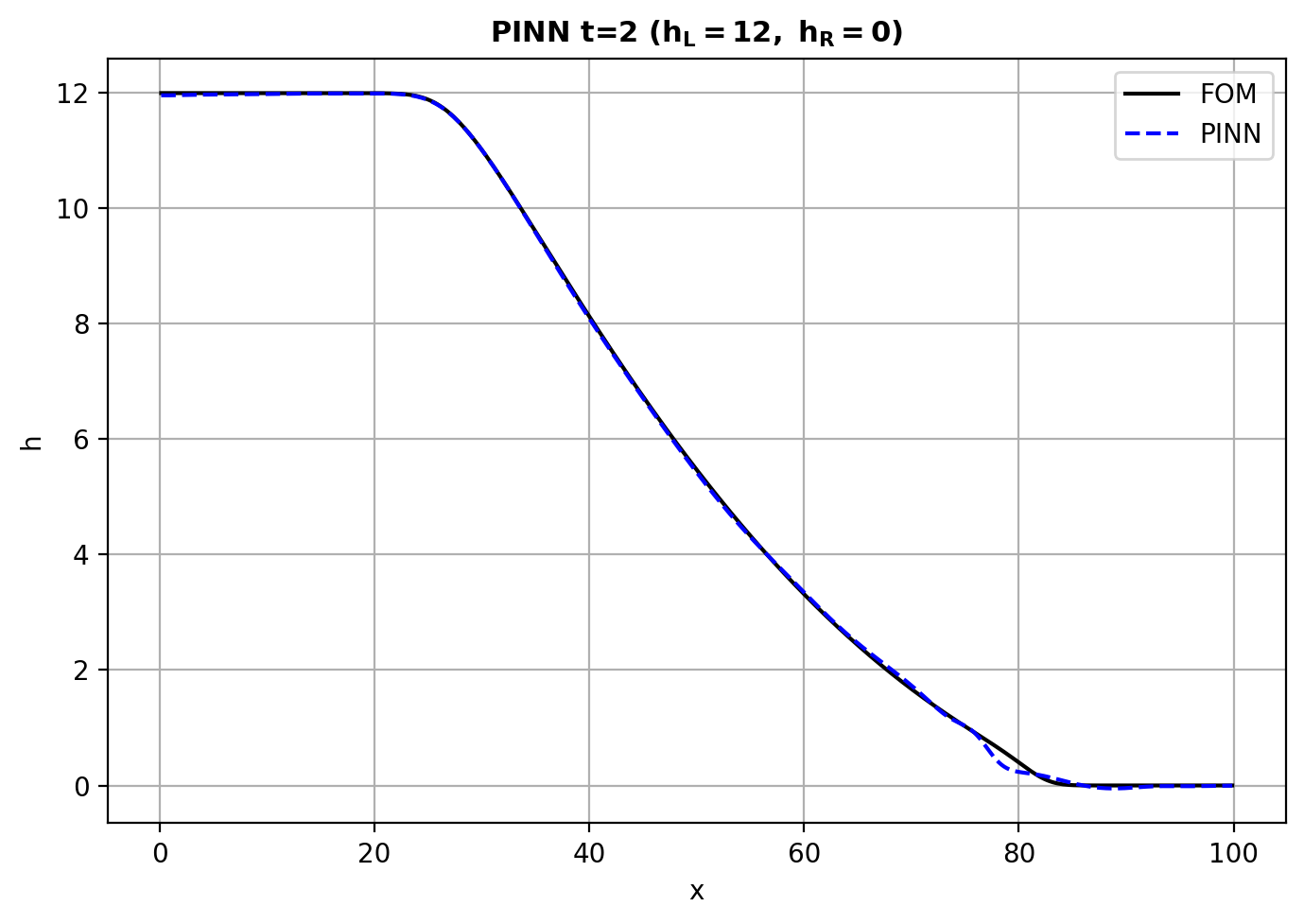}
\includegraphics[width=\figwidth\textwidth,height=\figheight]{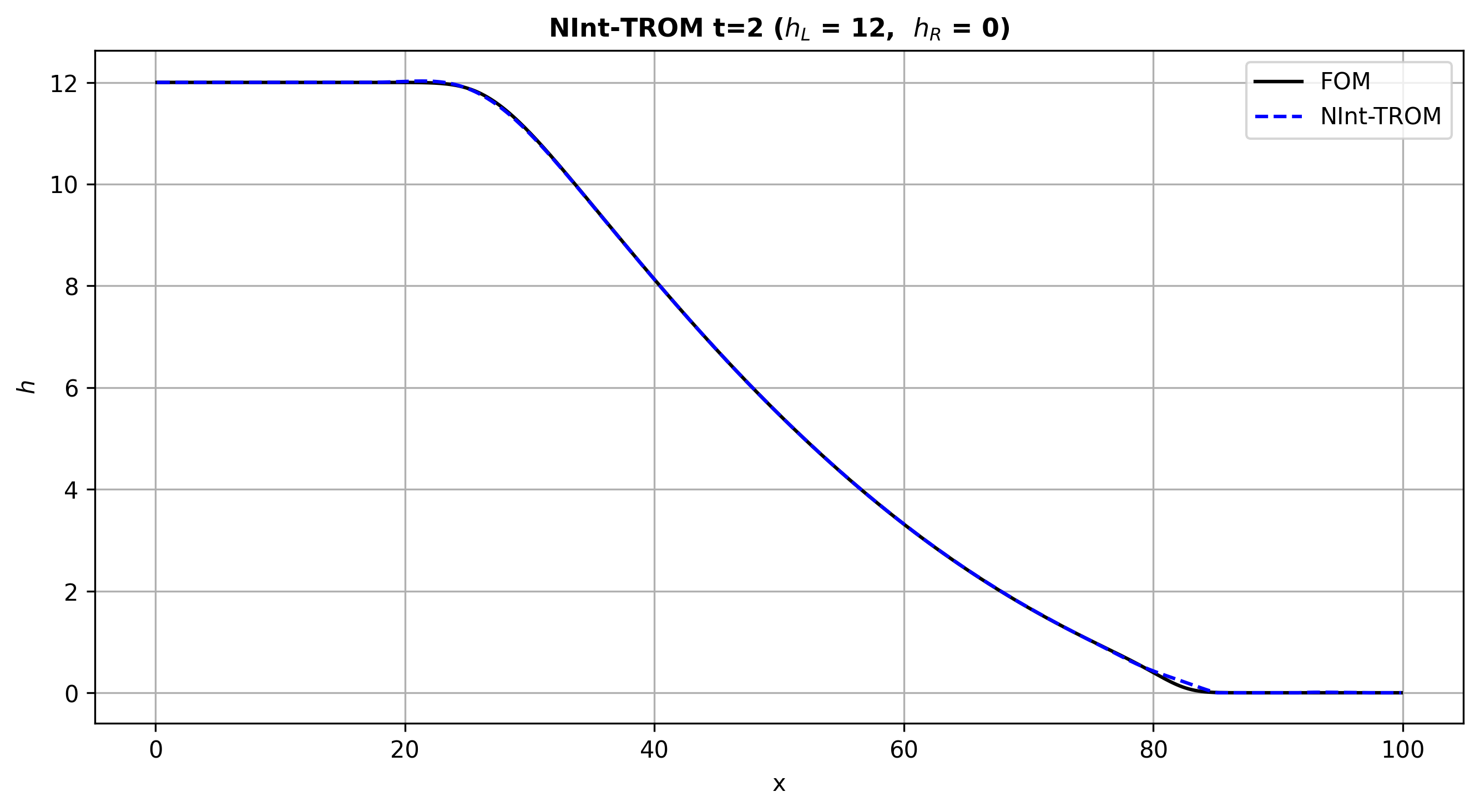}}
\centerline{\includegraphics[width=\figwidth\textwidth,height=\figheight]{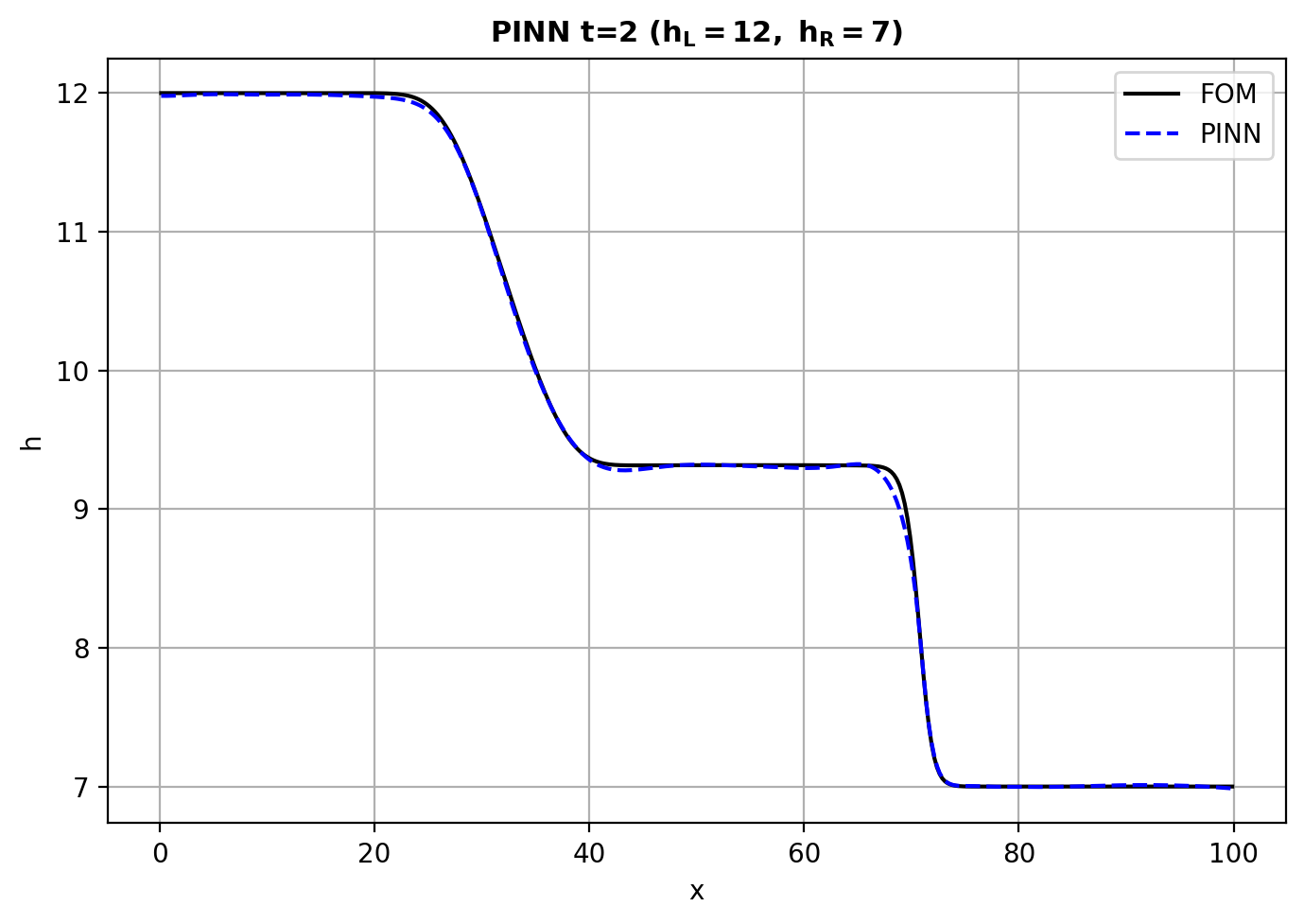}
\includegraphics[width=\figwidth\textwidth,height=\figheight]{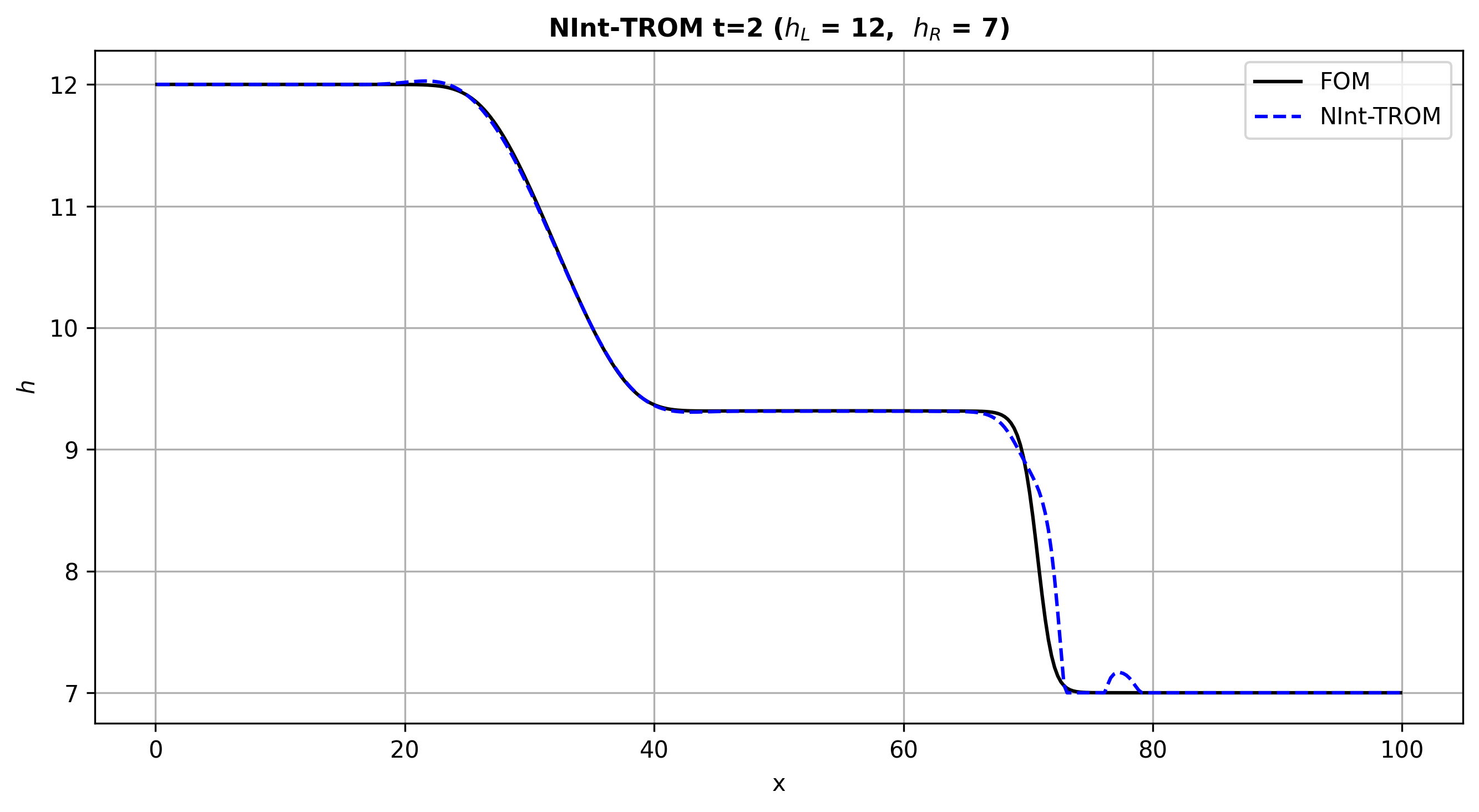}}
\centerline{\includegraphics[width=\figwidth\textwidth,height=\figheight]{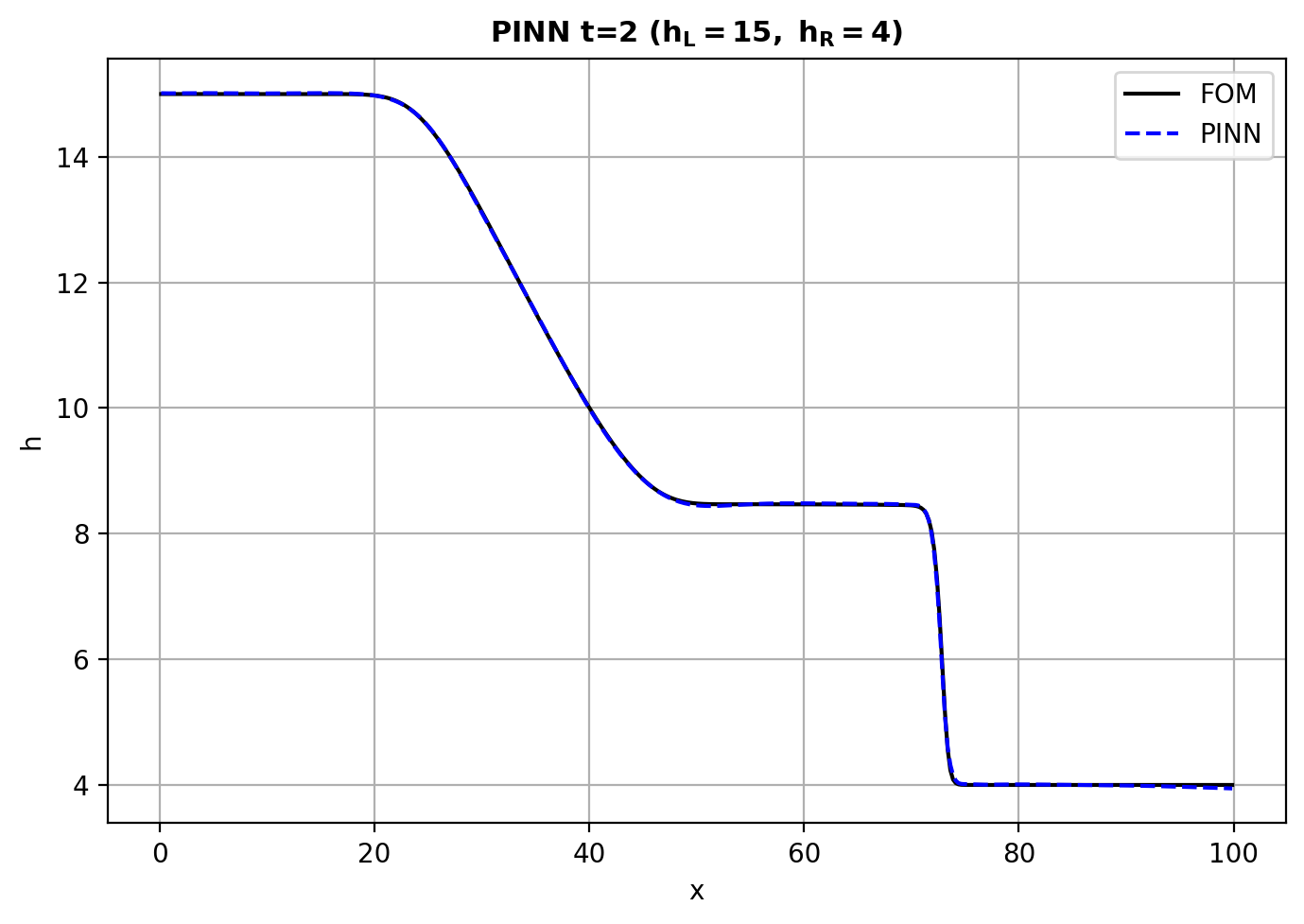}
\includegraphics[width=\figwidth\textwidth,height=\figheight]{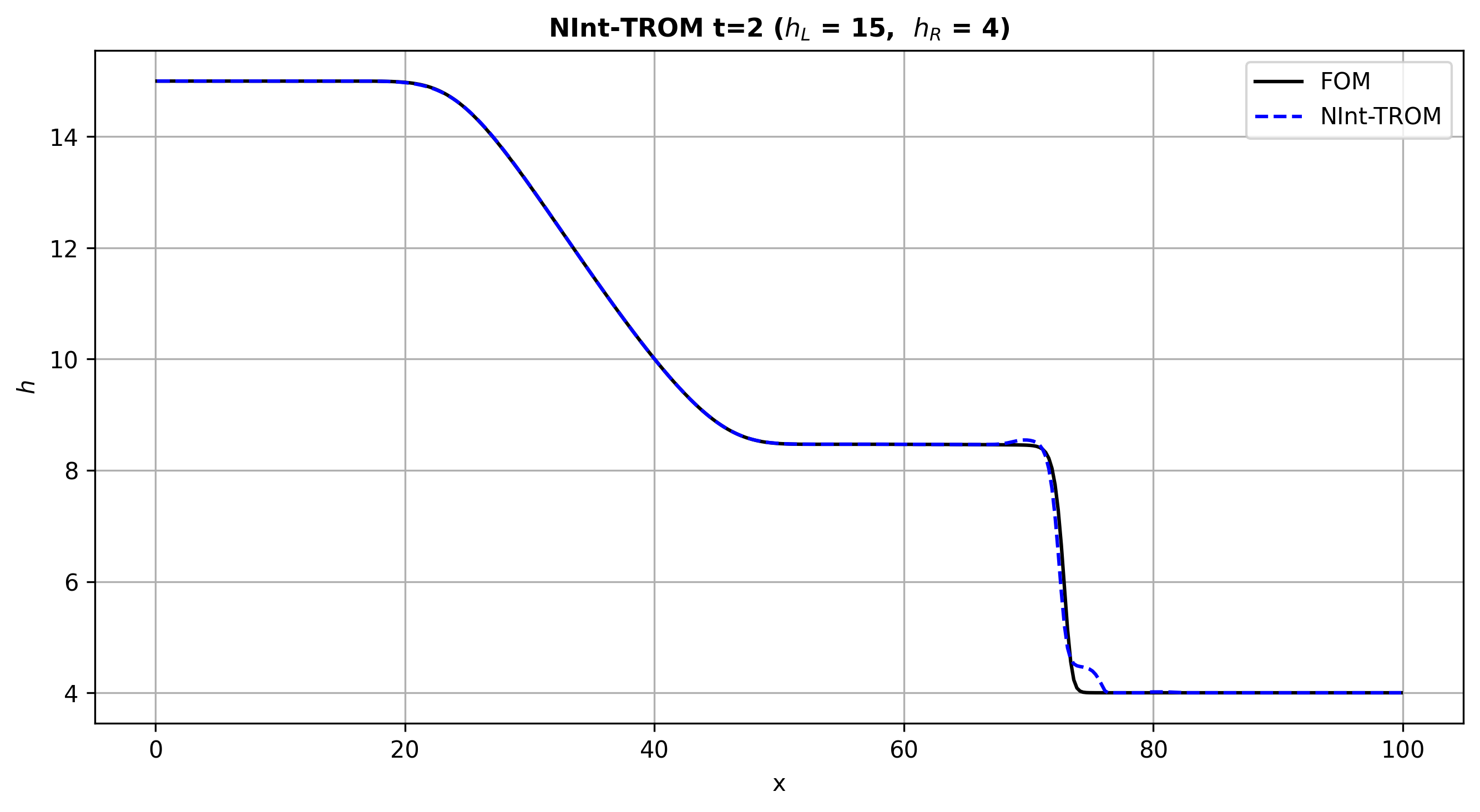}}
\caption{Comparison of $h(x,t;\muvec)$ from the FOM simulations and predictions of the PINN (left column) and the TROM (right column) for $\muvec=(12, 0), (12, 7), (15, 4)$ (top to bottom) at time $t=2$.}
\label{fig:profh1}
\end{figure}
\begin{figure}[H]  
\centerline{\includegraphics[width=\figwidth\textwidth,height=\figheight]{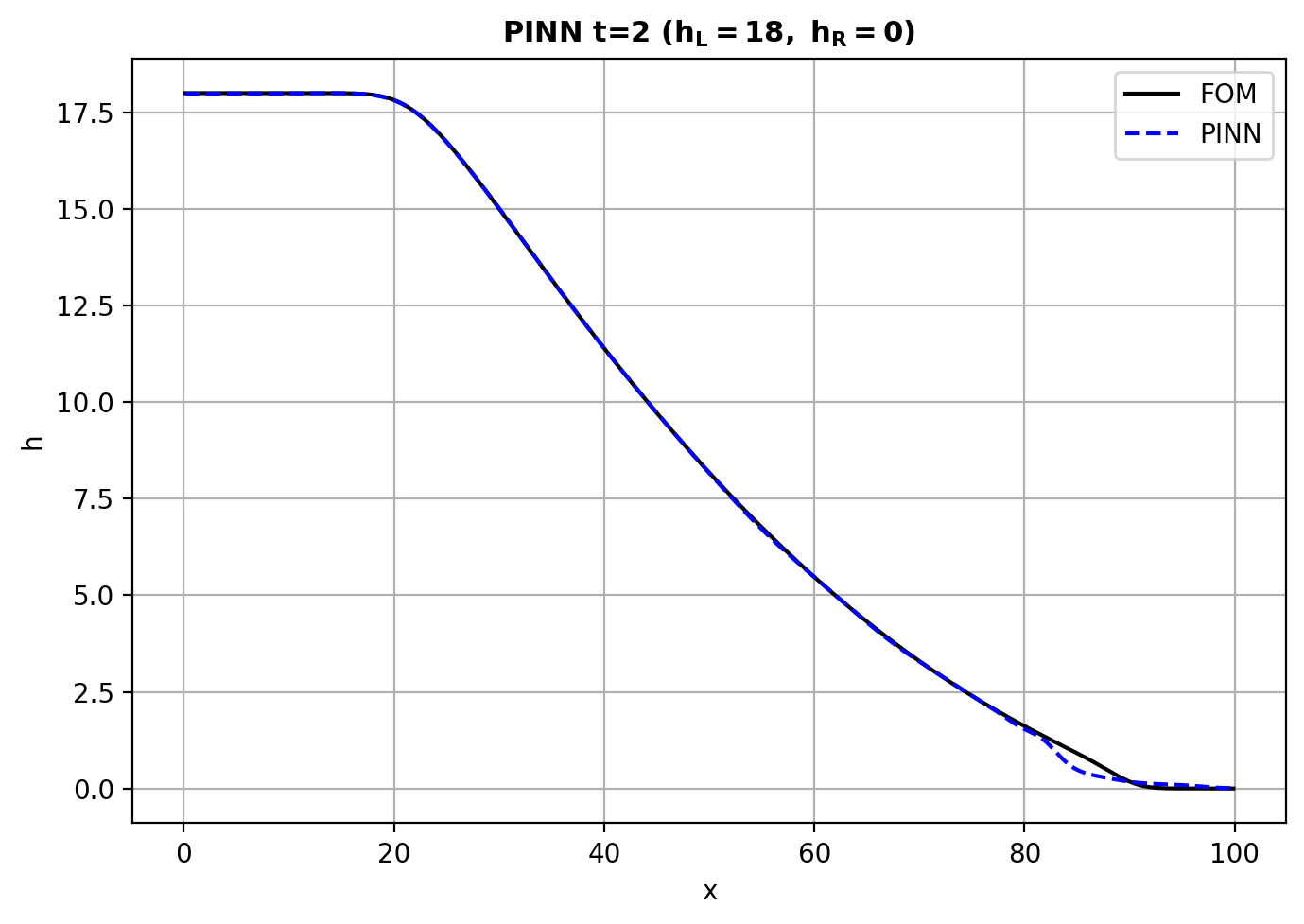}
\includegraphics[width=\figwidth\textwidth,height=\figheight]{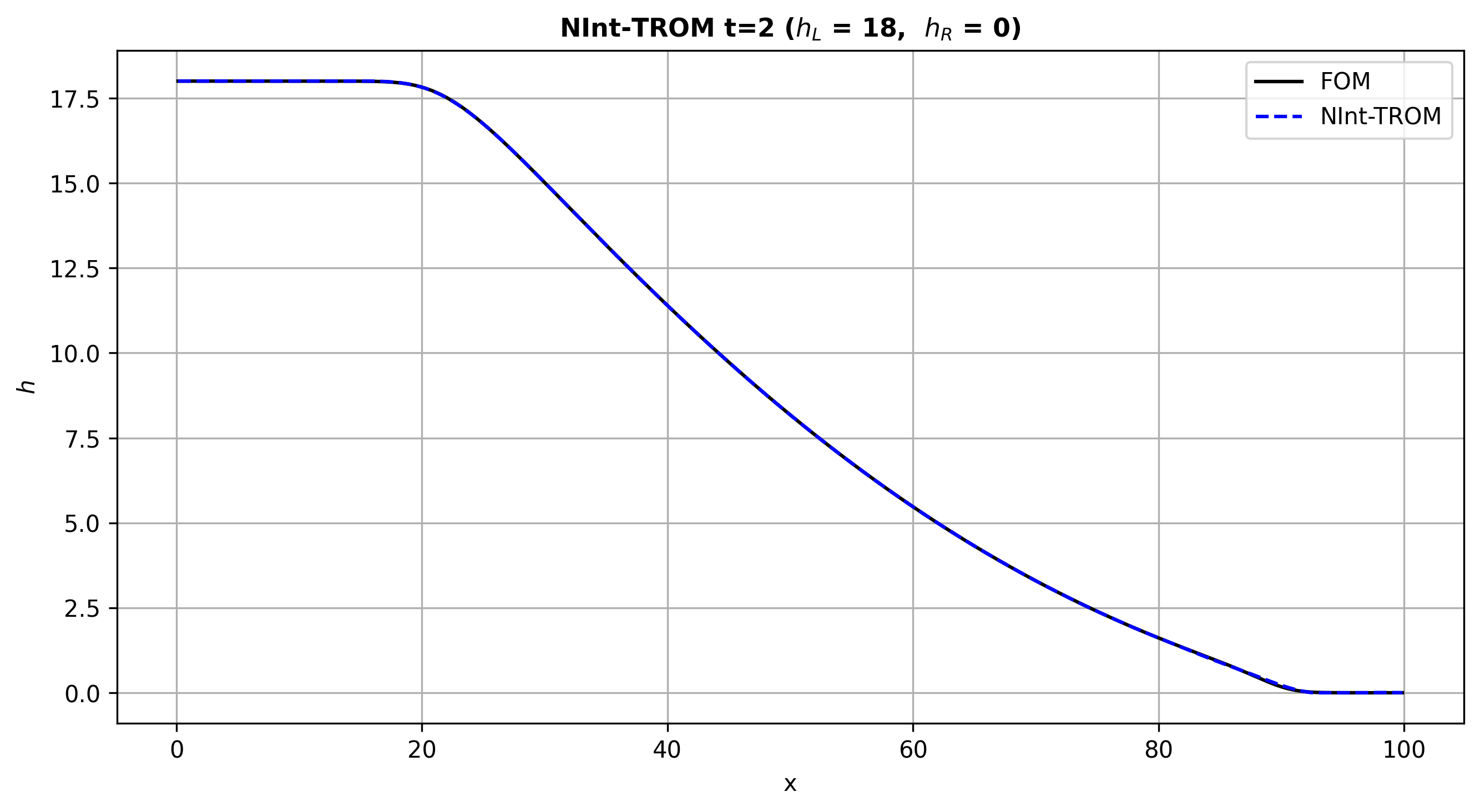}}
\centerline{\includegraphics[width=\figwidth\textwidth,height=\figheight]{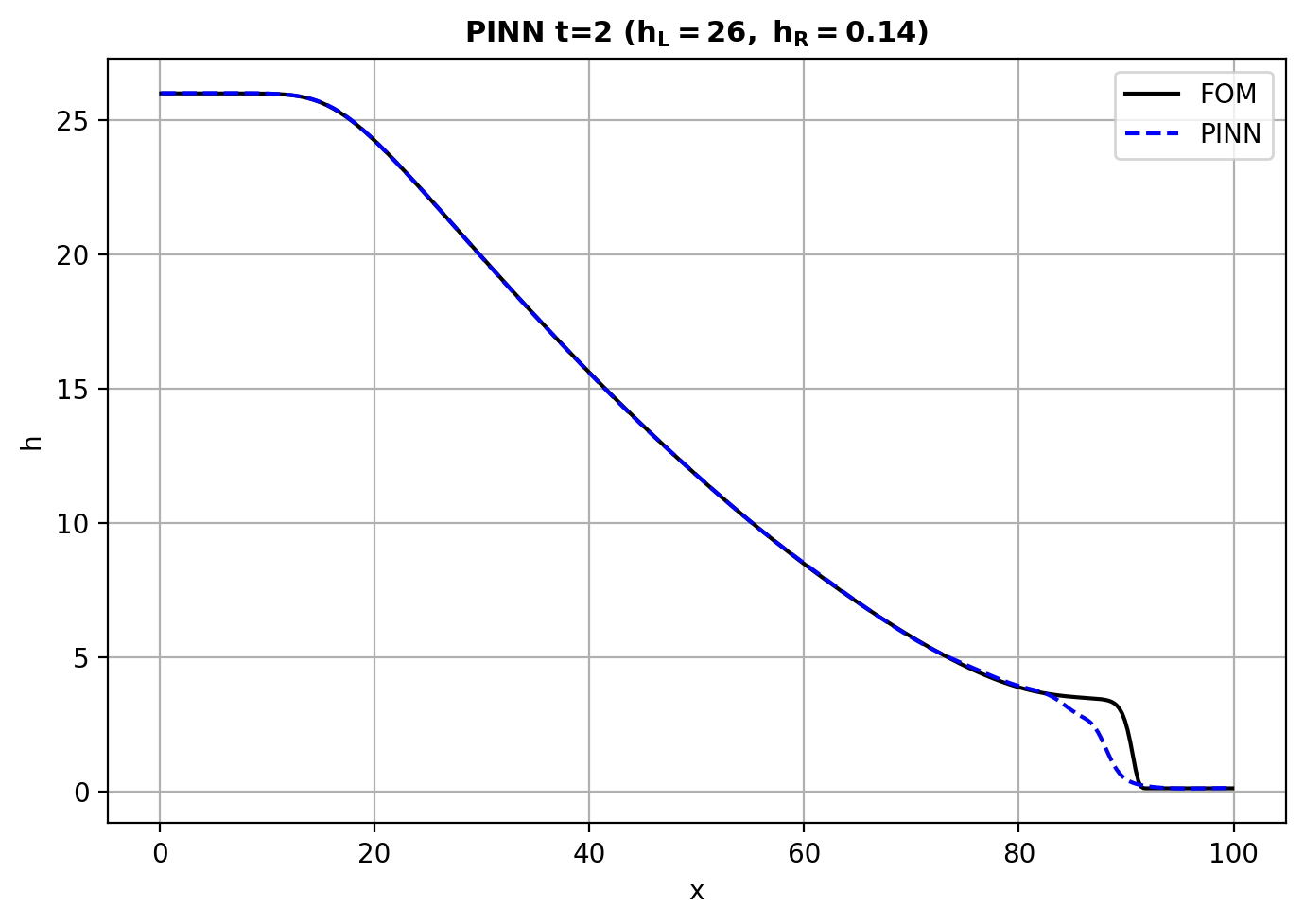}
\includegraphics[width=\figwidth\textwidth,height=\figheight]{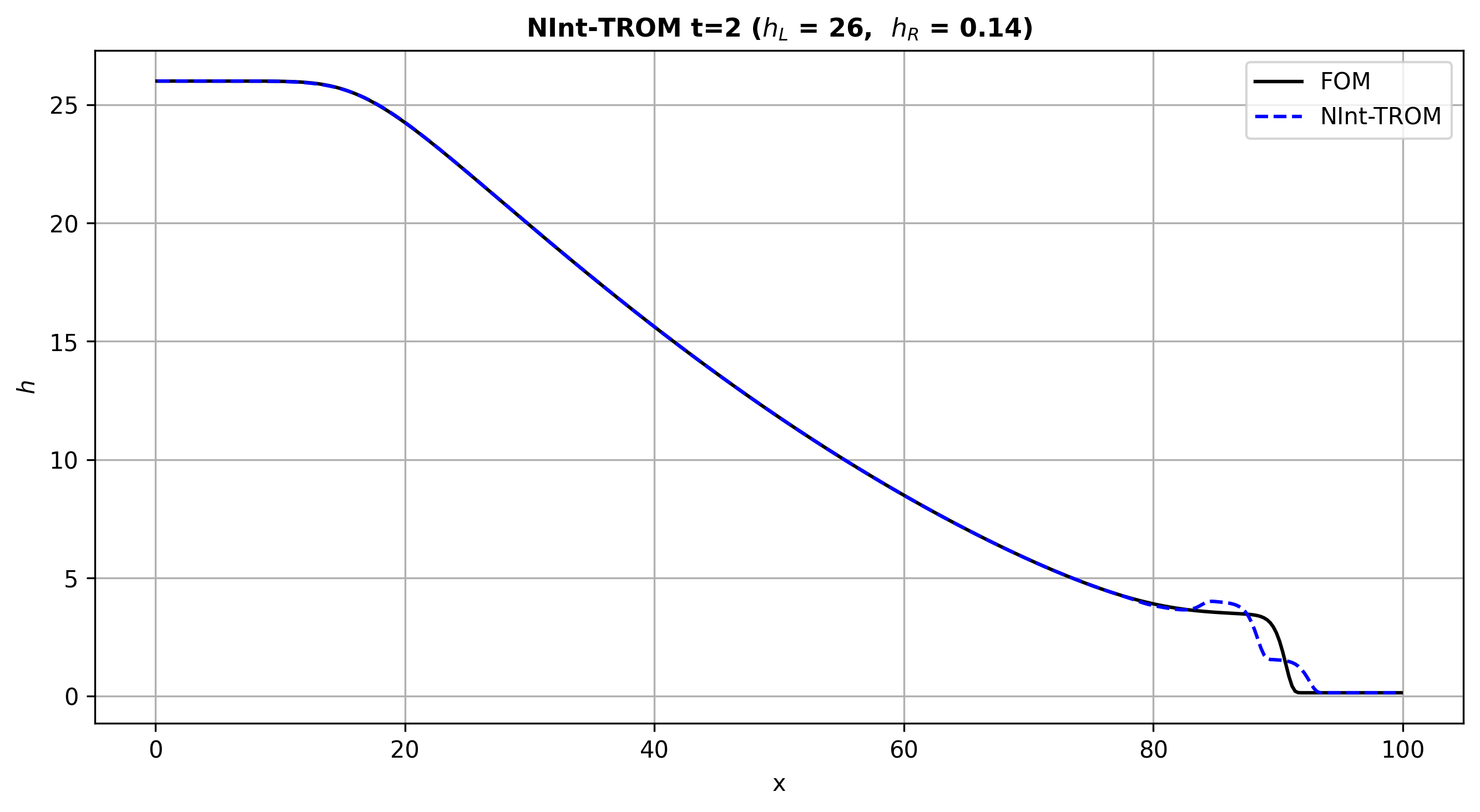}}
\centerline{\includegraphics[width=\figwidth\textwidth,height=\figheight]{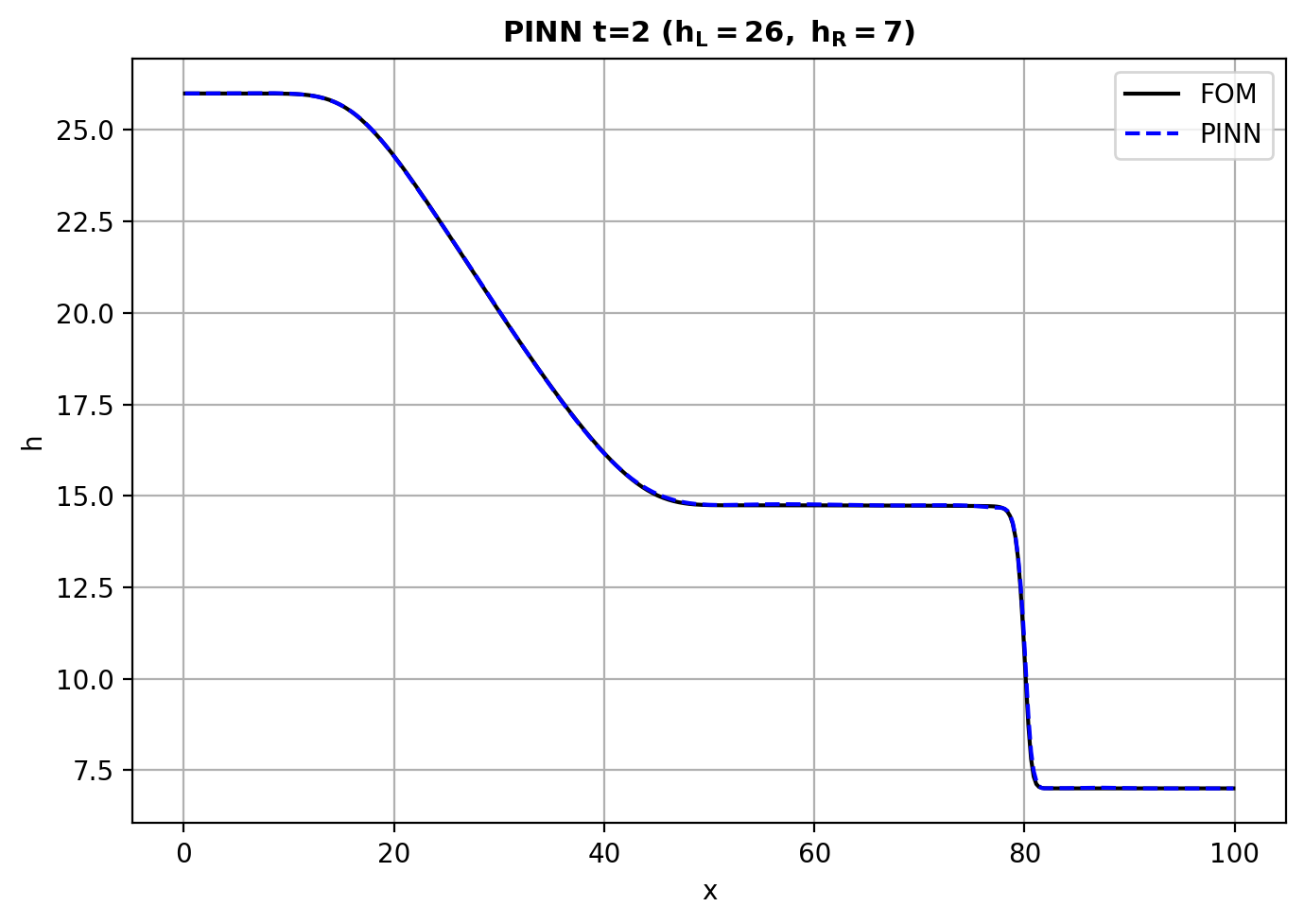}
\includegraphics[width=\figwidth\textwidth,height=\figheight]{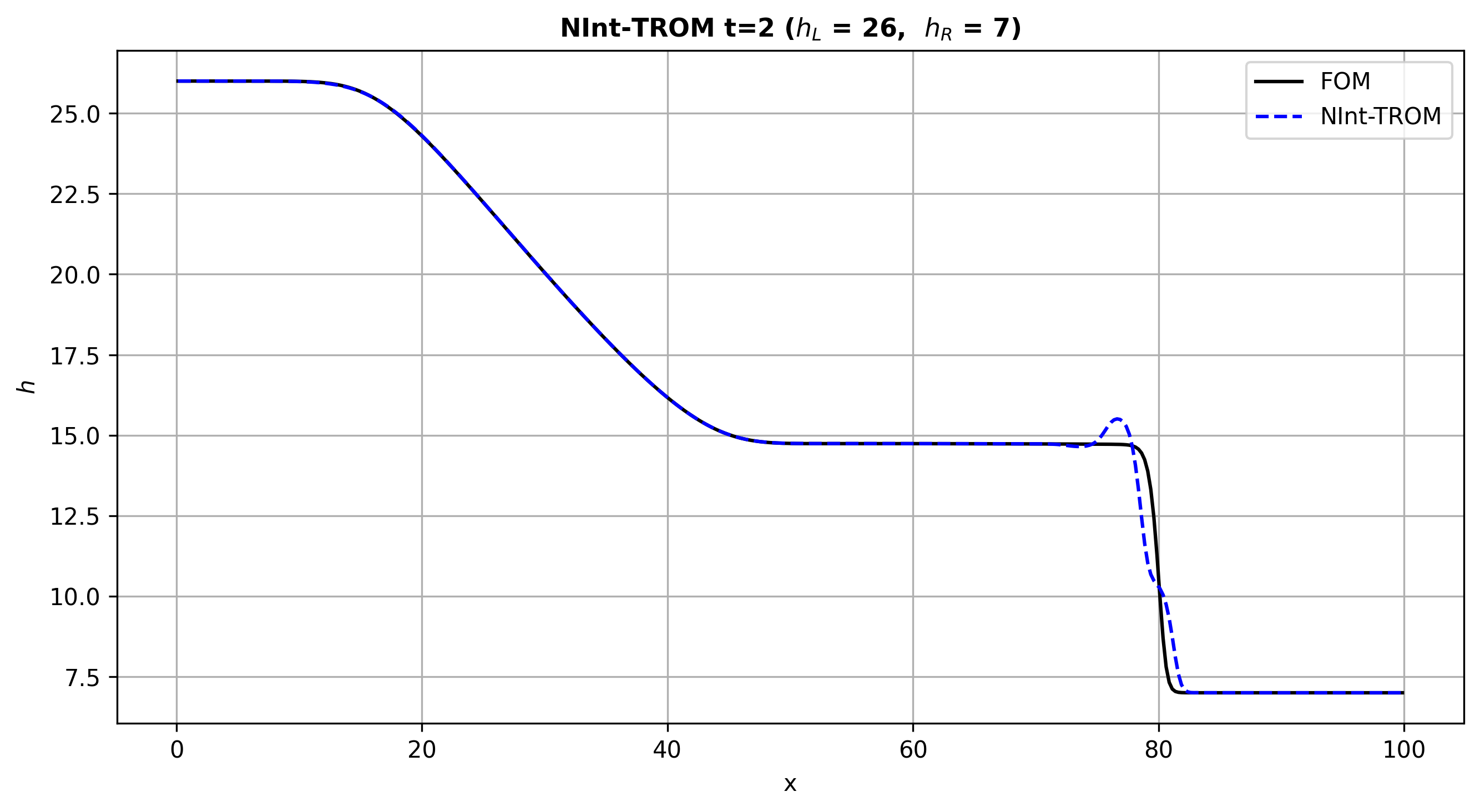}}
\caption{Comparison of $h(x,t;\muvec)$ from the FOM simulations and predictions of the PINN (left column) and the TROM (right column) for $\muvec=(18, 0), (26, 0.14), (26, 7)$ (top to bottom) at time $t=2$.}
\label{fig:profh2}
\end{figure}

\begin{figure}[H]  
\centerline{\includegraphics[width=\figwidth\textwidth,height=\figheight]{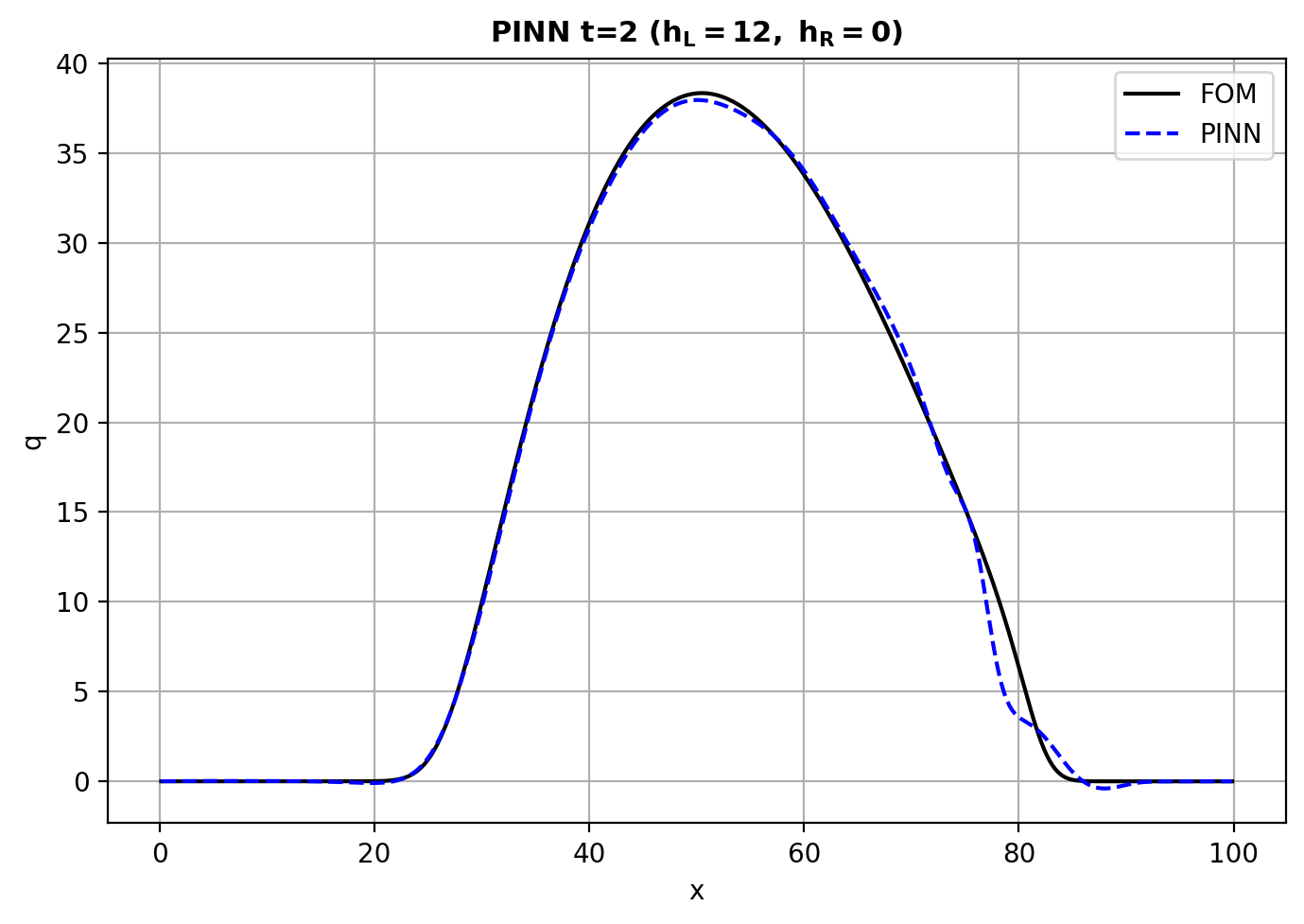}
\includegraphics[width=\figwidth\textwidth,height=\figheight]{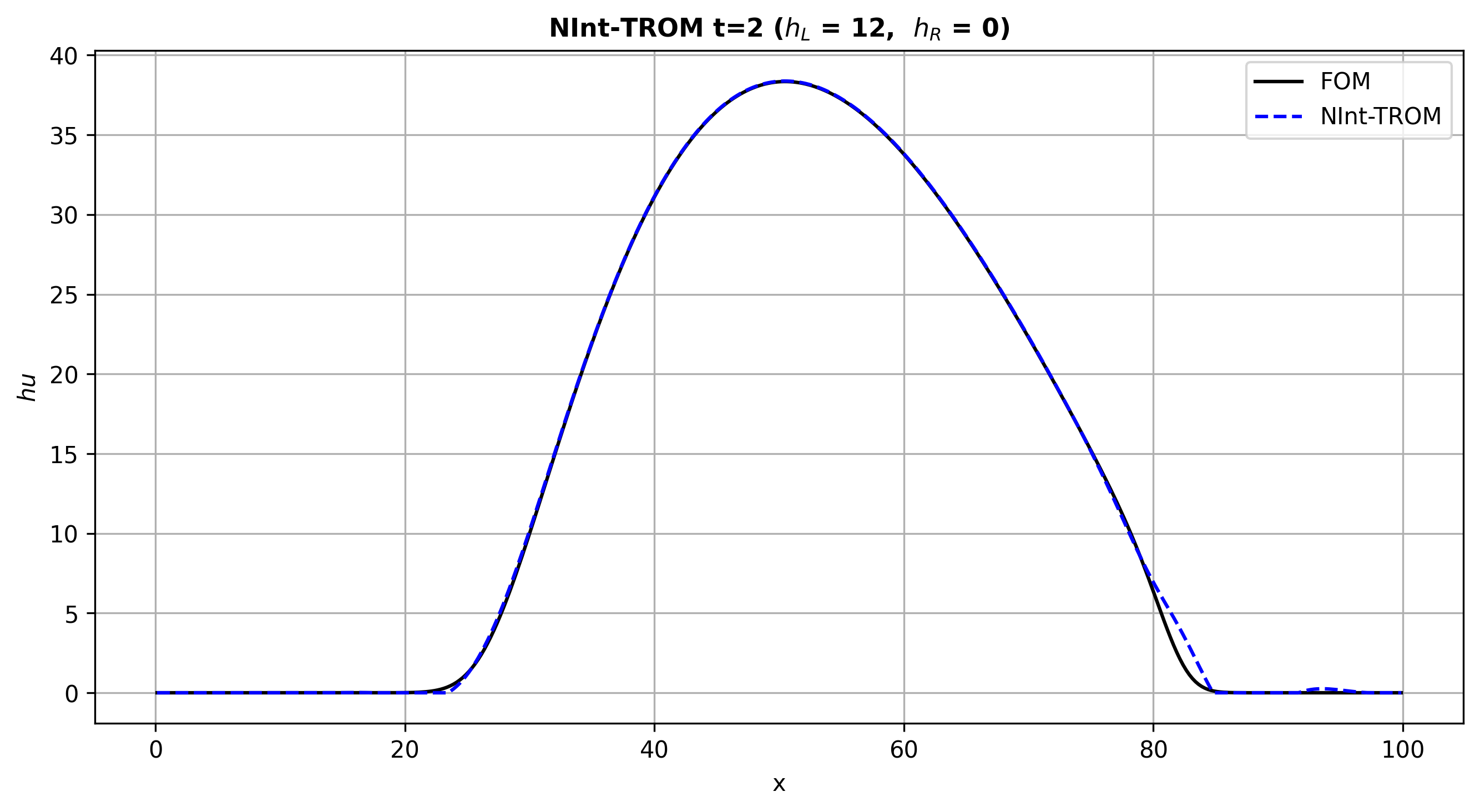}}
\centerline{\includegraphics[width=\figwidth\textwidth,height=\figheight]{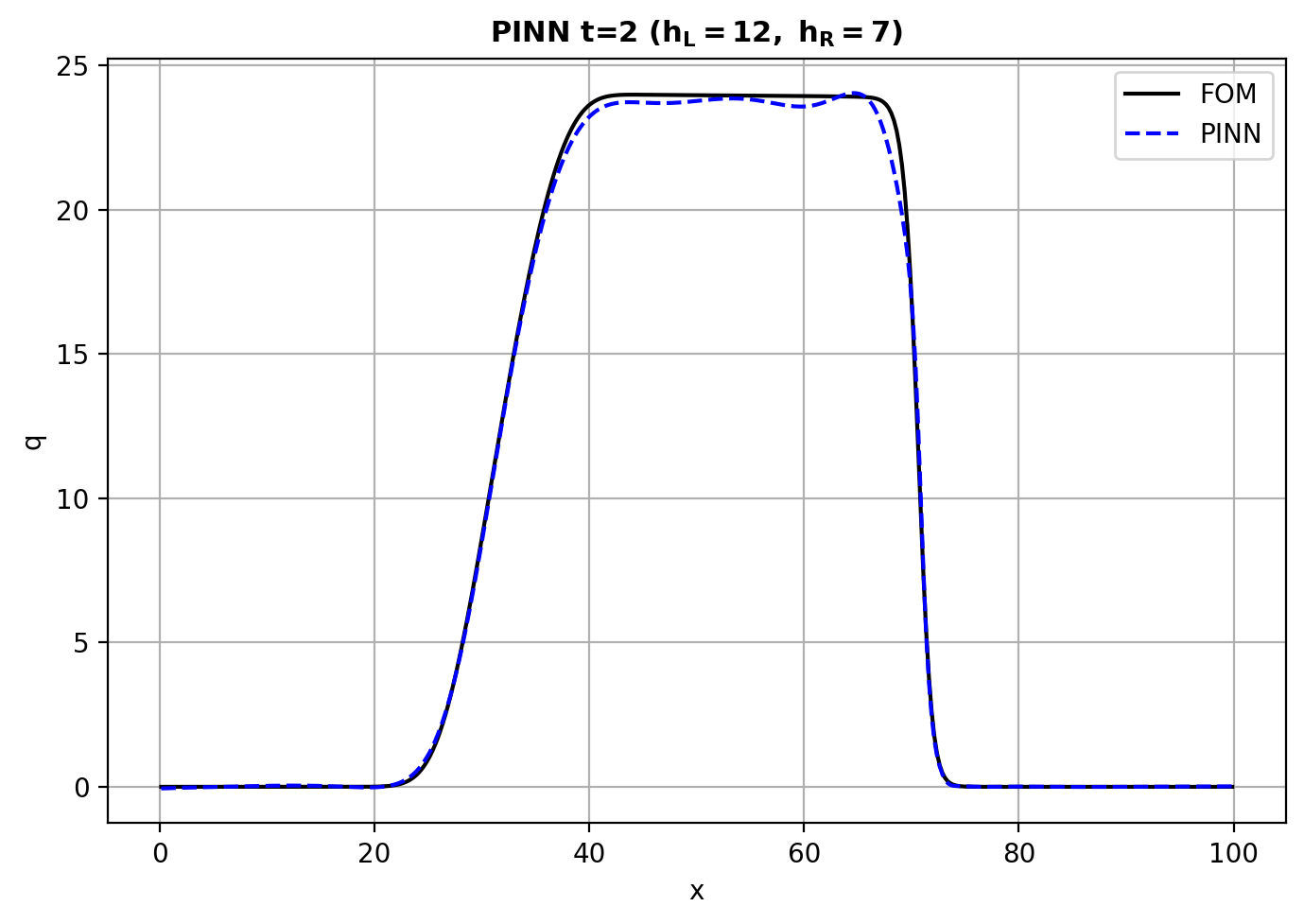}
\includegraphics[width=\figwidth\textwidth,height=\figheight]{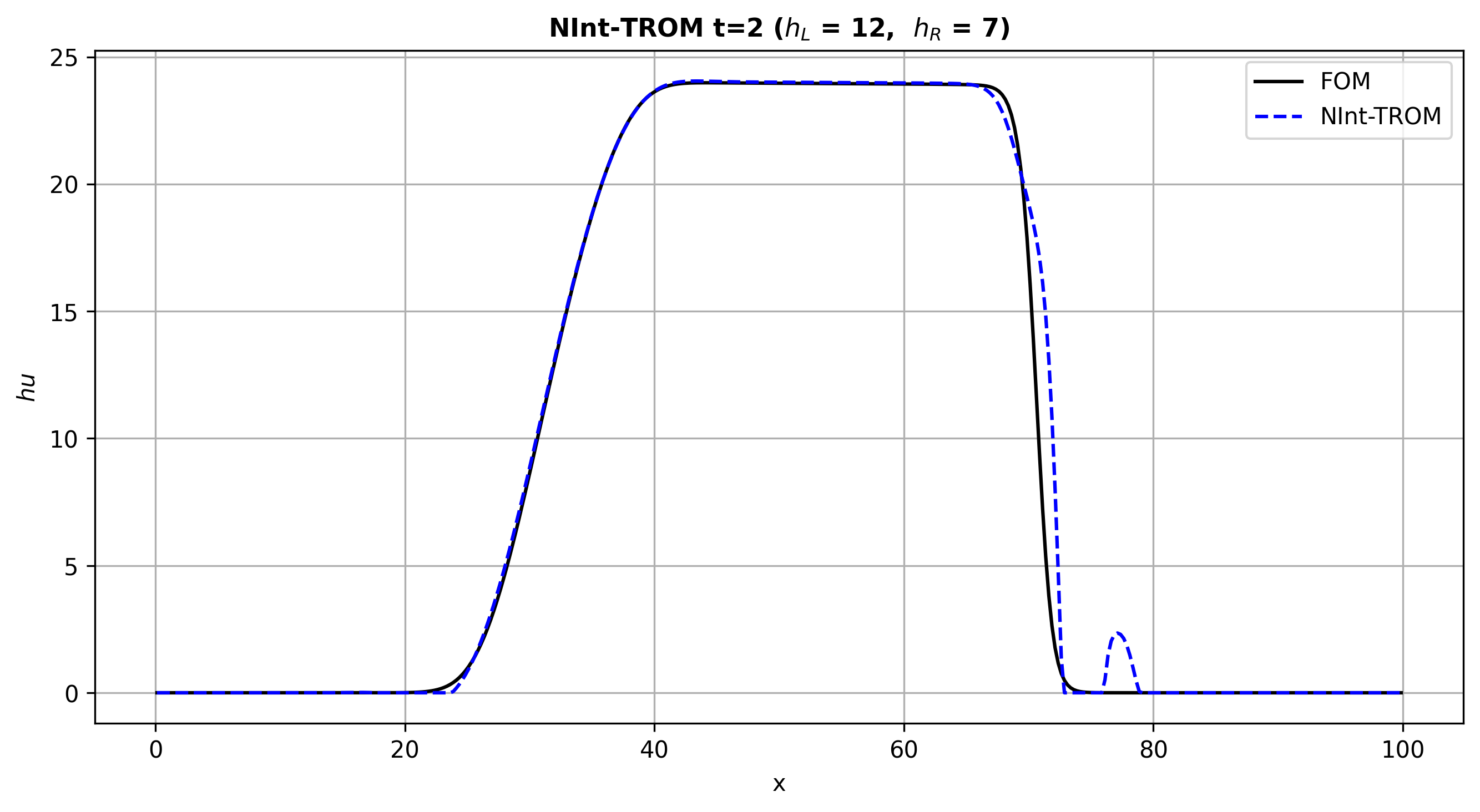}}
\centerline{\includegraphics[width=\figwidth\textwidth,height=\figheight]{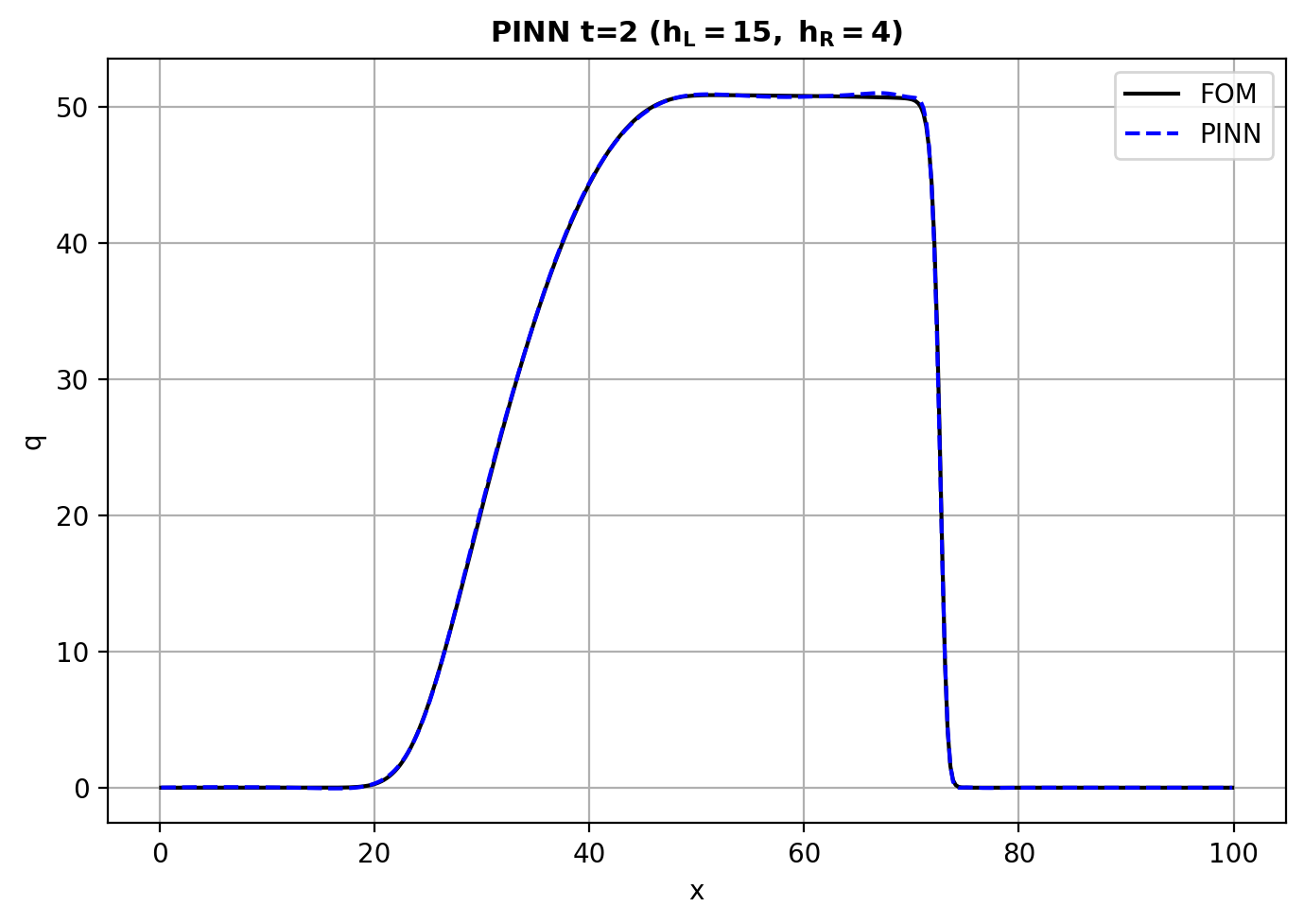}
\includegraphics[width=\figwidth\textwidth,height=\figheight]{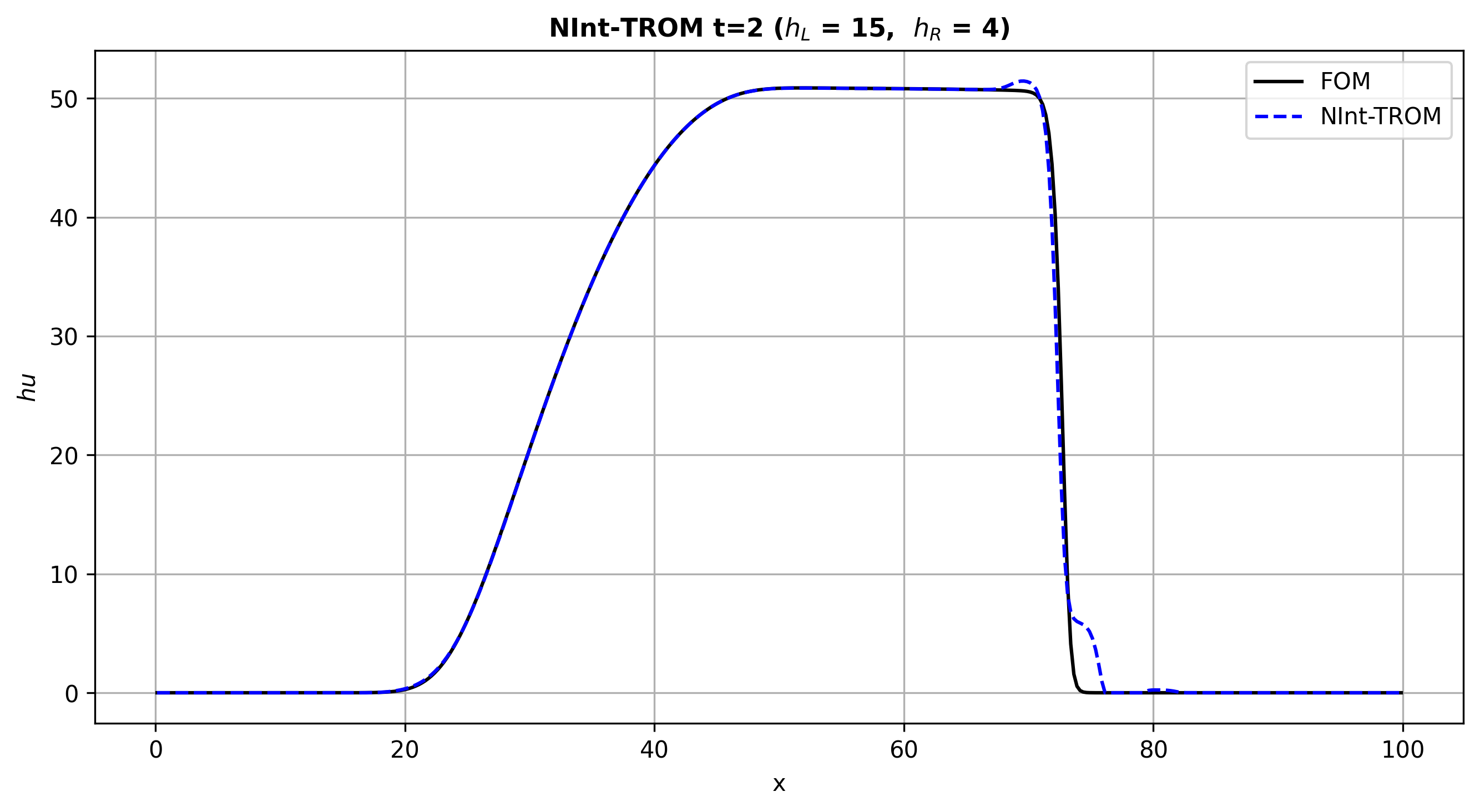}}
\caption{Comparison of $q(x,t;\muvec)$ from the FOM simulations and predictions of the PINN (left column) and the TROM (right column) for $\muvec=(12, 0), (12, 7), (15, 4)$ (top to bottom) at time $t=2$.}
\label{fig:profq1}
\end{figure}
\begin{figure}[H]  
\centerline{\includegraphics[width=\figwidth\textwidth,height=\figheight]{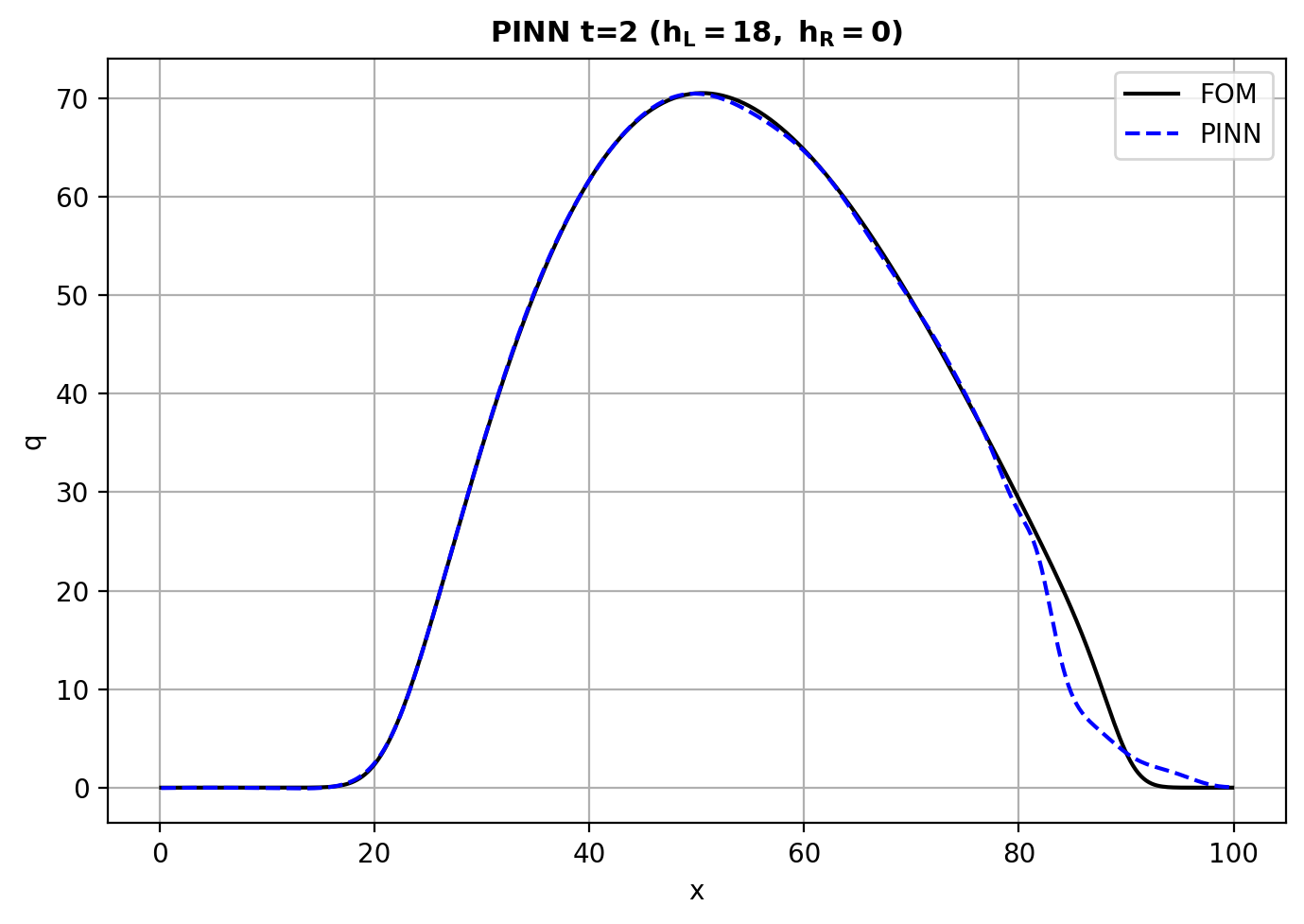}
\includegraphics[width=\figwidth\textwidth,height=\figheight]{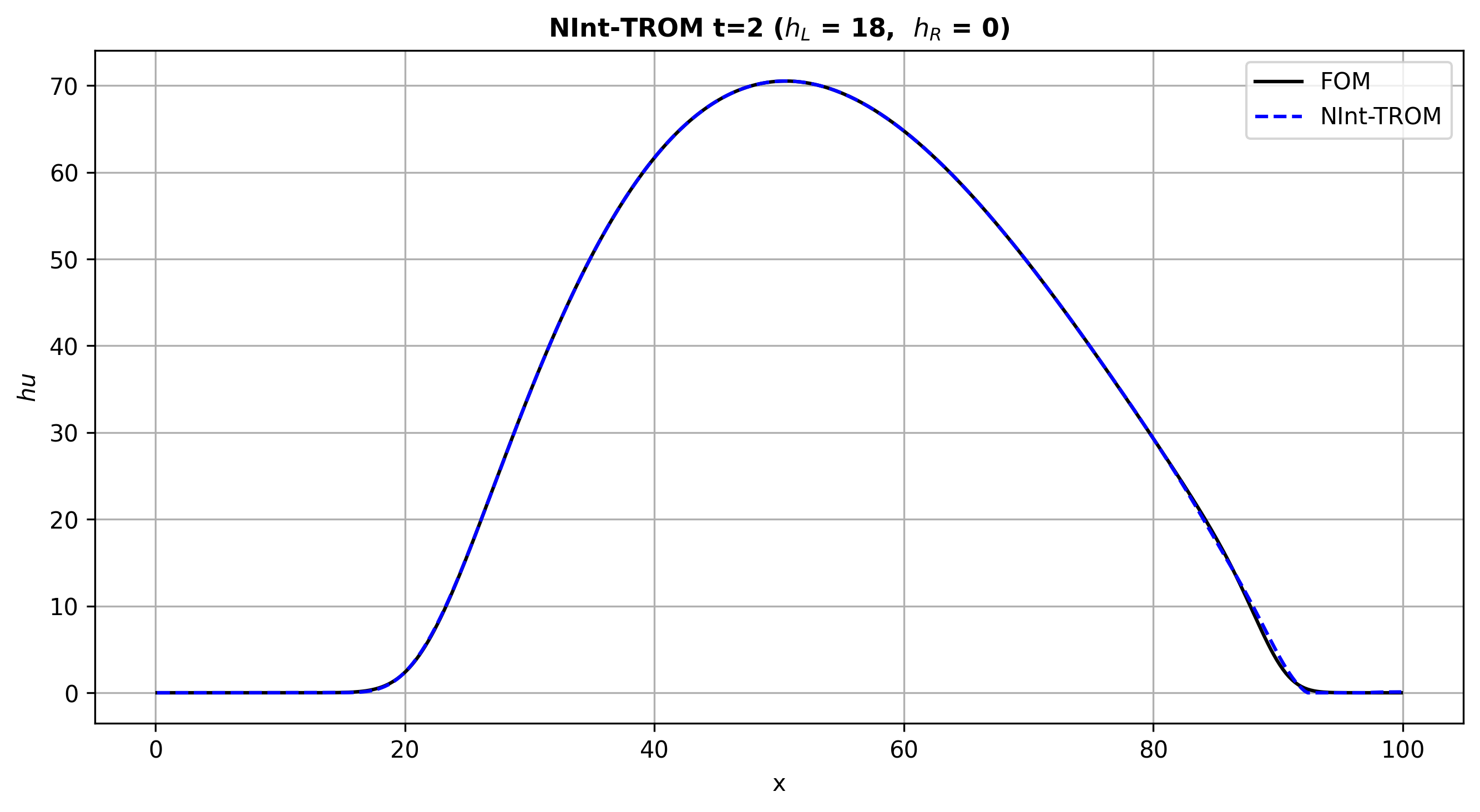}}
\centerline{\includegraphics[width=\figwidth\textwidth,height=\figheight]{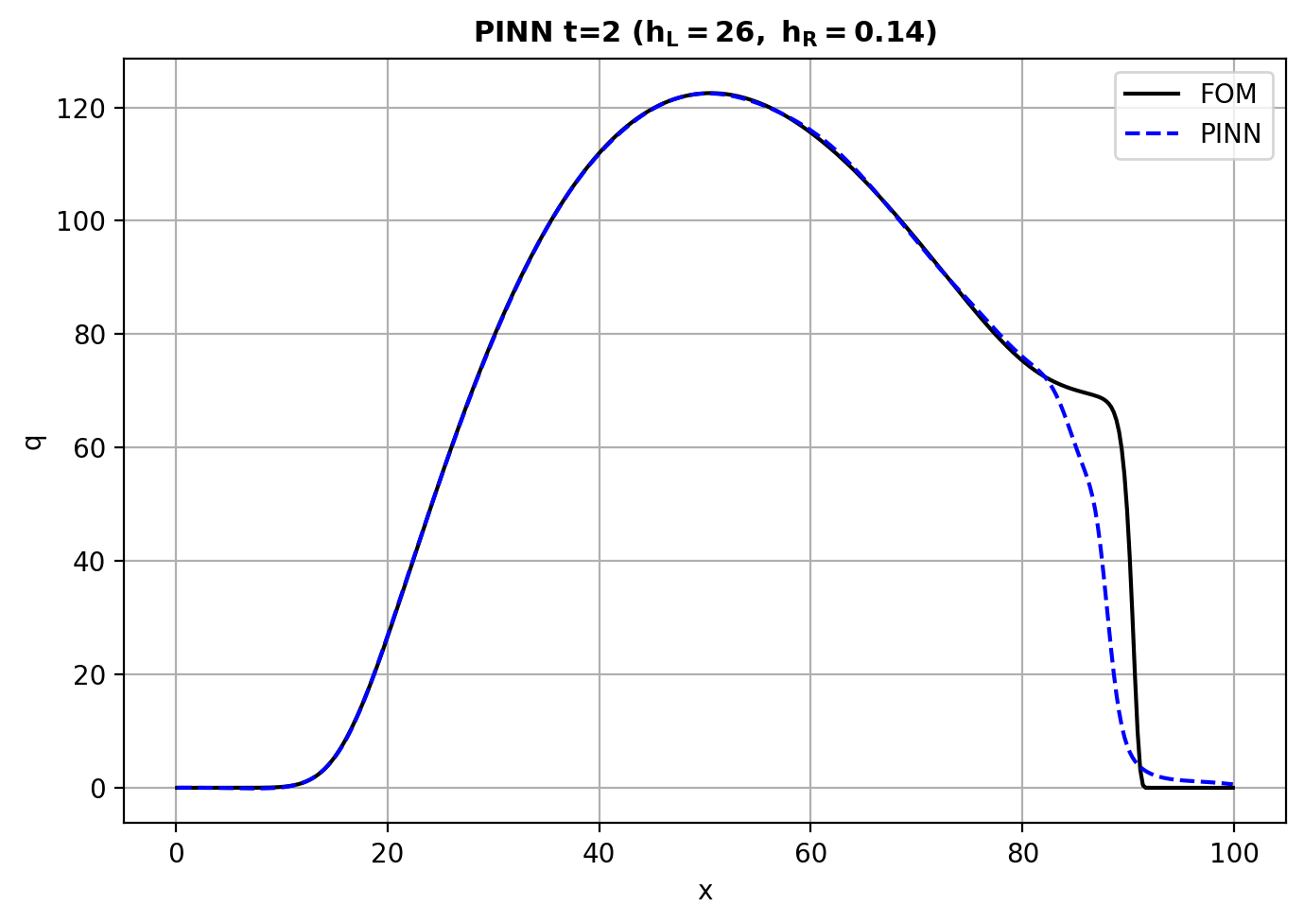}
\includegraphics[width=\figwidth\textwidth,height=\figheight]{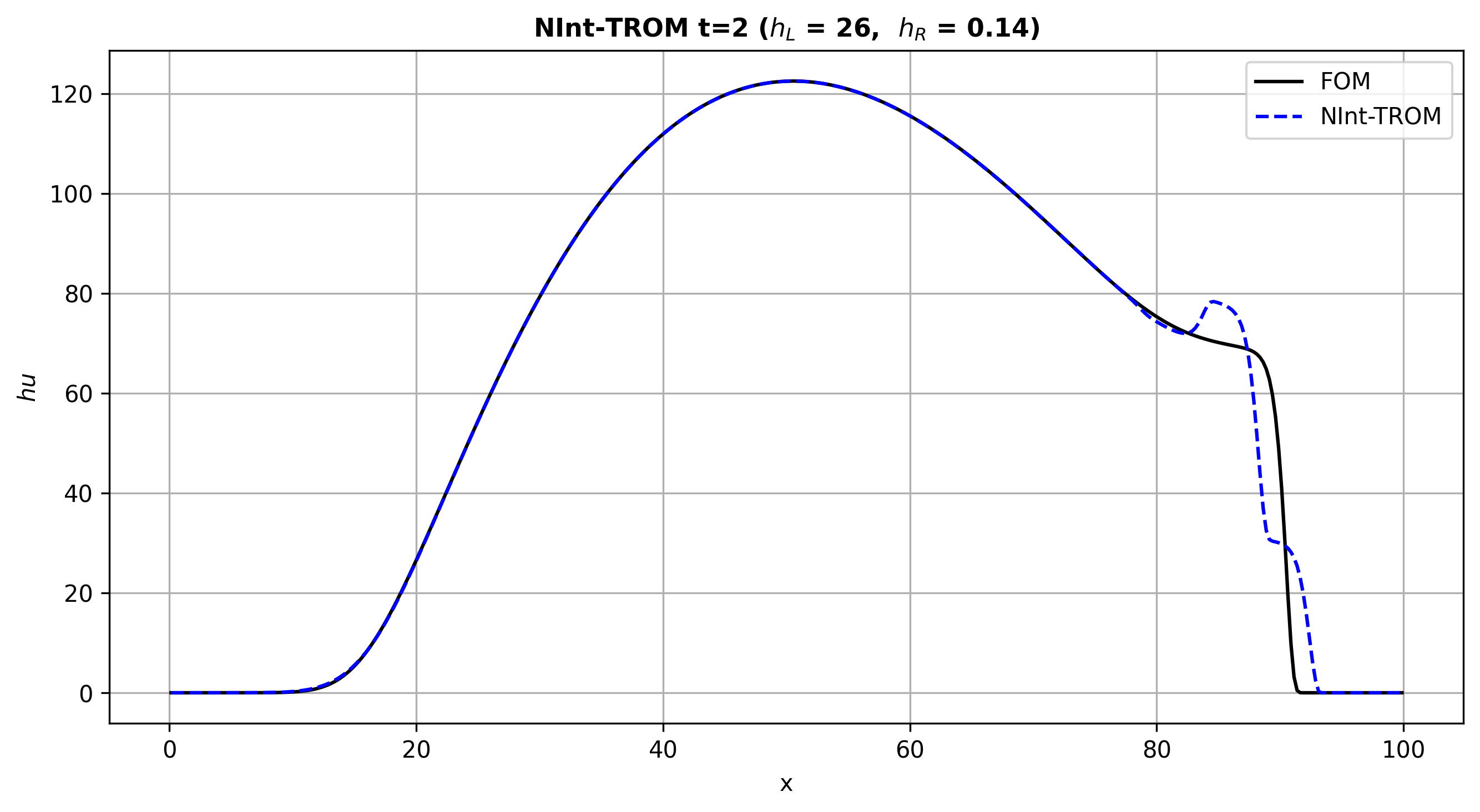}}
\centerline{\includegraphics[width=\figwidth\textwidth,height=\figheight]{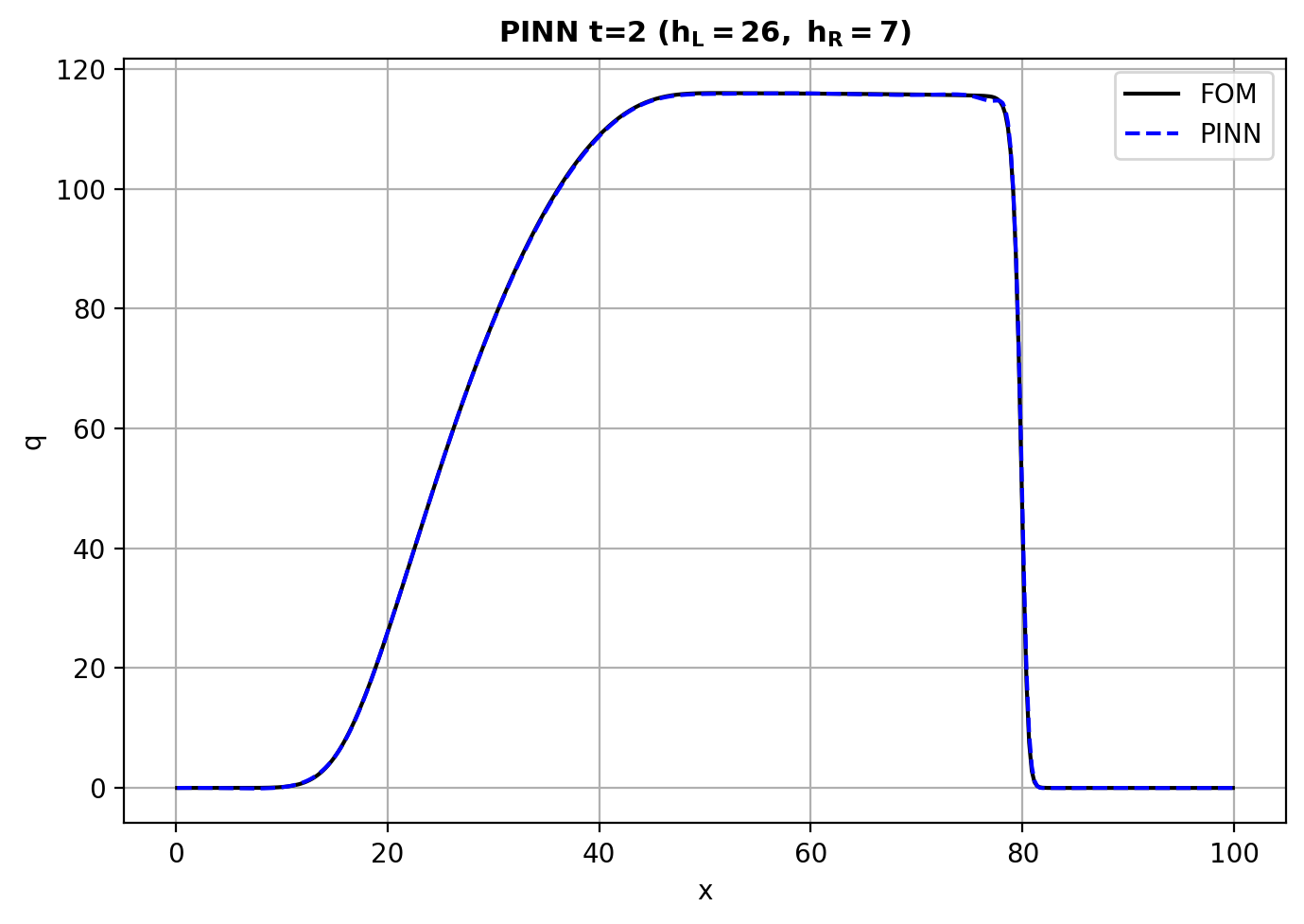}
\includegraphics[width=\figwidth\textwidth,height=\figheight]{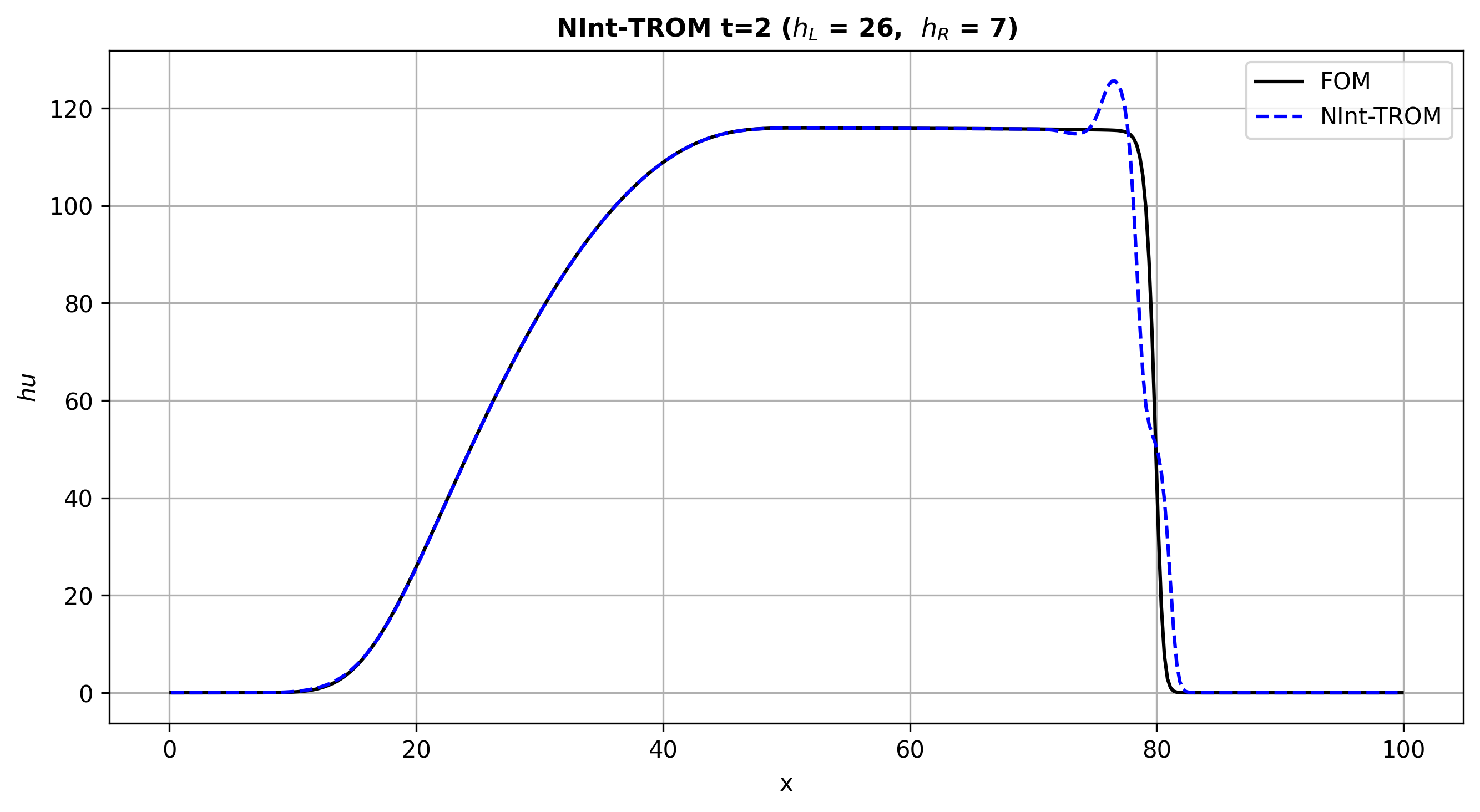}}
\caption{Comparison of $q(x,t;\muvec)$ from the FOM simulations and predictions of the PINN (left column) and the TROM (right column) for $\muvec=(18, 0), (26, 0.14), (26, 7)$ (top to bottom) at time $t=2$.}
\label{fig:profq2}
\end{figure}

Both reduced models provide accurate predictions across all six parameter regimes and reproduce the shock location with good accuracy. The largest discrepancies usually occur near the moving discontinuities, where small errors in both profiles can produce noticeable pointwise differences. Nevertheless, the relative \(L^2\) errors remain low overall, consistent with the close agreement between the predicted and the FOM profiles of \(h(x,t)\) and \(q(x,t)\) at \(t=2\) shown in Figures~\ref{fig:profh1}--\ref{fig:profq2}. The most challenging case for both models is the near-dry-bed regime \(\muvec=(26,0.14)\), for which both methods exhibit larger errors, particularly in the vicinity of the shock. A clear regime-dependent trend is also observed: the TROM is generally more accurate for the dry-bed cases, whereas the parametric PINN consistently achieves lower errors in the wet-bed regimes.

\begin{figure}[H]  
\centerline{\includegraphics[width=\figwidth\textwidth,height=\figheight]{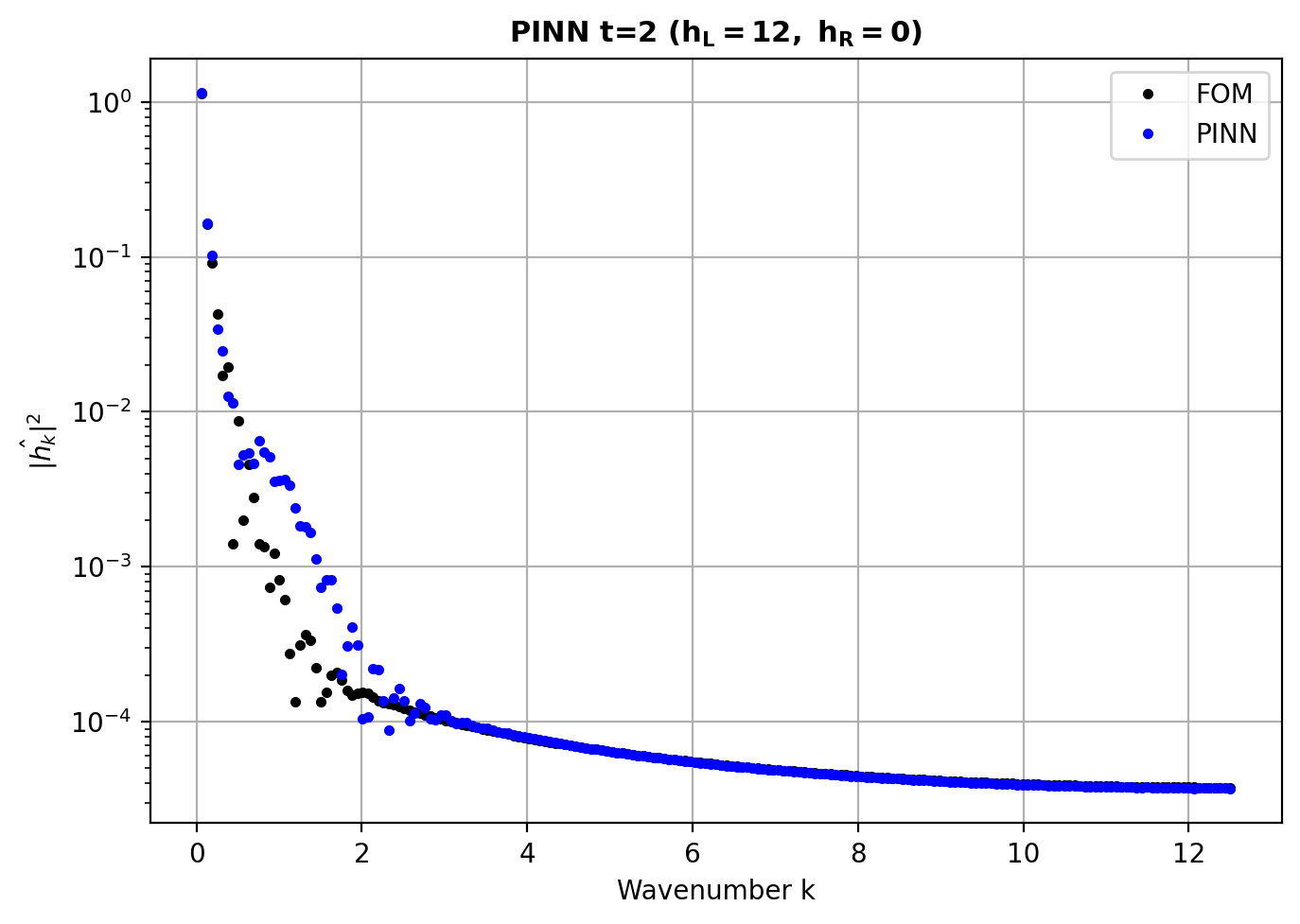}
\includegraphics[width=\figwidth\textwidth,height=\figheight]{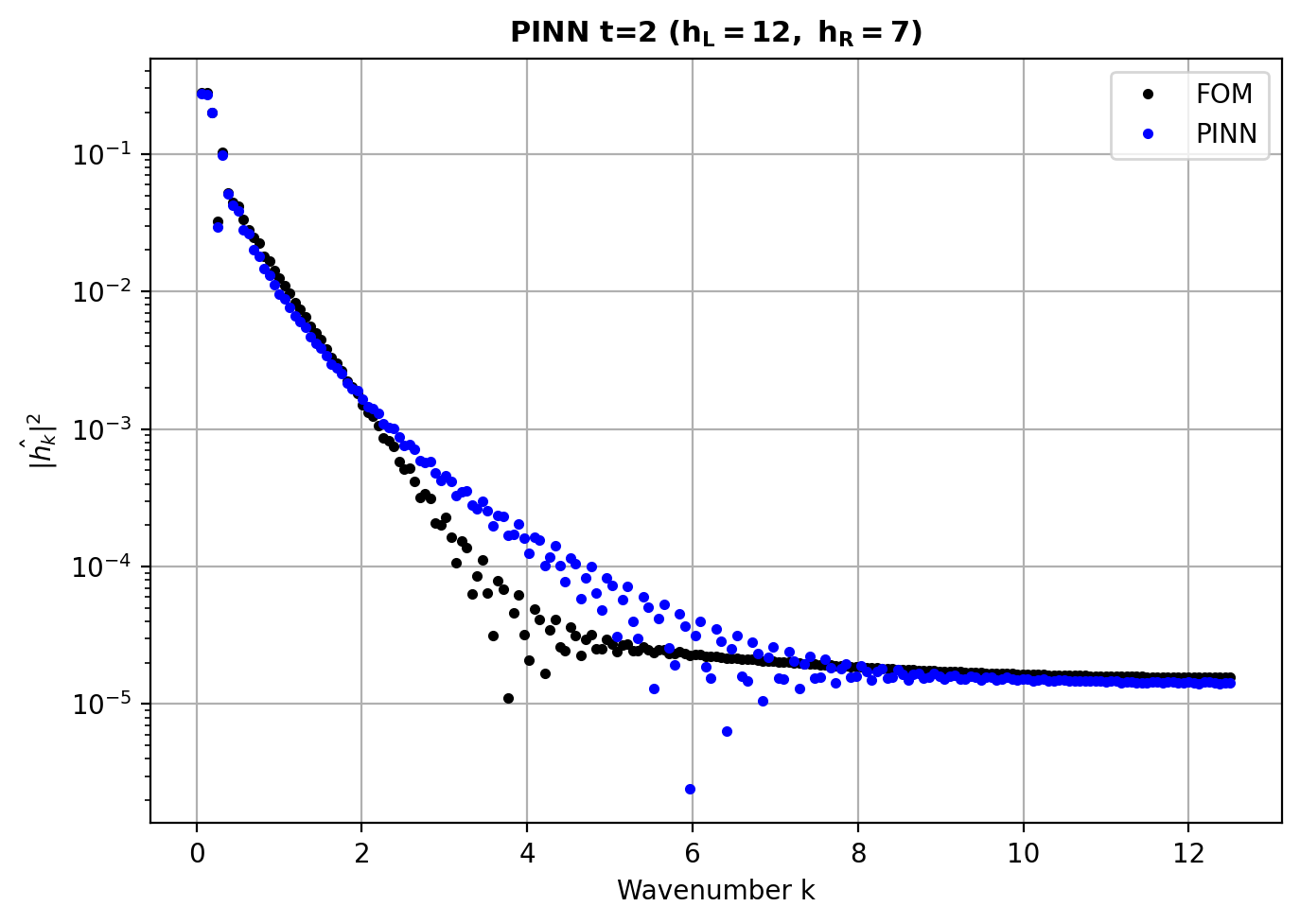}}
\centerline{\includegraphics[width=\figwidth\textwidth,height=\figheight]{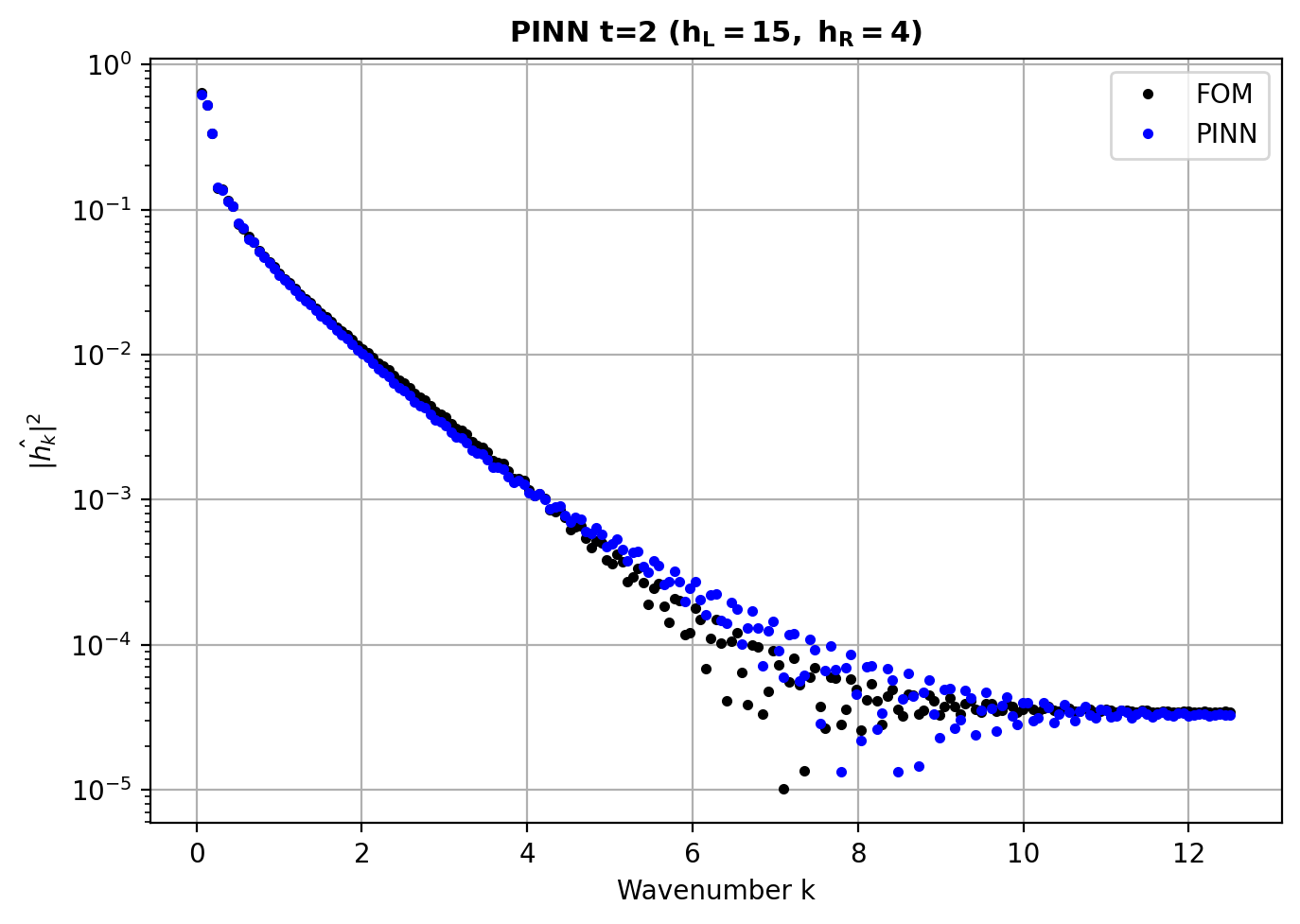}
\includegraphics[width=\figwidth\textwidth,height=\figheight]{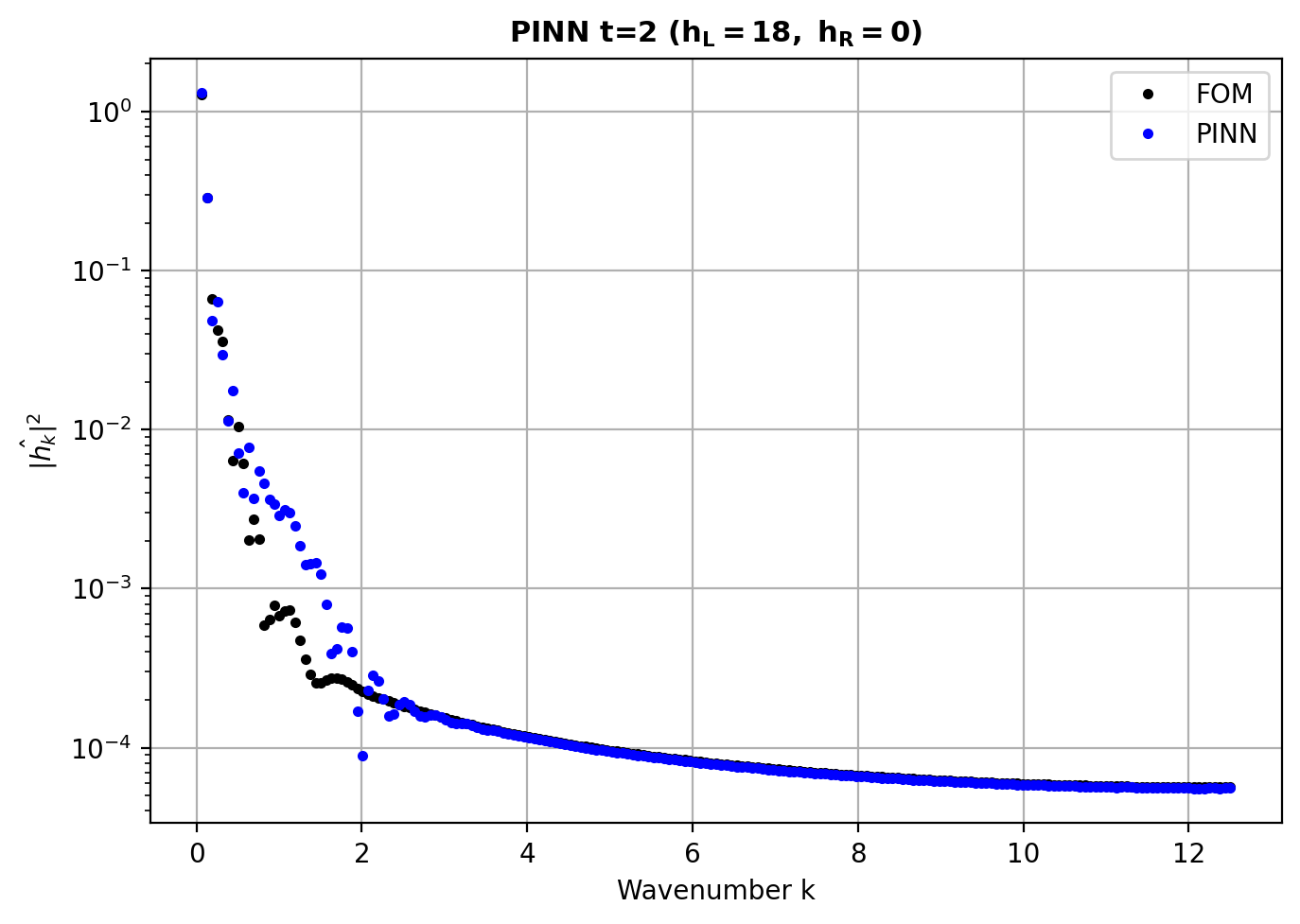}}
\centerline{\includegraphics[width=\figwidth\textwidth,height=\figheight]{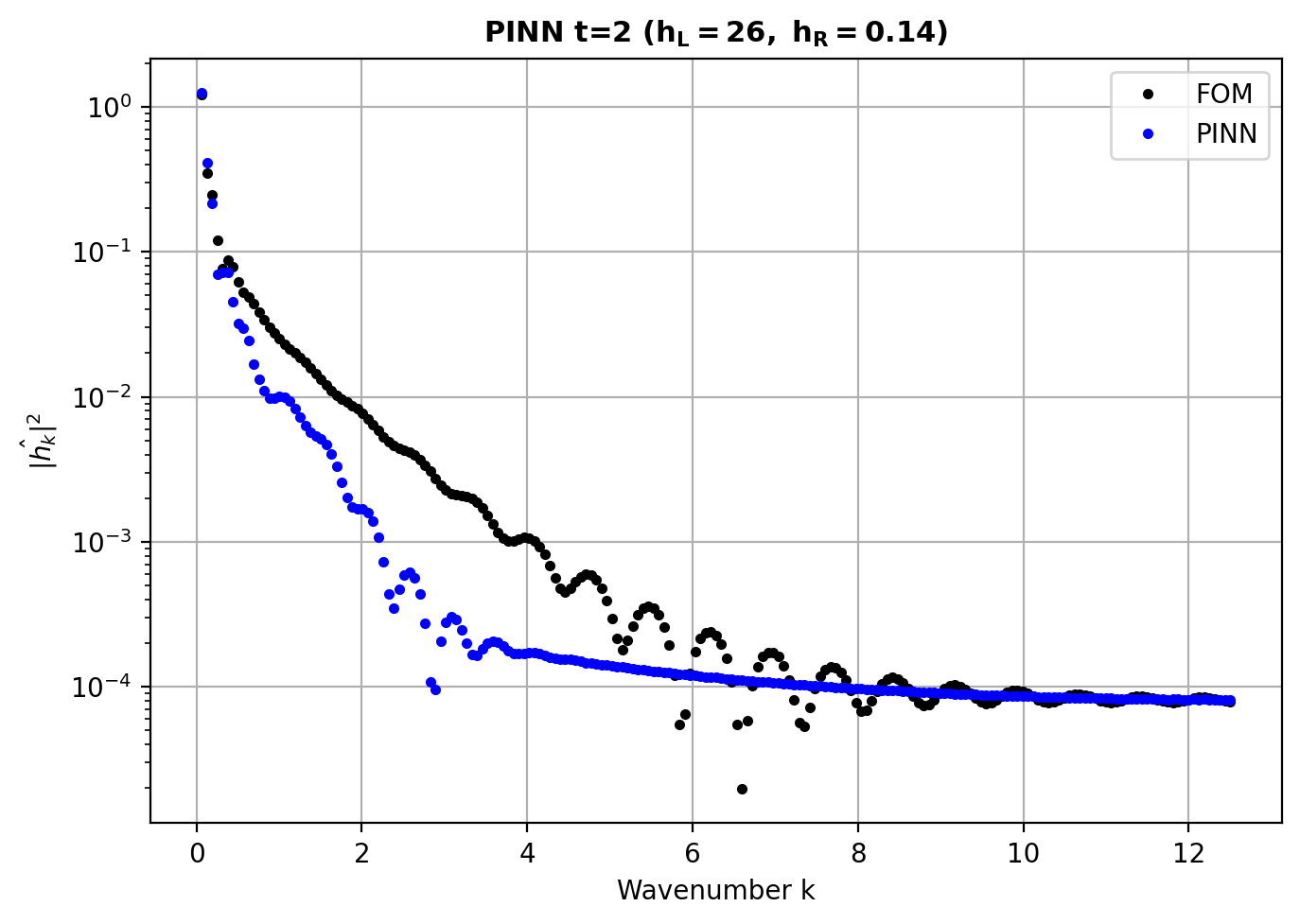}
\includegraphics[width=\figwidth\textwidth,height=\figheight]{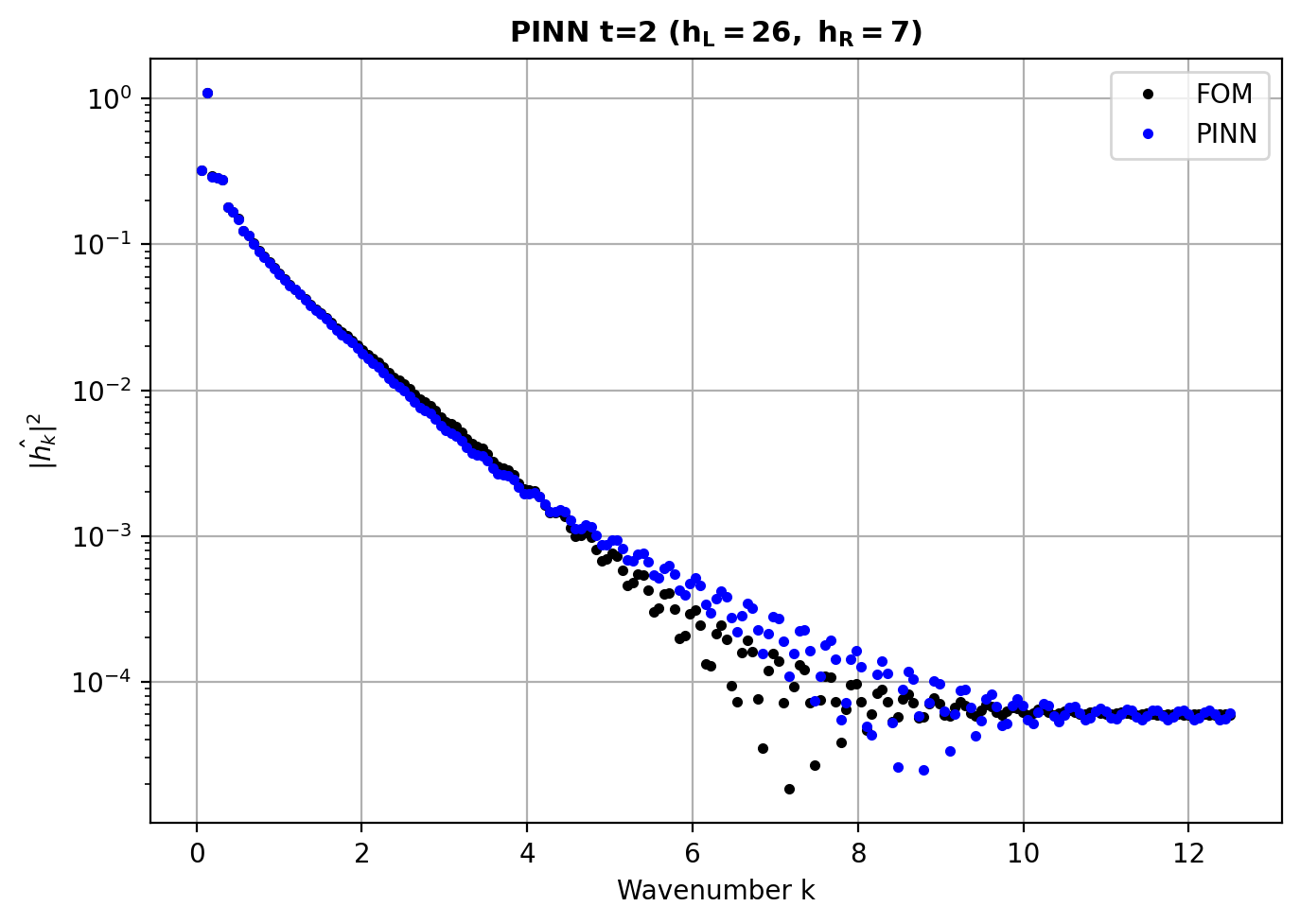}}
\caption{Fourier-amplitude spectra of the water-depth profiles at $t=2$ for $\muvec=(12,0)$ (top left), (12,7) (top right), 
(15,4) (middle left), (18,0) (middle right), (26,0.14) (bottom left), (26,7) (bottom right).}
\label{fig:spech}
\end{figure}

A well-known limitation of deep neural networks is their spectral bias (or Frequency Principle) \cite{Rahaman2019Spectral,Xu2020FPrinciple,Xu2022Overview}, when the low-frequency (large-scale) components of a target function (the solution map in our case) are learned earlier and more rapidly than high-frequency components, making sharp gradients and discontinuities more difficult to approximate accurately.
To further assess the performance of the parametric PINN model, we depict the magnitude of the spectra of the predicted water depth \(h(x,t=2)\) in Figure \ref{fig:spech}. To this end, we first remove the linear trend connecting the boundary values, \((h(0,t),h(L,t))=(h_L,h_R)\), and then compute the magnitude of the fast Fourier transform.
The largest spectral discrepancies occur at high wavenumbers, particularly for the dry-bed and near-dry-bed regimes, which contain the sharpest solution features. Among the tested cases, \(\muvec=(26,0.14)\) presents the greatest challenge for reproducing high-frequencies of the solution; here the spectra in the predicted solution decays faster that the spectra in the reference FOM solution. 
\rev{For all other cases, the predicted spectra show good agreement with the reference FOM over a broad range of wavenumbers, while the largest discrepancies remain concentrated at high wavenumbers. Thus, the spectral agreement is stronger for larger values of \(h_R\), whereas the dry-bed and near-dry-bed regimes remain more challenging.}
%

Tables \ref{tab:l2errt1} and \ref{tab:l2errt2} indicate that the relative errors for both the water depth and the discharge remain of the same order in most cases, although some increase with time. In addition, the relative errors for the discharge are slightly larger and more regime-dependent than those for the water depth, reflecting the greater sensitivity of $q=hu$ to local mismatches near the shock front. Similar to the water-depth errors, the dry-bed and near-dry regimes remain the most challenging, while the wet-bed cases achieve the lowest discharge errors.
Overall, these results demonstrate that the proposed parametric PINN is capable of learning the \rev{FOM parameter-to-solution map} with high accuracy across a diverse range of dam-break configurations. Despite the challenges posed by hyperbolic dynamics and moving discontinuities, it achieves low relative errors and accurate shock localization, making it a robust and competitive surrogate model for parametrized shallow-water flows.

\begin{figure}[H]  
\centerline{\includegraphics[width=\figwidth\textwidth,height=\figheight]{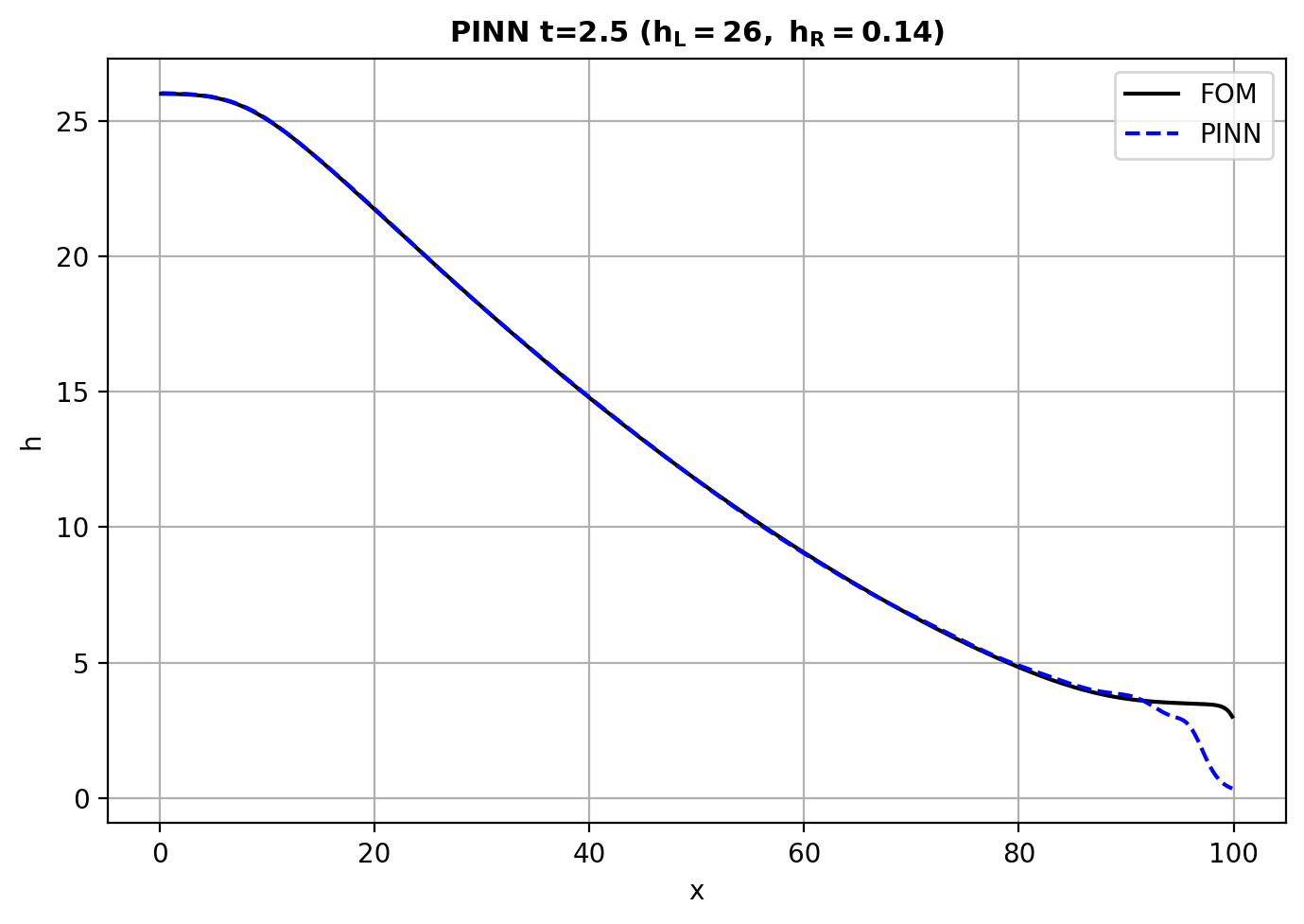}
\includegraphics[width=\figwidth\textwidth,height=\figheight]{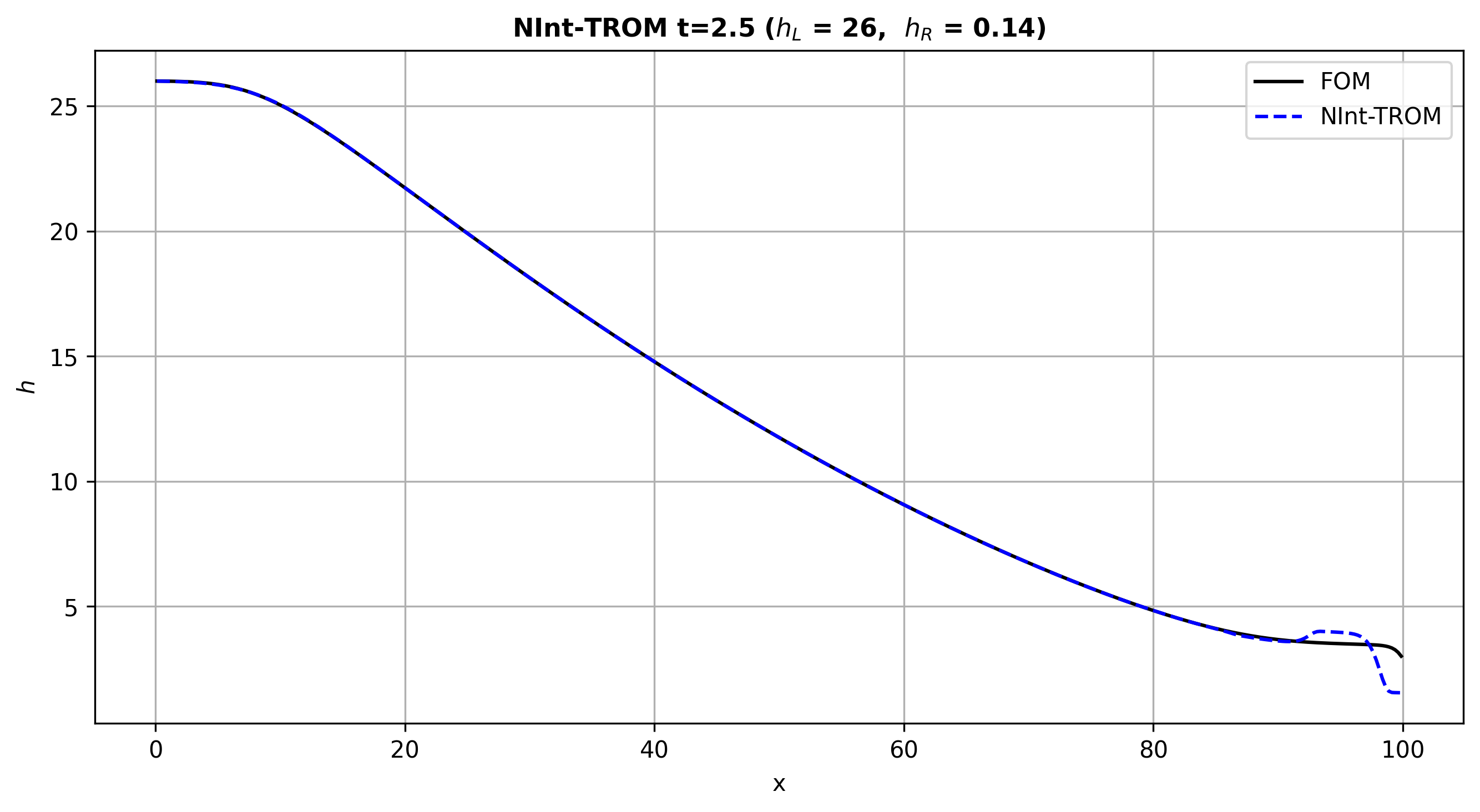}}
\centerline{\includegraphics[width=\figwidth\textwidth,height=\figheight]{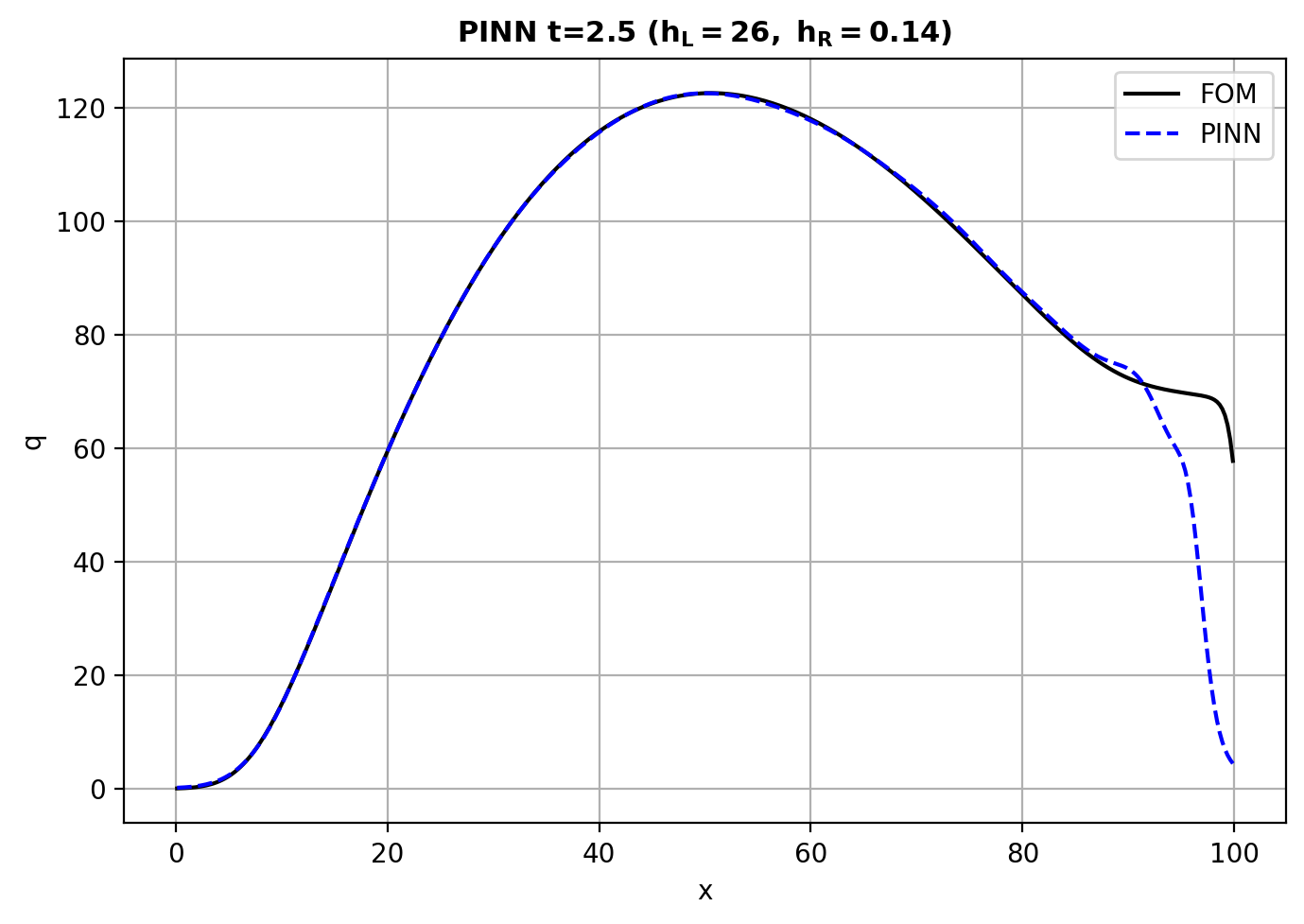}
\includegraphics[width=\figwidth\textwidth,height=\figheight]{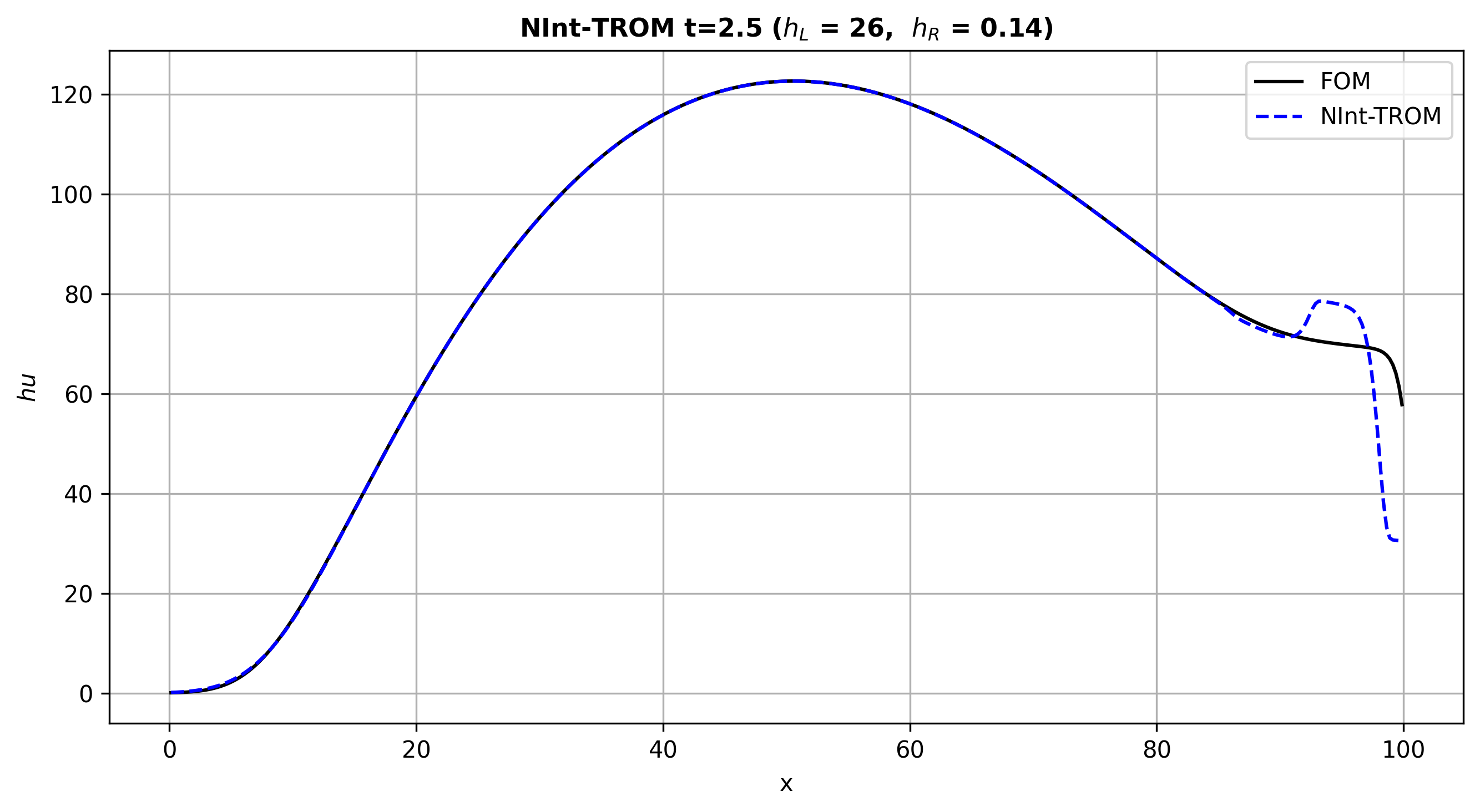}}
\caption{Comparison of the water-depth profiles $h(x,t)$ (top) and discharge profiles $q(x,t)$ (bottom) at $t=2.5$ predicted by the PINN (left column) and TROM (right column) for $\muvec=(26,0.14)$.}
\label{fig:t2.5}
\end{figure}

For the dam-break problem, a right-propagating shock forms when $h_R>0$, whereas the dry-bed case is characterized by a moving rarefaction wave. Since the shock propagation speed increases as $h_R$ decreases, the shock exits the computational domain by $t=2.5$ for the near-dry-bed regime $\muvec=(26,0.14)$. 
Figure~\ref{fig:t2.5} shows the predictions of the PINN and the TROM for the near-dry-bed regime $\muvec=(26,0.14)$ at time $t=2.5$. Although the training data cover the time interval up to $t=2.5$, neither model accurately captures the FOM solution at time $t=2.5$ near the right boundary in this regime. However, the discrepancies are localized near the right boundary, and both models accurately predict the solution over the rest of the computational domain.

Our results indicate that both models reproduce the FOM simulations with high accuracy for parameter regimes within the training domain. The errors of both the PINN and TROM are primarily localized near moving fronts, \rev{while the global wave structure is preserved. 
For the PINN, the spectral comparisons also show good agreement with the FOM over low and intermediate wavenumbers in most cases.}
The dry-bed and near-dry-bed configurations are the most challenging for both models, with the discharge being more sensitive than the water depth. In addition, both models face difficulties when predicting FOM solutions at times when the shock exits the computational domain.

\subsection{Computational Time Comparison}
\rev{Generation of the FOM training dataset over the 65 training parameter instances required $42.19\pm4.3$~s, with individual FOM run times ranging from $0.64$ to $0.66$~s.}

\rev{
For the PINN, shock detection and construction of the shock trajectories for the full training dataset required $0.049\pm0.018$~s. Thus, the shock-aware collocation strategy introduces negligible preprocessing cost and does not require any additional FOM simulations beyond those already used to construct the supervised training dataset.

Training the final PINN configuration for $20{,}000$ epochs required $18.5\pm3.0$~h, including the validation evaluations performed during training. These timing results are reported as the mean and sample standard deviation over five independent training runs. Once trained, evaluation of the PINN for a single parameter–time query over the full spatial grid requires approximately $1.6\times10^{-3}$~s. Hence, the dominant computational expense of the PINN is the offline neural-network training, whereas FOM dataset generation, shock-aware preprocessing, and subsequent online evaluations constitute comparatively small costs.
}

\rev{
Construction of the Tucker decomposition of the resulting snapshot tensors required $2.5$~s ($0.7$~s for the $h$-component, $1.8$~s for the $hu$-component). Thus, FOM data generation accounts for approximately 94\% of the offline cost, with the Tucker decomposition contributing the remaining 6\%.

Once constructed, the TROM required approximately $6.3\times10^{-4}$~s to reconstruct a solution for a single online query, compared to approximately $0.65$~s for a corresponding FOM run, a speedup of roughly $1049$ times. Therefore, the dominant computational cost of the TROM lies in the offline construction stage, primarily FOM data generation rather than the decomposition step, while online evaluation is negligible by comparison.
}

\section{Extrapolation experiments}
\label{sec:extrap}

To further assess the performance of the parametric PINN and the TROM, we consider the extrapolation benchmark, where we evaluate the accuracy of both models for parameter values outside of the training range  $(h_L,h_R) \in [10,28]\times [0,8]$ for $\calP = \calP_1 \times \calP_2$ with parameters for the training mesh $\calP_{1}$, $\calP_{2}$
in \eqref{eq:paramgrid}. This test highlights how well each model generalizes beyond the training regime.
Specifically, we consider six parameter values $(h_L,h_R) = (9,4),$
$(20,8.5)$, $(29,4)$, $(32,4)$, $(29,8.5)$, and $(32,10)$ and predict the solution at time $t=2$.
Note that the value $h_R=4$ is inside the training range for $h_R$, but the corresponding values for $h_L$ are outside of $[10,28]$.

\begin{figure}[H]  
\centerline{\includegraphics[width=\figwidth\textwidth,height=\figheight]{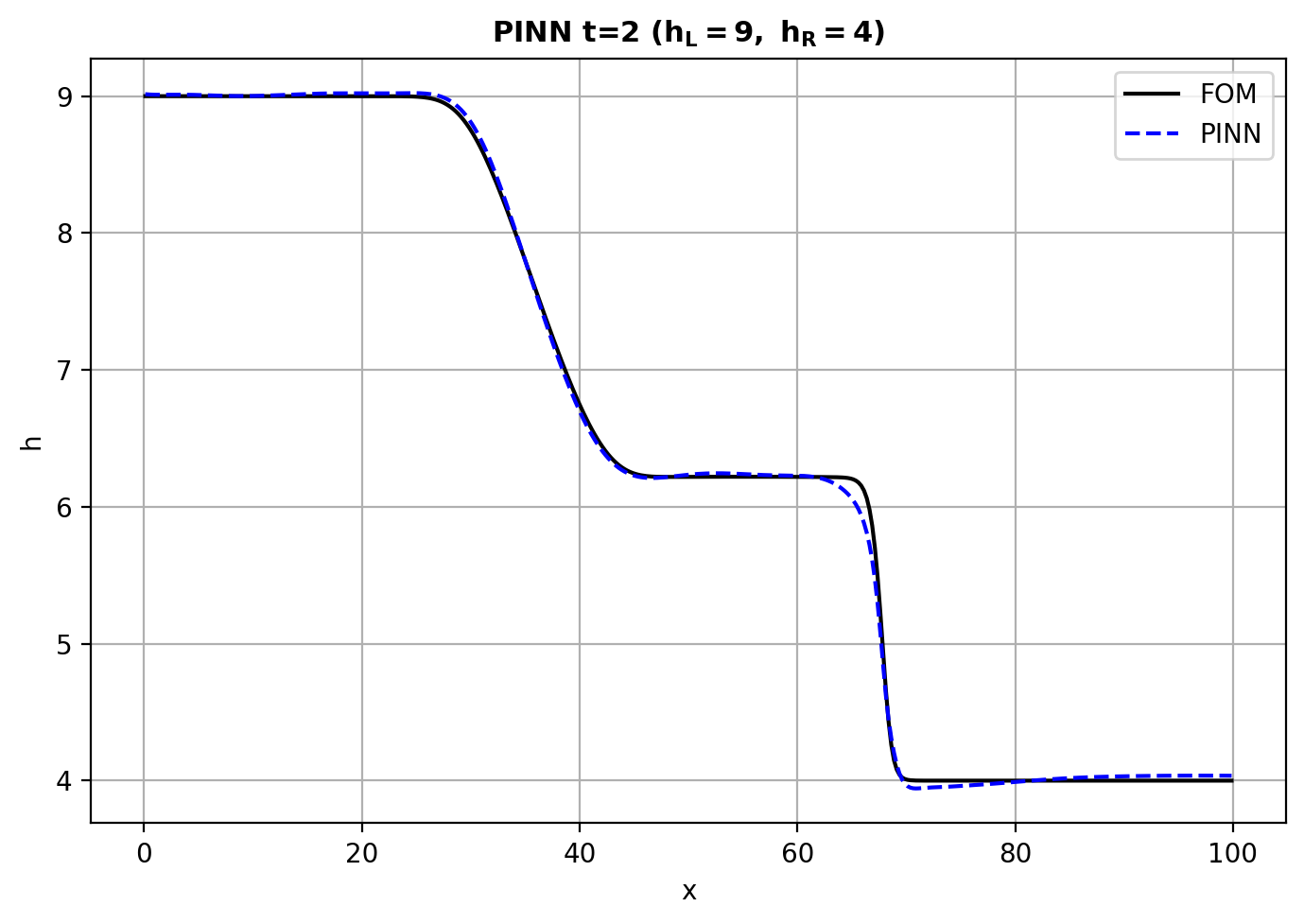}
\includegraphics[width=\figwidth\textwidth,height=\figheight]{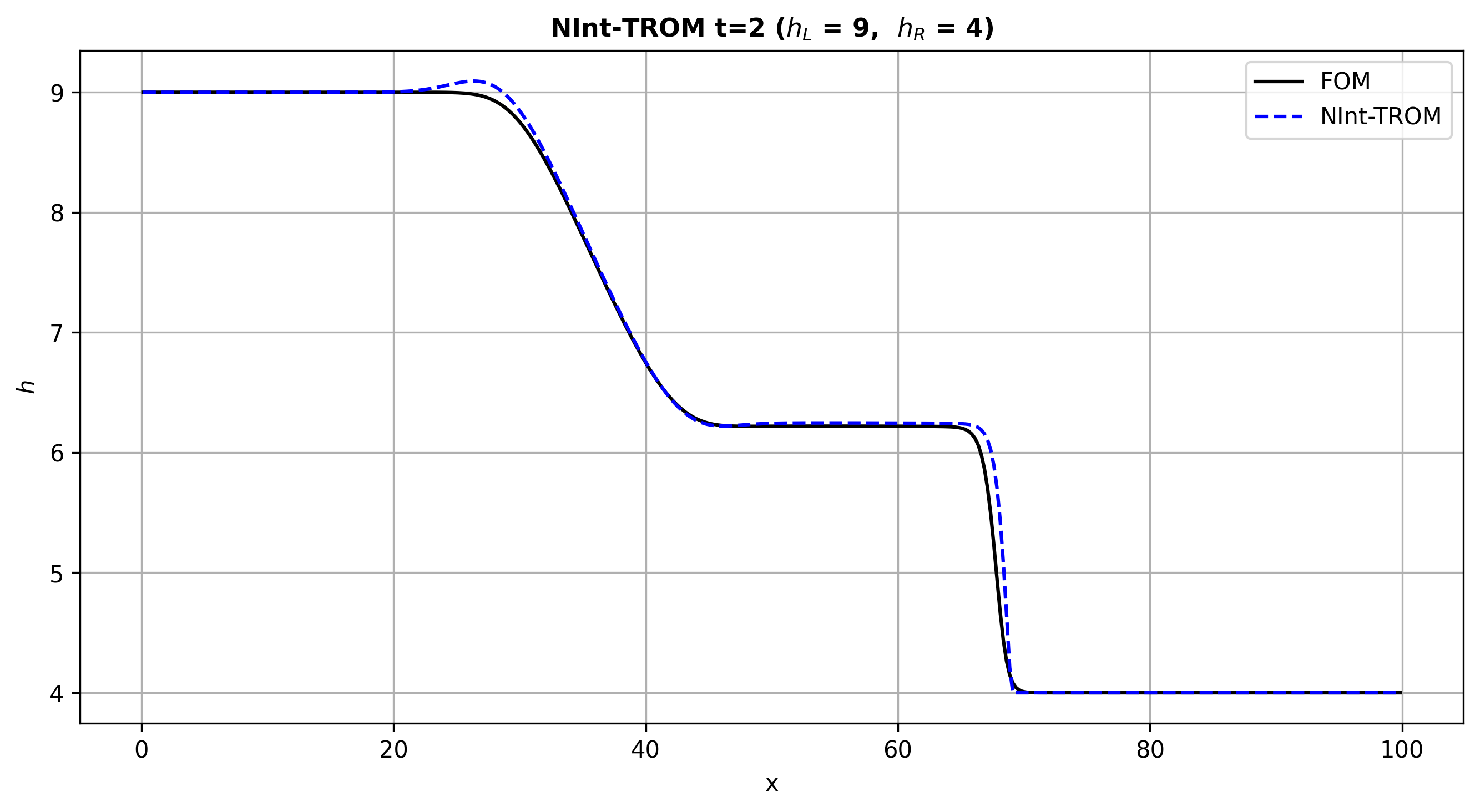}}
\centerline{\includegraphics[width=\figwidth\textwidth,height=\figheight]{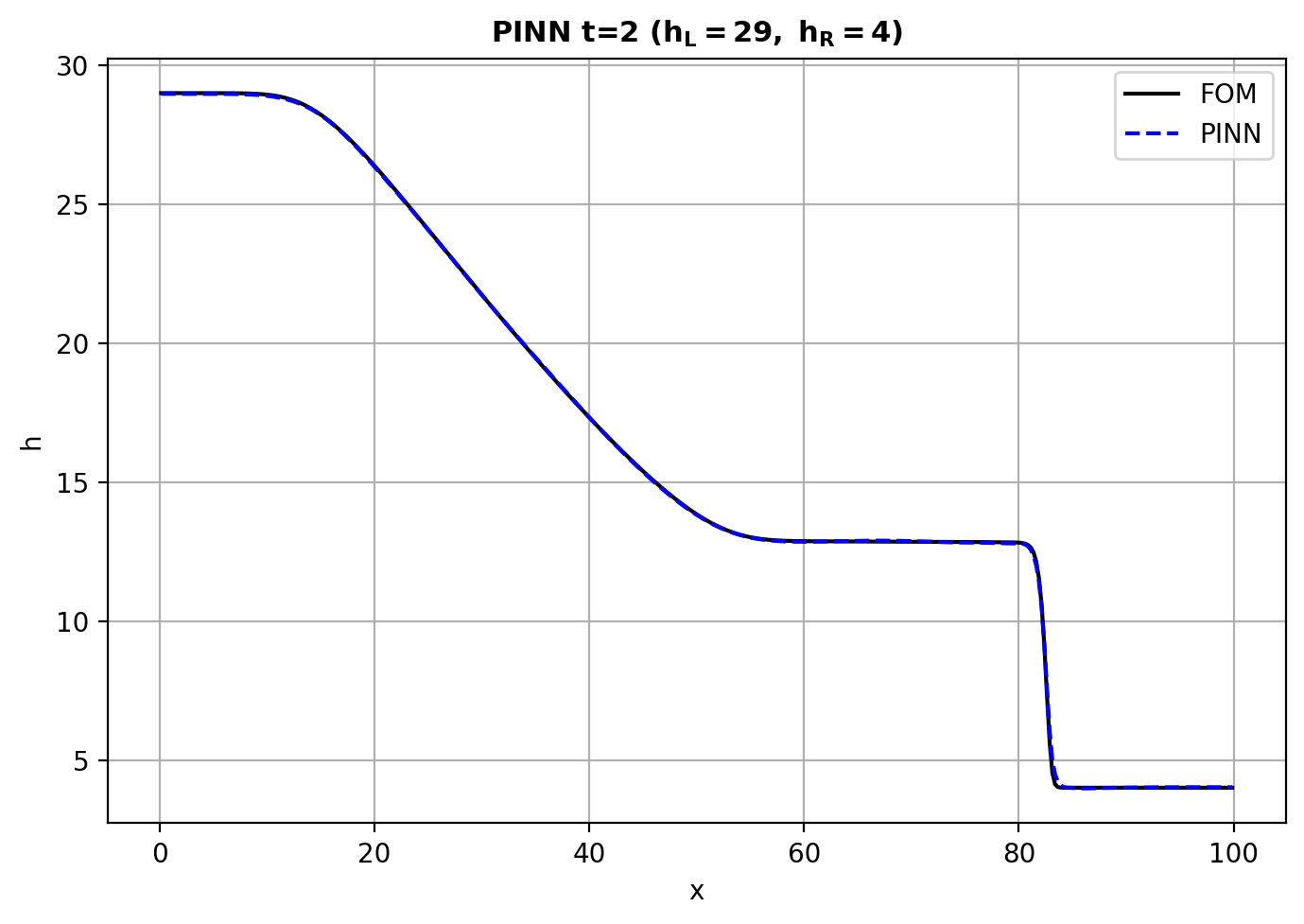}
\includegraphics[width=\figwidth\textwidth,height=\figheight]{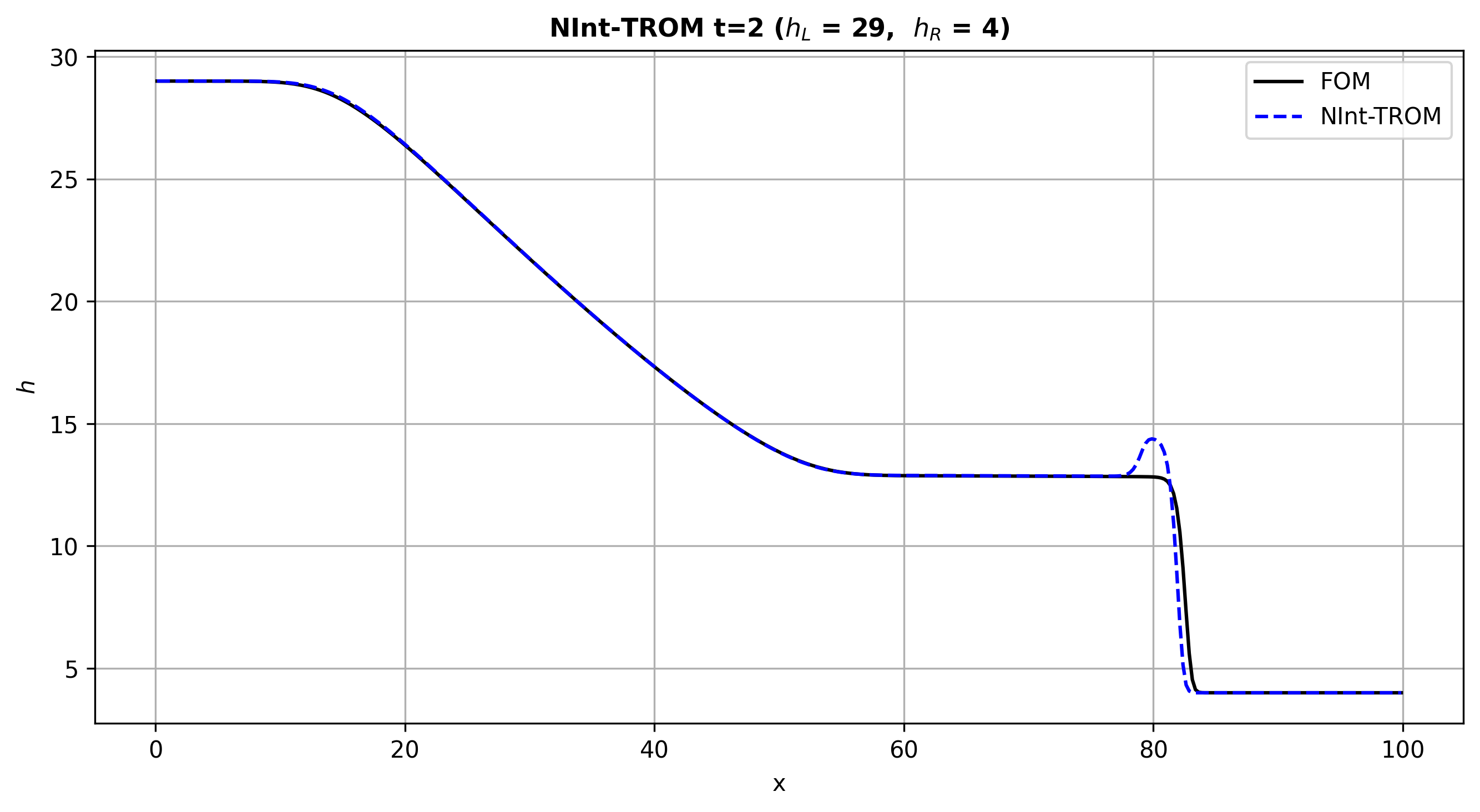}}
\centerline{\includegraphics[width=\figwidth\textwidth,height=\figheight]{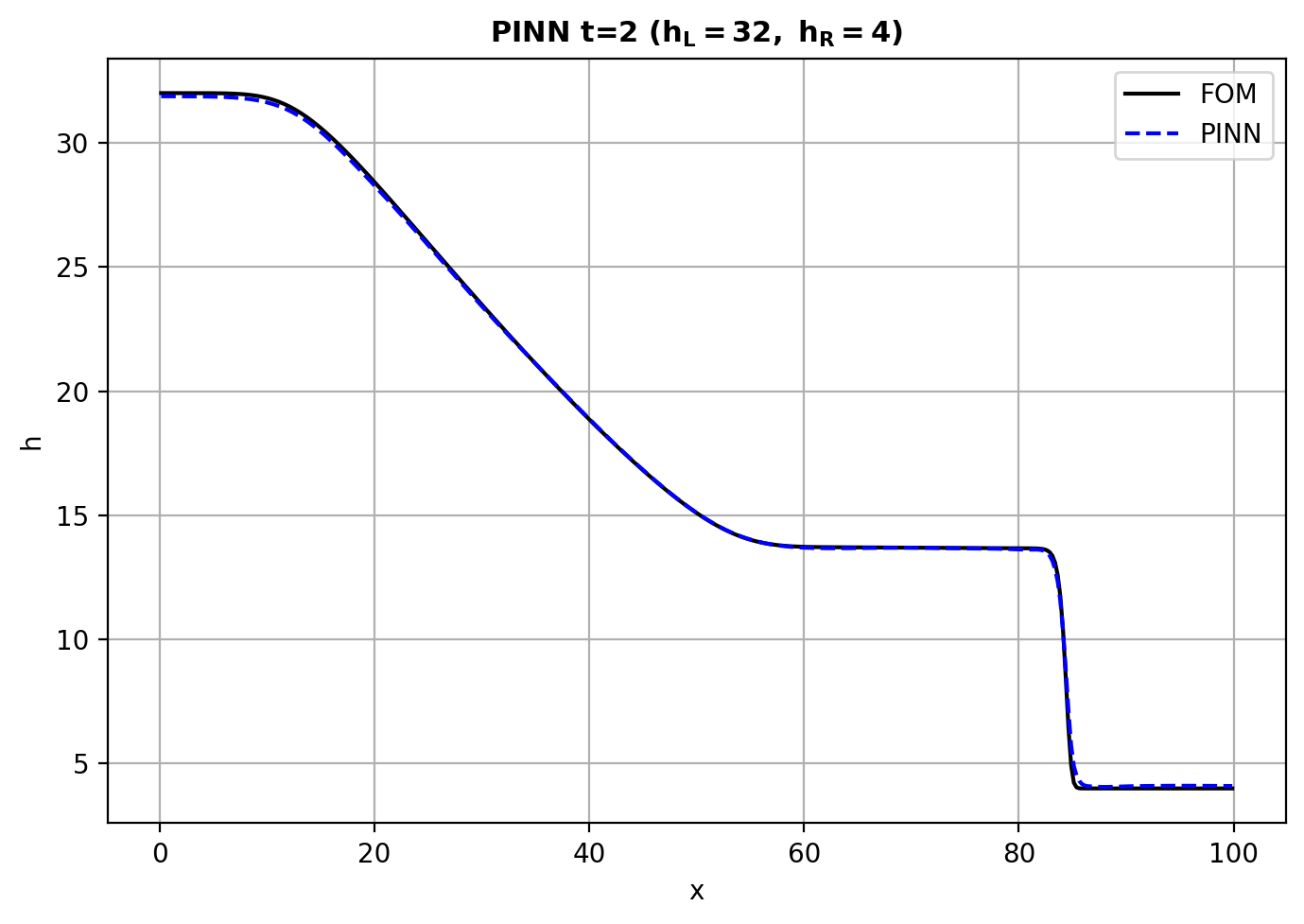}
\includegraphics[width=\figwidth\textwidth,height=\figheight]{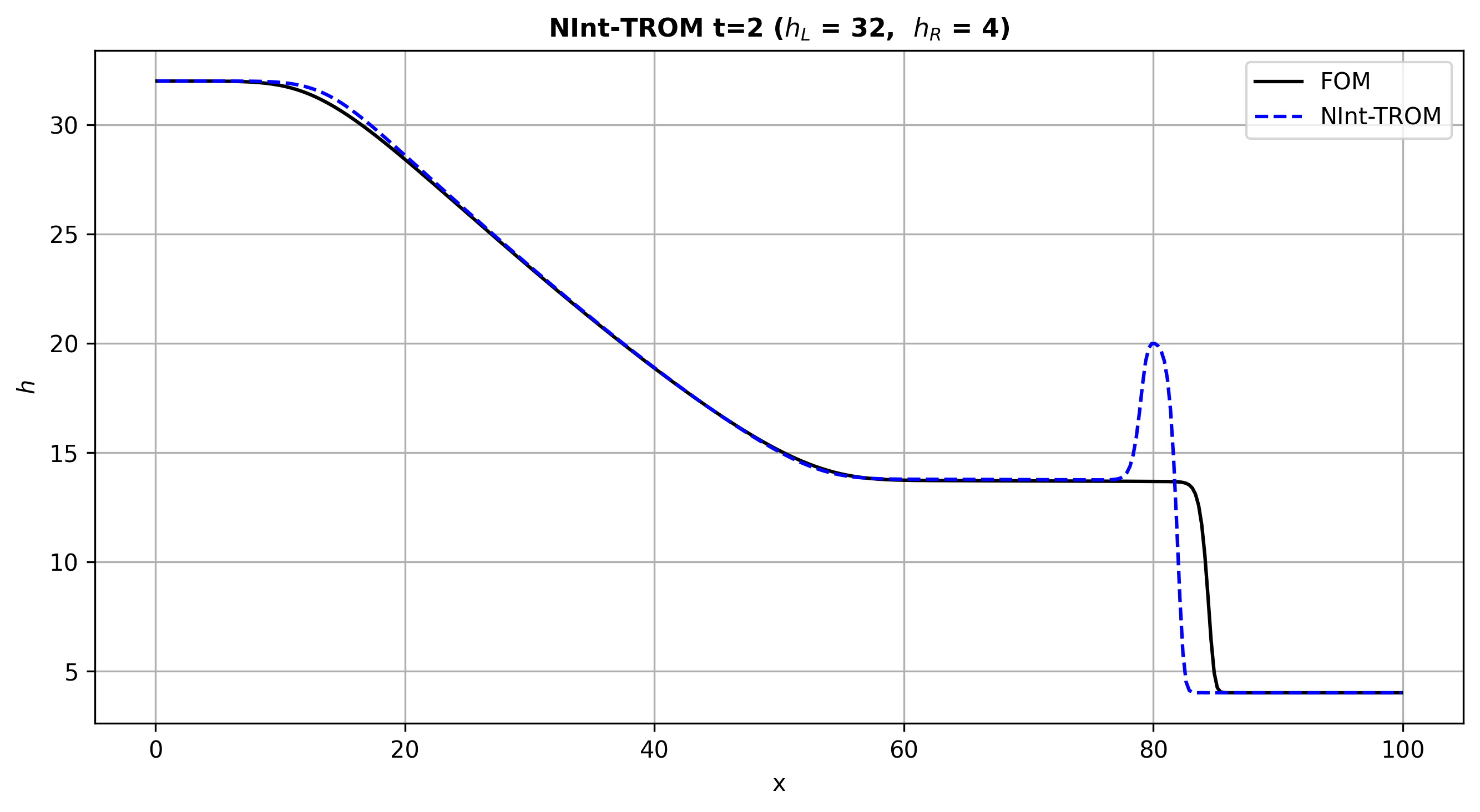}}
\caption{Comparison of $h(x,t;\muvec)$ in FOM simulations and predictions of PINN (left column) and TROM (right column) for $\muvec=(9, 4), (29, 4), (32, 4)$ (top to bottom) at time $t=2$.}
\label{fig:extraph1}
\end{figure}

\begin{figure}[H]  
\centerline{\includegraphics[width=\figwidth\textwidth,height=\figheight]{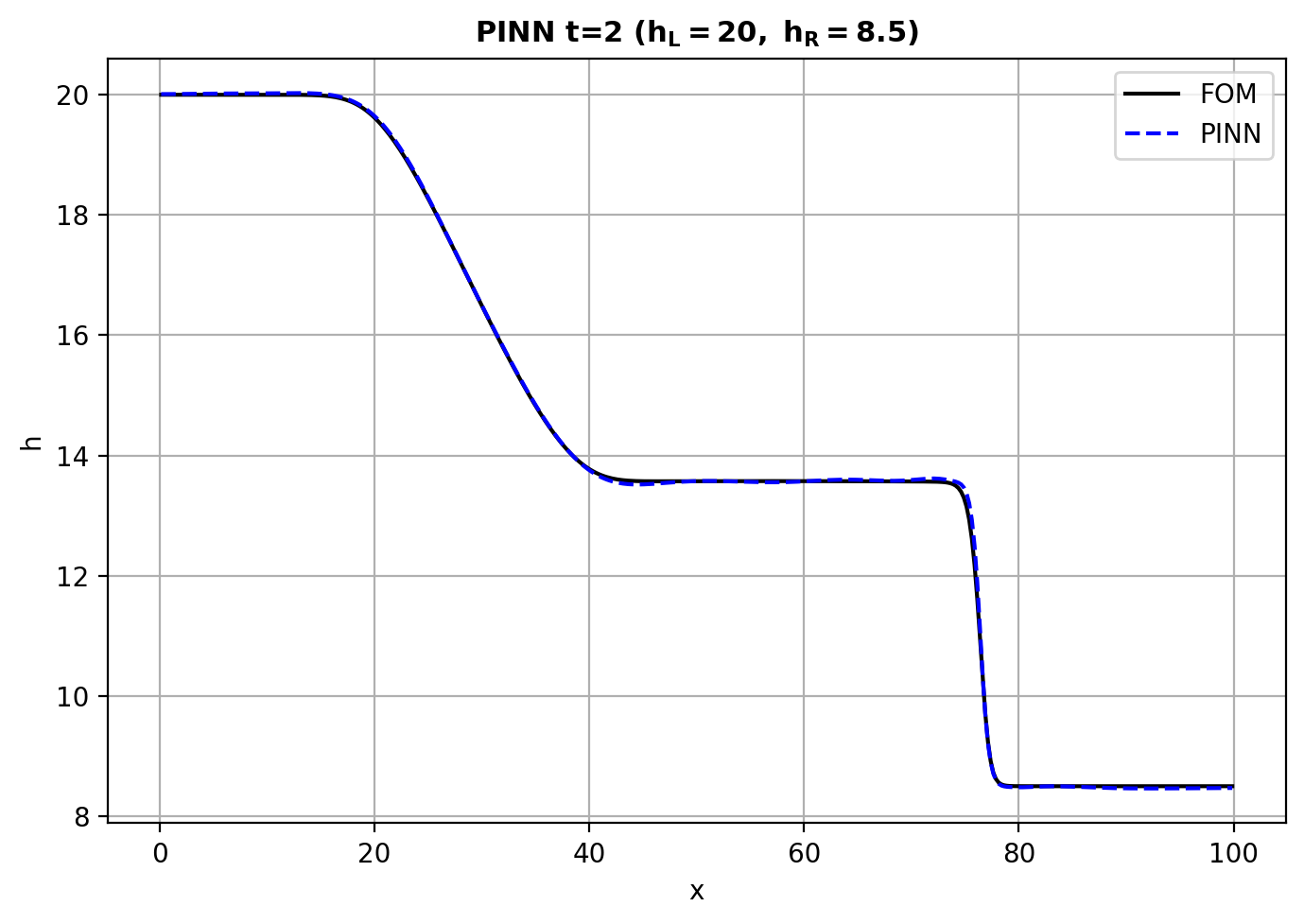}
\includegraphics[width=\figwidth\textwidth,height=\figheight]{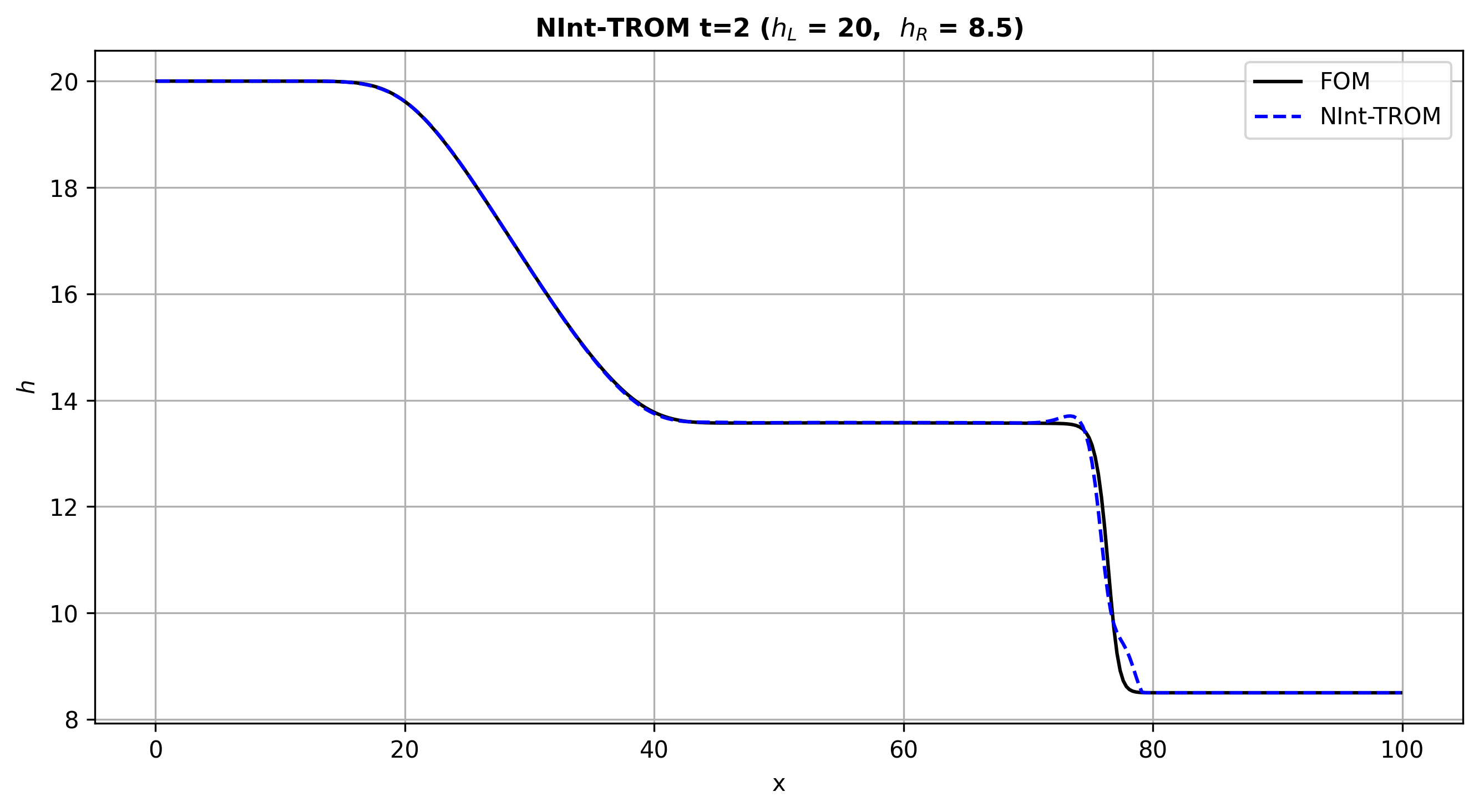}}
\centerline{\includegraphics[width=\figwidth\textwidth,height=\figheight]{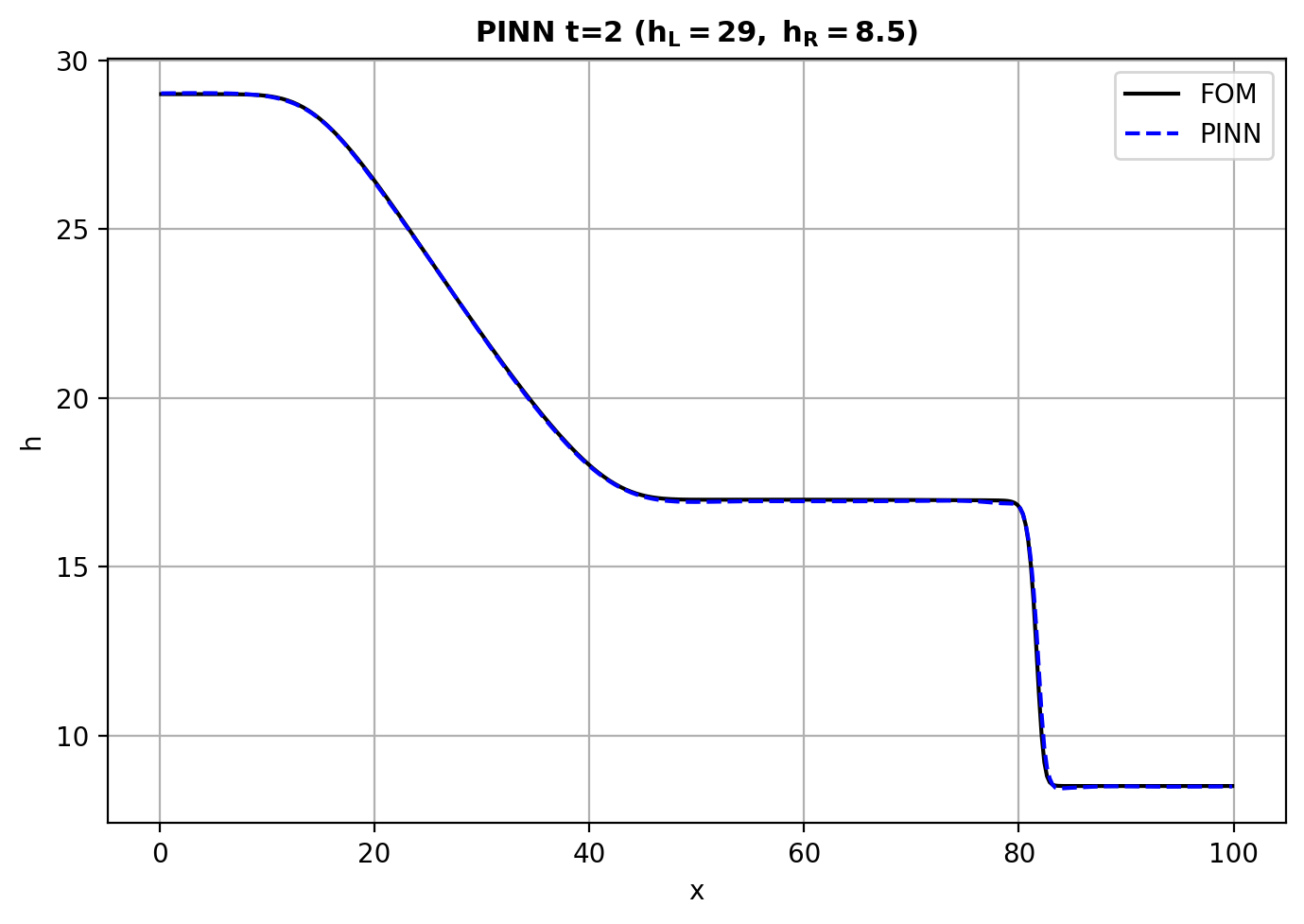}
\includegraphics[width=\figwidth\textwidth,height=\figheight]{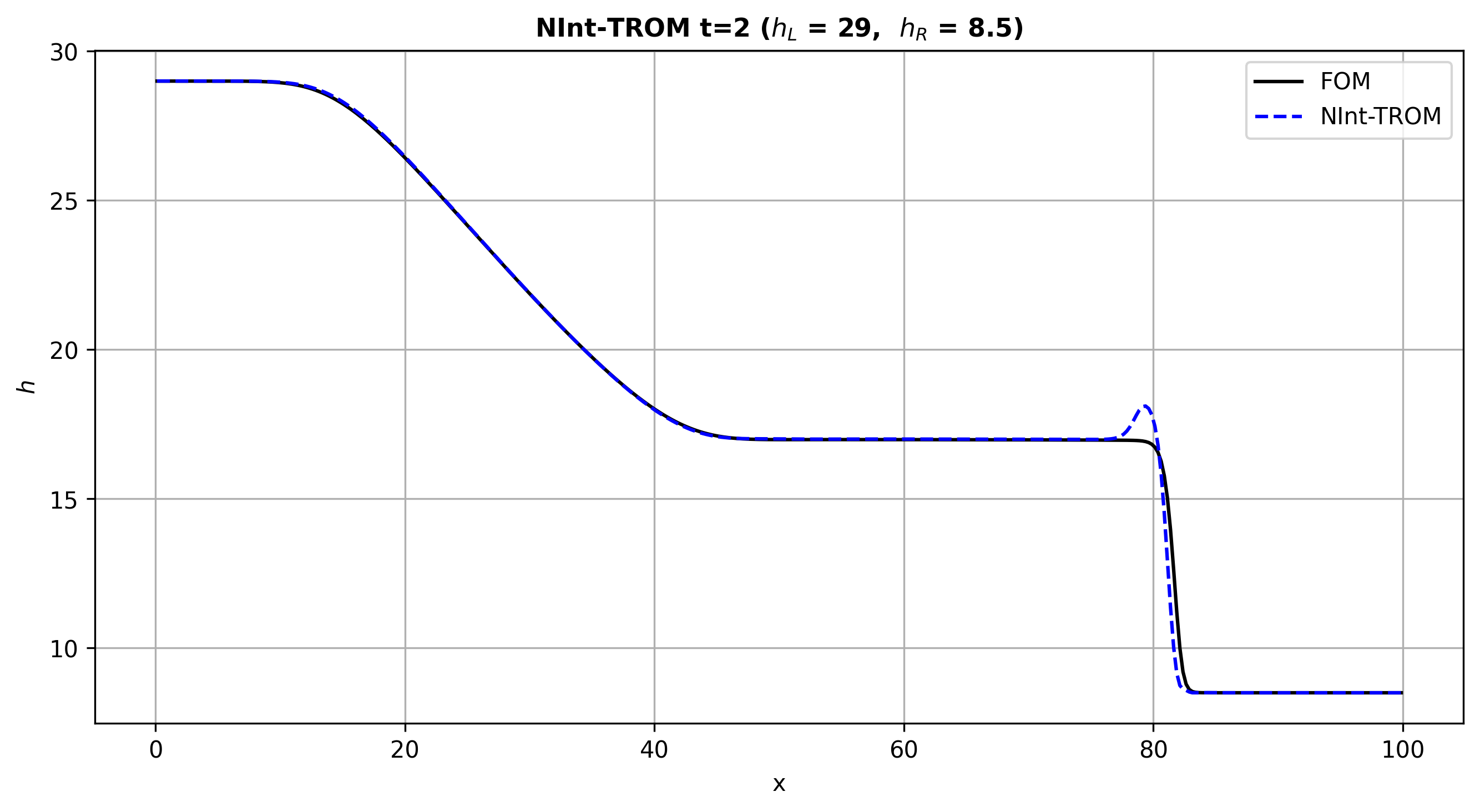}}
\centerline{\includegraphics[width=\figwidth\textwidth,height=\figheight]{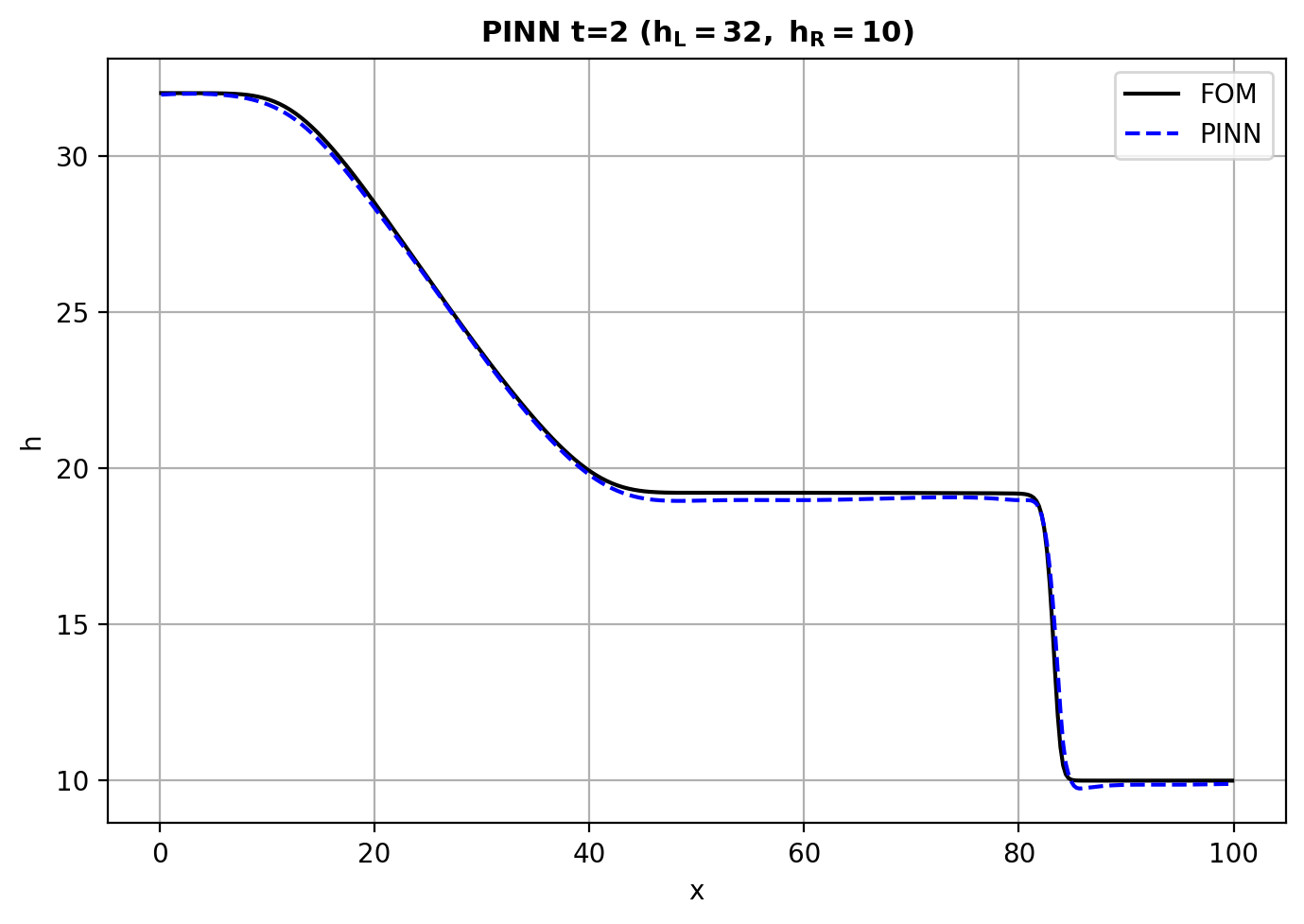}
\includegraphics[width=\figwidth\textwidth,height=\figheight]{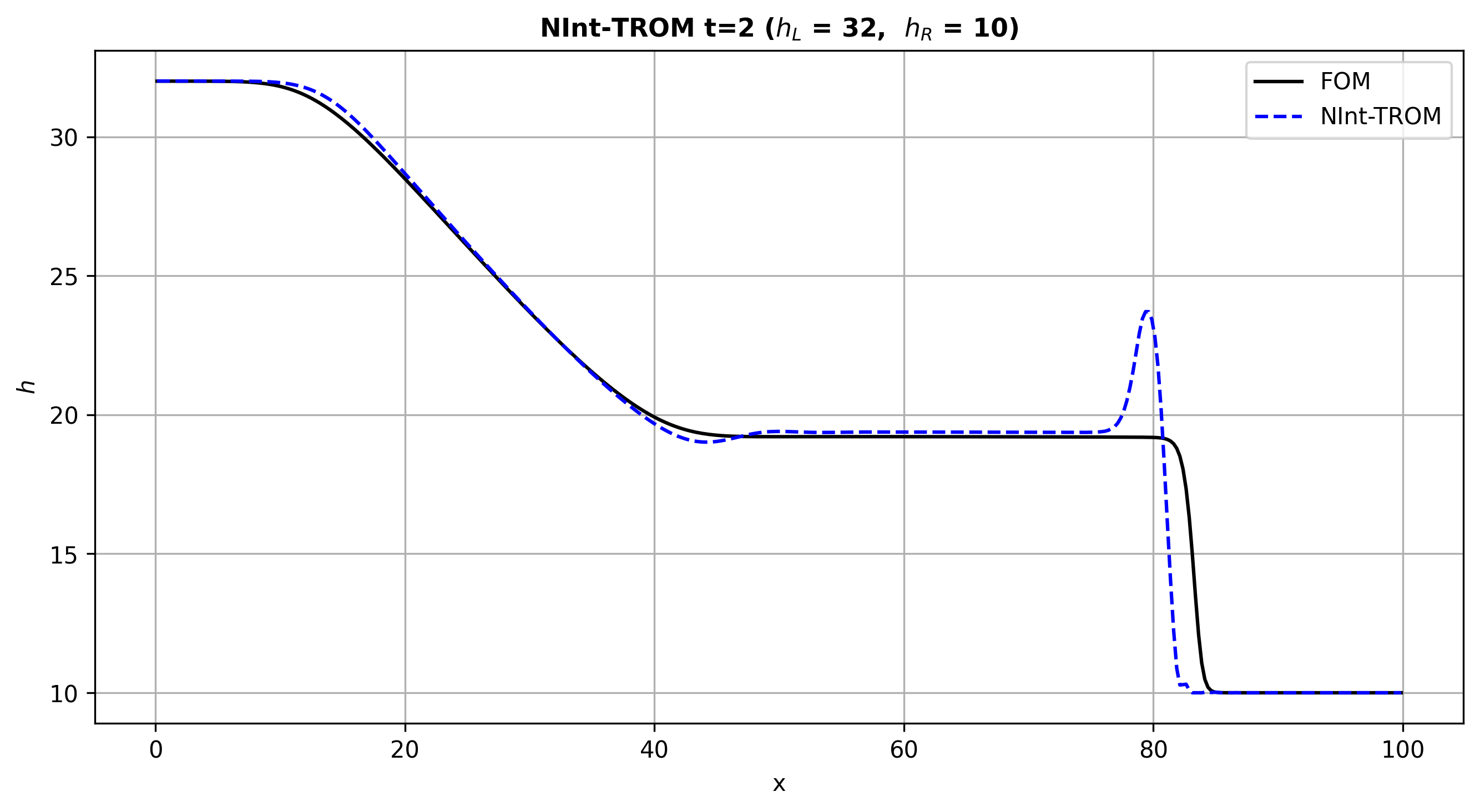}}
\caption{Comparison of $h(x,t;\muvec)$ from the FOM simulations and predictions of the PINN (left column) and the TROM (right column) for $\muvec=(20, 8.5), (29, 8.5), (32, 10)$ (top to bottom) at time $t=2$.}
\label{fig:extraph2}
\end{figure}

Figures \ref{fig:extraph1} and \ref{fig:extraph2} compare profiles for $h(x,t)$ and Table \ref{tab:extrapl2} presents relative $L^2$ errors at the final time $t=2$.
Profiles for the discharge $q(x,t)$ follow a similar trend and are not presented here for brevity.

\begin{table}[H]
\centering
\caption{Relative $L^2$ Errors in \eqref{eq:exp_rel_l2} for the TROM and the PINN parameter extrapolation predictions at time $t=2$.}
\label{tab:extrapl2}
\begin{tabular}{lcc|cc}
\toprule
$\muvec$ & $\varepsilon_h^{\text{TROM}}$ & $\varepsilon_h^{\text{PINN}}$ &
$\varepsilon_q^{\text{TROM}}$ &
$\varepsilon_q^{\text{PINN}}$ \\
\midrule
(9,4)      
& $1.28\times 10^{-2}$ & $6.87\times 10^{-3}$ 
& $6.73\times 10^{-2}$ & $2.3\times 10^{-2}$  \\
(20,8.5)      
& $7.25\times 10^{-3}$ & $2.88\times 10^{-3}$ 
& $3.36\times 10^{-2}$ & $0.91\times 10^{-2}$ \\
(29,4)      
& $1.97\times 10^{-2}$ & $2.58\times 10^{-3}$ 
& $6.35\times 10^{-2}$ & $0.35\times 10^{-2}$ \\
(32,4)      
& $6.79\times 10^{-2}$ & $5.64\times 10^{-3}$ 
& $2.06\times 10^{-1}$ & $0.59\times 10^{-2}$ \\
(29,8.5)   
& $1.37\times 10^{-2}$ & $3.64 \times 10^{-3}$ 
& $4.76\times 10^{-2}$ & $0.71\times 10^{-2}$ \\
(32,10)      
& $4.92\times 10^{-2}$ & $9.46\times 10^{-3}$ 
& $1.65\times 10^{-1}$ & $1.74\times 10^{-2}$ \\
\bottomrule
\end{tabular}
\end{table}

Overall, the performance of both models in the extrapolation regime is good; both models reproduce the large-scale solution structure and correctly track the propagating shock. 
However, the parametric PINN provides more accurate predictions in the vicinity of the shock and consistently yields lower relative errors than the TROM for all extrapolation cases considered here.
In particular, for $h_L>28$, the TROM prediction develops an artificial bump near the shock location, as seen in the middle and bottom rows of Figure~\ref{fig:extraph1} and the bottom row of Figure~\ref{fig:extraph2}. This localized artifact becomes more pronounced as the extrapolation distance in $h_L$ increases, whereas the PINN preserves a more accurate representation of the front.

These results indicate that both models are able to predict solutions for parameter values outside the training range, but the parametric PINN model provides more accurate extrapolations, particularly for larger values of $h_L$.

The spectral plots for the PINN model (not shown here for brevity) 
\rev{show good agreement with the FOM at low and intermediate wavenumbers, while the remaining mismatch is concentrated primarily in the high-wavenumber tail.}
Thus, even outside the training domain in parameter space, the learned PINN model continues to capture the dominant scales of the flow with good accuracy.


\section{Ablation Study}
\label{sec:ablation}

As discussed in Section~\ref{sec:pinn}, the design of the PINN model involves several nontrivial choices, some of which are specific to hyperbolic conservation laws. In this section, we perform a controlled, one-at-a-time ablation study around a fixed reference configuration to assess the effect of individual design choices. The goal is to identify which design decisions are crucial to the performance of the PINN model under the training and evaluation setup considered in this paper. As the reference configuration, we use the model described in Section~\ref{sec:pinn}, with two minor modifications introduced to reduce the training time: \rev{15,000 training epochs and a longer 50\% ramp-up period.
All ablation networks were initialized with random seed 42.}

Each ablation then turns off one parameter relative to this reference configuration, while the rest of the setup is kept unchanged.

Table~\ref{tab:ablation_summary} reports the mean space-time relative $L^2$ errors for $h$ and $q$ over a Monte Carlo benchmark consisting of $MC=99$ parameter instances $\mu=(h_L,h_R)$. To construct the benchmark, we first discretize the parameter domain $\mathcal{P}$ defined in~\eqref{eq:param_domain_box} using $20$ values of $h_L$ and $25$ values of $h_R$, resulting in a tensor-product set of $500$ candidate parameter instances. We then sample $MC$ instances uniformly without replacement from this set. The same sampled benchmark set is used to evaluate all ablation configurations. For a field $y\in\{h,q\}$, the space-time relative $L^2$ error corresponding to a parameter instance $\muvec$ is defined as
\[
E_{L^2}^{y}(\muvec)
=
\frac{
\left(
\int_0^T \int_{\Omega}
|y^{\mathrm{pred}}(x,t;\muvec)-y^{\mathrm{ref}}(x,t;\muvec)|^2
\,dx\,dt
\right)^{1/2}
}{
\left(
\int_0^T \int_{\Omega}
|y^{\mathrm{ref}}(x,t;\muvec)|^2
\,dx\,dt
\right)^{1/2}
}.
\]
The corresponding Monte Carlo means and the aggregate score are
\begin{equation}
\label{eq:S}
\overline{E}_{L^2}^{y}
=
\frac{1}{MC} \sum_{i=1}^{MC} E_{L^2}^{y}(\muvec_i), \;\; y\in \{h,q\}, \qquad 
S :=
\frac{1}{2} \left( \overline{E}_{L^2}^{h} + \overline{E}_{L^2}^{q} \right).
\end{equation}
The aggregate errors above are used as a ranking statistic for the ablation study. 

\begin{table}[h]
\centering
\small
\caption{Controlled ablation study around the reference PINN configuration. Each row modifies one component of the baseline while keeping the remaining setup fixed. Reported values are the aggregate errors defined in \eqref{eq:S}.}
\label{tab:ablation_summary}
\begin{tabular}{llccc}
\toprule
Variant & Change from baseline & $\overline{E}_{L^2}^{h}$ & $\overline{E}_{L^2}^{q}$ & $S$ \\
\midrule
Baseline & reference configuration & $1.358\times 10^{-2}$ & $3.629\times 10^{-2}$ & $2.494\times 10^{-2}$ \\
No PDE & $\lambda_{\mathrm{pde}}=0$ & $1.528\times 10^{-2}$ & $3.904\times 10^{-2}$ & $2.716\times 10^{-2}$ \\
Reduced $q$ weight & $\lambda_u=1$ & $1.306\times 10^{-2}$ & $4.138\times 10^{-2}$ & $2.722\times 10^{-2}$ \\
No ramp-up & $\lambda_{\mathrm{pde}}$, $\lambda_{\mathrm{bc}}$, and $\lambda_{\mathrm{ic}}$ are constant & $1.264\times 10^{-2}$ & $3.464\times 10^{-2}$ & $2.364\times 10^{-2}$ \\
$L^1$ PDE loss & $L^2 \rightarrow L^1$ in PDE term & $2.122\times 10^{-2}$ & $4.052\times 10^{-2}$ & $3.087\times 10^{-2}$ \\
Shock-aware off & uniform collocation & $4.031\times 10^{-2}$ & $1.165\times 10^{-1}$ & $7.839\times 10^{-2}$ \\
PDE gate off & no residual gating & $2.142\times 10^{-2}$ & $5.253\times 10^{-2}$ & $3.697\times 10^{-2}$ \\
No residual scaling & EMA scaling removed & $2.603\times 10^{-2}$ & $6.092\times 10^{-2}$ & $4.347\times 10^{-2}$ \\
\bottomrule
\end{tabular}
\end{table}
\FloatBarrier

The dominant trend in Table~\ref{tab:ablation_summary} is that the largest degradation occurs when shock-aware collocation (see Eq. \eqref{eq:shock_excl_band})
is removed. Replacing the shock-aware strategy with uniform collocation increases the aggregate score from $2.494\times 10^{-2}$ to $7.839\times 10^{-2}$, more than tripling the baseline value. The degradation is significant for both fields, but is especially pronounced for the discharge $q$, whose mean relative error rises from $3.629\times 10^{-2}$ to $1.165\times 10^{-1}$. It shows that careful treatment of the moving front region is essential for obtaining good accuracy, especially for the discharge $q$.

The next most important components are PDE gating 
(see Eq. \eqref{eq:pdegate})
and residual scaling (see Eq. \eqref{eq:pde_residual_scaling}). Disabling the soft PDE gate increases the aggregate score to $3.697\times 10^{-2}$, while removing EMA residual scaling yields $4.347\times 10^{-2}$. Both changes degrade $h$ and $q$, with a particularly clear effect on the discharge variable. Taken together, these results indicate that, once the collocation points are placed in a shock-aware manner, the way the strong-form residual is weighted and balanced remains critical. In other words, this ablation study shows that both (i) \emph{where} (in space) the PDE is enforced and (ii) \emph{how} it is enforced are important.

Removing the PDE loss altogether also worsens performance, increasing the aggregate score from $2.494\times 10^{-2}$ to $2.716\times 10^{-2}$. This confirms that physics-informed supervision provides measurable benefit even in the presence of substantial paired training data. However, the increase in error is much smaller than for the shock-aware, gating, and residual-scaling ablations. This suggests that the main challenge for this shock-dominated benchmark is not the inclusion of a PDE loss in itself, but rather enforcing that PDE loss
remains meaningful near shocks.

The discharge weighting parameter $\lambda_u$ (see eqs. \eqref{eq:data_loss} and \eqref{eq:bc_loss_method}) also has a visible effect. Reducing it from $5$ to $1$ leaves the mean depth error nearly unchanged, but increases the mean discharge error from $3.629\times 10^{-2}$ to $4.138\times 10^{-2}$. This is consistent with the observation made throughout the paper that predicting $q$ is more difficult than predicting $h$, and that the discharge component benefits from additional emphasis during optimization.

By contrast, removing ramp-up produces only a small change in the aggregate score. However, it appears to act more as a stabilizing or preventive training choice, which may still be useful in practice even if its direct effect on the aggregate error is limited in this sweep.

Overall, the ablation study shows that the most important components of the present training design are those related to shock-aware collocation, PDE gating, and residual scaling. By comparison, removing the PDE term or changing the discharge weight leads to a more moderate degradation, while removing ramp-up has little effect on the reported aggregate metric. For this benchmark, the main gains therefore come from approaches that improve how the residual is enforced near sharp fronts.

\section{Conclusion}
\label{sec:conc}
Fast and accurate parametric reduced-order models are important for a broad range of practical applications, including inverse problems, shape optimization, and digital twins. Such applications often require efficient evaluation of the underlying physical system across many parameter configurations. Developing reduced models for unified representation of parametric solutions is a challenging task because the behavior of solutions can change drastically for different parameter values, even for parameters in the training dataset.

In this paper, we compared and contrasted two parametric reduced models—the PINN and the TROM. Both models were constructed using the same dataset \rev{generated by a FOM solver} and evaluated at the same parameter values.
\rev{In this work, the accuracy of the surrogate models is therefore understood as their fidelity to the common FOM parameter-to-solution map generated by the FOM LLF solver, rather than to the exact solution of the SWE.
Consequently, both surrogate models might inherit some of the discretization characteristics of the underlying FOM data, including the numerical diffusion introduced near discontinuities.}
\rev{Neither model requires time integration during online evaluation, allowing solutions to be reconstructed rapidly; the PINN can reconstruct solutions at arbitrary times, whereas the TROM reconstructs solutions at the snapshot times stored in the training data.}. In our simulations, both models were more than two orders of magnitude (approximately \(500\) to \(700\) times) faster than the full-order model at \(t=2\). Moreover, the inference cost of either model is essentially independent of the queried time \(t\), since neither approach uses time stepping to predict the solution.
This substantial acceleration potentially allows the efficient treatment of realistic, computationally demanding problems.

Both the TROM and the PINN 
perform well in the interpolation cases considered here,
producing small relative errors overall, with the largest discrepancies concentrated near the moving shock. The TROM is more accurate in the dry-bed regime. Both models perform comparably in the challenging near-dry-bed regime at $h_R=0.14$, where they exhibit some loss of accuracy near the moving shock. In contrast, the PINN performs better in the wet-bed regimes with larger values of \(h_R\) considered here.
\rev{In addition, the PINN 
spectra show good agreement with the FOM over a broad range of wavenumbers in most cases, with the largest discrepancies occurring at high wavenumbers in the near-dry-bed regime \(h_R=0.14\).}
Overall, the near-dry-bed regime with $h_R=0.14$ remains the most challenging for both models because of the sharp, small-scale features associated with the propagating shock.

In addition, both models encounter difficulties when the shock exits the computational domain \([0,L]\), as occurs in the near-dry-bed regime \(\muvec=(26,0.14)\) at \(t=2.5\). Although \(t=2.5\) is included in the training dataset, neither model accurately reproduces the solution profile near the right boundary \(x=L\). Away from the right boundary, however, both the water depth and the discharge are predicted accurately.
One possible strategy for mitigating this limitation is to initialize a short-time FOM simulation using the reduced-model prediction at an earlier time, thereby resolving the shock as it exits the computational domain.

A more demanding test of a reduced model is its ability to generalize to parameter regimes outside the range represented in the training data.
Demonstrating predictive accuracy under such conditions extends the model’s applicability for parameter values outside of the training range and increases its value for practical use.
The PINN performs better in the extrapolation regimes considered here. This result is consistent with the construction of the TROM, which relies on interpolation (in parameter space) among nearby solutions in the training dataset to approximate the solution at a new parameter value. Because interpolation-based approximations can lose accuracy outside the sampled parameter domain, this limitation naturally carries over to the TROM when it is applied outside of the training range.
While the PINN reproduces the solutions with high accuracy in all extrapolation cases, the TROM exhibits larger discrepancies, concentrated primarily near the shock front.
If predictions outside the initial training range are required, the extrapolation performance of the TROM can be improved by augmenting the training set with full-order solutions at selected new parameter values and updating the low-rank representation.

To compare the offline stages of the two approaches, we first note that the TROM can efficiently handle large datasets generated on fine parameter grids. In contrast, PINN training is more computationally intensive and requires iterative batch processing of solution samples. Consequently, its training cost can grow substantially with the size of the dataset.
In addition, the offline stage of the TROM is considerably faster than PINN training. 
\rev{
As quantified in Section~\ref{sec:num}, the offline construction of the TROM is considerably less expensive than PINN training.
}

The offline phase of the TROM consists of computing a low-rank representation of the full snapshot tensor, and the singular-value threshold in \eqref{eq:rank_selection} is the only parameter controlling rank selection. In contrast, the PINN has several tunable hyperparameters and requires shock tracking and shock-aware collocation during training. 
The complexity of the PINN construction described in this paper reflects the challenges inherent in approximating hyperbolic problems with shocks.

The ablation study further demonstrates that shock-aware collocation, PDE gating, and residual scaling are essential components of the proposed PINN framework. Removing any of these components substantially increases the prediction errors, whereas removing the PDE loss itself produces a smaller degradation. These results indicate that, for shock-dominated problems, where and how the strong-form residual is enforced may be more important than simply including it in the loss function.

One of the major advantages of PINNs is their ability to work with incomplete or irregularly sampled data in parameter space. In contrast, TROMs requires a rectangular tensor-product grid in parameter space and therefore do not readily support local parameter grid refinement or irregular parameter domains. This restriction follows directly from the tensor representation of the full snapshot dataset. 
When data are missing, modern tensor-completion methods can be incorporated into a TROM framework; however, they increase computational complexity and introduce an additional source of approximation error.
The requirement that TROMs use rectangular tensor-product grids can lead to very large datasets, particularly for high-resolution or high-dimensional parameter discretizations. PINNs are more flexible in this respect because they can be trained using solution snapshots sampled nonuniformly or irregularly in parameter space.
Incorporating new data into a PINN is relatively straightforward because an existing model can be efficiently updated through transfer learning or fine-tuning. In contrast, TROMs rely on rectangular tensor-product grids, so a single additional FOM solution generally cannot be incorporated directly. Expanding the parameter grid may require additional FOM simulations to complete the new tensor.

Both TROMs and PINNs can be extended to two- and three-dimensional problems. However, extending TROMs is relatively straightforward; see, for example, \cite{Mizan2025TROM,Mizan2026TROM}. In contrast, extending PINNs is more challenging because of the higher computational cost of neural-network training and the possible need to track multidimensional shocks. These questions will be investigated in future work.

Both models favor a sufficiently regular dependence of the solution on the parameters. This dependence is explicit in the TROM, which interpolates between neighboring solutions in parameter space, and is also embedded in the PINN through its smooth neural-network representation. Near the dry-bed limit $h_R=0$, however, the solution depends rapidly and non-smoothly on $h_R$: the wet-bed regime $h_R>0$ contains a propagating shock, whereas the dry-bed regime $h_R=0$ is dominated by a rarefaction wave. This loss of parametric regularity potentially explains why the near-dry-bed regime $h_R=0.14$ is particularly challenging for both models.

\rev{Finally, the accuracy of both surrogate models in this study is assessed relative to the common FOM parameter-to-solution map generated by the Local Lax-Friedrichs solver. Consequently, both models may inherit some of the discretization characteristics of the reference data, including the numerical diffusion near discontinuities. This dependence is direct for the TROM, which is constructed entirely from the FOM snapshots. The PINN also uses these snapshots through its supervised loss, but additionally incorporates the governing SWE, initial conditions, and boundary conditions through the physics-informed loss. Thus, the PINN is not constrained solely by the numerical data. The same TROM and PINN formulations can, in principle, be applied to snapshot data generated by a higher-resolution shock-capturing solver, although the models would need to be reconstructed and, in the case of the PINN, retrained. 
Therefore, both surrogate models are evaluated against the numerical LLF solution of the FOM, which allows the surrogate approximation error to be assessed separately from the discretization error of the underlying full-order scheme.}

Overall, the two reduced-order modeling frameworks considered here should be viewed as complementary rather than competing approaches. This study highlights the strengths and potential limitations of each approach when applied to hyperbolic problems. TROMs offer a fast and relatively straightforward construction when solution data are available on a structured parameter grid, whereas PINNs provide greater flexibility in handling irregularly sampled data and extrapolating beyond the training regime. These advantages come with different costs: TROMs rely on compressed tensor-product representations, while PINNs require substantially more computational effort and careful training to capture shock-dominated dynamics. The choice between the two approaches should therefore depend on the available data, the structure of the parameter domain, and the intended application. Combining elements of both frameworks may provide a promising direction for developing more accurate, flexible, and efficient reduced-order models for complex hyperbolic systems.

\subsubsection*{Data Availability}
The numerical code and data supporting the parametric PINN results are publicly available at
\url{https://github.com/Saizt/parametric-swe-pinn}.

\appendix
\section{PINN Design}
\label{app:pinn}
\paragraph{Estimating shock position.}
For each regime \(\muvec=(h_L,h_R)\) and time \(t\), we estimate a dominant shock location \(x_{\mathrm{sh}}(t;\muvec)\) from the reference FOM snapshots. Let
\[
y(\cdot,t;\muvec)\in\{h(\cdot,t;\muvec),\,u(\cdot,t;\muvec)\}
\]
be the chosen indicator field on the uniform spatial grid \(\{x_i\}_{i=1}^{N_x}\). We first apply a moving-average filter,
\begin{equation*}
  \widetilde{y}_i
  :=
  \frac{1}{m}
  \sum_{k=-(m-1)/2}^{(m-1)/2} y_{i+k},
\end{equation*}
where \(m\ge 3\) is an odd window size. We then approximate the spatial derivative by finite differences:
\begin{equation*}
  (\partial_x \widetilde{y})_i
  \approx
  \begin{cases}
    \dfrac{\widetilde{y}_{i+1}-\widetilde{y}_i}{\Delta x}, & i=1,\\[1.2ex]
    \dfrac{\widetilde{y}_{i+1}-\widetilde{y}_{i-1}}{2\Delta x}, & 2\le i\le N_x-1,\\[1.2ex]
    \dfrac{\widetilde{y}_{i}-\widetilde{y}_{i-1}}{\Delta x}, & i=N_x,
  \end{cases}
\end{equation*}
where \(\Delta x\) is the spatial mesh size. The shock location is then defined as the grid point at which the gradient magnitude is maximal over an interior subset \(\Omega_{\mathrm{int}}\subset\Omega\), i.e.,
\begin{equation*}
  x_{\mathrm{sh}}(t;\muvec)
  \approx
  x_{i^\ast},
  \qquad
  i^\ast
  \in
  \arg\max_{x_i\in\Omega_{\mathrm{int}}}
  \left|(\partial_x \widetilde{y})_i\right|,
\end{equation*}
where \(\Omega_{\mathrm{int}}\) excludes a small boundary buffer near \(x=0\) and \(x=L\) in order to avoid boundary artifacts.

This produces shock-position estimates at the discrete snapshot times stored in the FOM dataset. Since collocation times are sampled continuously in \([0,T]\), the shock location at an arbitrary time \(t\) is obtained by interpolating between the neighboring discrete shock-position estimates.

\rev{
For shock-aware collocation, $\mu$ is sampled uniformly from the discrete set of 65 training parameter instances $P$, while $t$ is sampled continuously from $[0,T]$. Hence, the shock trajectory is interpolated only in time; no interpolation in parameter space is performed. Only FOM snapshots corresponding to $\mu\in P$ are used to construct these trajectories, and no shock-location information from the validation or test parameter values enters the training procedure. For the final model, the moving-average window is $m=7$, the outer $2\%$ of the spatial domain is excluded from the shock search, and shock exclusion is disabled for $h_R<\epsilon_{\rm dry}=0.5$.

The shock trajectories are extracted from the same FOM snapshots used to construct the supervised PINN training data and therefore require no additional simulations; shock detection and temporal interpolation constitute only an additional preprocessing step.}

\paragraph{PDE Residual Rescaling.}
Recall PDE residual scaling in equation \eqref{eq:pde_residual_scaling}.
At training iteration \(n\), we first compute the batch RMS values
\begin{align}
  r_{i,\mathrm{rms}}^{(n)}
  &:=
  \left(
    \frac{1}{N_f}\sum_{j=1}^{N_f}\mathcal{R}_i(z_j^{(f)};\theta^{(n)})^2
    + \varepsilon_{\mathrm{scale}}
  \right)^{1/2},
\end{align}
and then update the running scales according to
\begin{equation}
  s_i^{(n)}=\beta s_i^{(n-1)}+(1-\beta)\,r_{i,\mathrm{rms}}^{(n)},\qquad s_i^{(1)} := r_{i,\mathrm{rms}}^{(1)},
\end{equation}
where \(\beta\in(0,1)\) is the EMA parameter and \(\varepsilon_{\mathrm{scale}}>0\) is a small constant for numerical stability. 

\paragraph{PDE residual gating.}
For the \(j\)-th collocation point we define
\begin{equation}
  w_j^{\mathrm{grad}}
  =
  \exp\!\left(
    -\alpha\,\bigl|\partial_x h_\theta(x_j^{(f)},t_j^{(f)};\mu_j^{(f)})\bigr|
  \right),
\end{equation}
and
\begin{equation}
  w_j^{(h)}
  =
  \min\!\left\{
    1,\,
    \max\!\left\{
      0,\,
      \frac{h_\theta(x_j^{(f)},t_j^{(f)};\mu_j^{(f)})-h_{\min}}{h_{\min}}
    \right\}
  \right\},
\end{equation}
where \(\alpha>0\) controls sensitivity to steep gradients and \(h_{\min}>0\) defines a small-depth threshold. Then, the final gate is
\begin{equation}
\label{eq:pdegate}
  w_j
  =
    \max\!\left\{
      w_{\min},\,
      w_j^{\mathrm{grad}}\,w_j^{(h)}
    \right\},
\end{equation}
where \(w_{\min}\in(0,1)\) is a prescribed lower bound. 
Thus, PDE residual enforcement is weakened near sharp fronts and near-dry regions, while remaining nonzero everywhere else.
\rev{The gating weights $w_j$ are recomputed from the current network prediction at each loss evaluation, but are detached from the automatic-differentiation graph and therefore treated as fixed pointwise weights when differentiating the PDE loss.}


\end{document}